\documentclass[10pt]{article}
\usepackage[top=0.8in, bottom=1in, left=1in, right=1in]{geometry}

\usepackage[utf8]{inputenc}
\usepackage{amsthm,amsmath}
\numberwithin{equation}{section}
\usepackage{amssymb, mathrsfs, color}
\usepackage{enumerate}
\usepackage{verbatim}
\usepackage{tikz-cd}
\usepackage{enumitem}
\usepackage{indentfirst}
\usepackage{hyperref}
\hypersetup{
    colorlinks=true,
    linkcolor=blue,
    citecolor=brown,
    filecolor=magenta,
    urlcolor=blue,
}
\usepackage{titlesec}
\usepackage{fancyhdr}
\usepackage{mathtools}
\usepackage{CJKutf8}
\usepackage{enumitem}

\usetikzlibrary{decorations.pathmorphing}

\newtheorem{thm}{Theorem}[section]
\newtheorem{prob}[thm]{Problem}
\newtheorem{ques}[thm]{Question}
\newtheorem{cor}[thm]{Corollary}
\newtheorem{lem}[thm]{Lemma}
\newtheorem*{lem*}{Lemma}
\newtheorem{prop}[thm]{Proposition}

\newtheorem{innercustomgeneric}{\customgenericname}
\providecommand{\customgenericname}{}
\newcommand{\newcustomtheorem}[2]{%
  \newenvironment{#1}[1]
  {%
   \renewcommand\customgenericname{#2}%
   \renewcommand\theinnercustomgeneric{##1}%
   \innercustomgeneric
  }
  {\endinnercustomgeneric}
}

\newcustomtheorem{thm*}{Theorem}
\newcustomtheorem{prop*}{Proposition}
\newcustomtheorem{cor*}{Corollary}

\theoremstyle{definition}
\newcustomtheorem{defn*}{Definition}
\newtheorem{defn}[thm]{Definition}
\newtheorem{conv}[thm]{Convention}
\newtheorem{note}[thm]{Notation}

\theoremstyle{remark}
\newtheorem{remark}[thm]{Remark}

\newenvironment{remark*}[2][Remark]{\begin{trivlist}
\item[\hskip \labelsep {\bfseries #1}\hskip \labelsep {\bfseries #2.}]}{\end{trivlist}}

\newcommand{\bro}[1]{{\textcolor{brown}{#1}}}

\newcommand{\As}{\mathscr{A}}

\newcommand{\D}{\mathscr{D}}

\newcommand{\E}{\mathcal{E}}

\newcommand{\N}{\mathbb{N}}
\newcommand{\R}{\mathbb{R}}

\newcommand{\M}{\mathcal{M}}

\newcommand{\Sc}{\mathscr{S}}
\newcommand{\V}{\mathcal{V}}

\newcommand{\X}{\mathfrak{X}}
\newcommand{\Xs}{\mathscr{X}}
\newcommand{\Gr}{\mathbf{Gr}}

\newcommand{\Z}{\mathbb{Z}}
\newcommand{\fb}{\mathbf{f}}
\newcommand{\gb}{\mathbf{g}}
\newcommand{\Clog}{C^{\mathrm{LogL}}}
\newcommand{\clog}{c^{\mathrm{logL}}}

\newcommand{\mulog}{{\mu_{\log}}}

\newcommand{\Su}{\mathsf{S}}
\newcommand{\De}{\mathsf{P}}
\newcommand{\id}{\operatorname{id}}

\newcommand{\Lie}[1]{\mathrm{Lie}_{#1}}

\newcommand{\sgn}{\operatorname{sgn}}

\newcommand{\eps}{\varepsilon}

\newcommand{\Span}{\operatorname{Span}}
\newcommand{\Vol}{\operatorname{Vol}}

\newcommand{\loc}{\mathrm{loc}}
\newcommand{\supp}{\operatorname{supp}}
\newcommand{\dist}{\operatorname{dist}}
\newcommand{\divg}{\operatorname{div}}
\newcommand{\rank}{\operatorname{rank}}

\newcommand{\Coorvec}[1]{\frac\partial{\partial#1}}

\begin{document}
\author{Liding Yao}
\title{The Frobenius Theorem for Log-Lipschitz Subbundles}
\date{}
\maketitle
\begin{abstract}
We extend the definition of involutivity to non-Lipschitz tangent subbundles using generalized functions. We prove the Frobenius Theorem with sharp regularity estimate when the subbundle is log-Lipschitz: if $\mathcal V$ is a log-Lipschitz involutive subbundle of rank $r$, then for any $\varepsilon>0$, locally there is a homeomorphism $\Phi(u,v)$ such that $\Phi,\frac{\partial\Phi}{\partial u^1},\dots,\frac{\partial\Phi}{\partial u^r}\in C^{0,1-\varepsilon}$, and $\mathcal V$ is spanned by the continuous vector fields $\Phi_*\frac\partial{\partial u^1},\dots,\Phi_*\frac\partial{\partial u^r}$.
\end{abstract}

\tableofcontents

\section{Introduction}

Given a smooth manifold $M$ and a smooth tangent subbundle $\V\le TM$, we say $\V$ is \textbf{involutive}, if for any vector fields $X,Y$ that are sections of $\V$, their Lie bracket $[X,Y]$ is also a section of $\V$.



The celebrated Frobenius theorem states that a smooth involutive tangent subbundle is always integrable, in the sense that locally there is a  foliation of a smooth family of submanifolds, such that the original bundle equals to the union of tangent bundles of these submanifolds. In our formulation, it is the following:
\begin{prop}[Frobenius theorem \cite{Frobenius}]\label{Prop::Intro::SmoothFro}
Let $M$ be an $n$-dimensional smooth manifold and let $\V$ be a smooth subbundle of $TM$ with rank $r$. Let $(u,v)=(u^1,\dots,u^r,v^1,\dots,v^{n-r})$ be the standard coordinate system for $\R^r\times\R^{n-r}$.

Suppose $\V$ is involutive. Then for any $p\in M$ there is a neighborhood $\Omega\subseteq\R^r_u\times\R^{n-r}_v$ of $0$, and a smooth regular parameterization\footnote{For $k\ge1$, we say that $\Phi$ is a $C^k$-\textbf{regular parameterization}, if $\Phi$ is defined on an open subset of $\R^n$ and $\Phi$ is a $C^k$-diffeomorphism onto its image. It is equivalent to say that if $\Phi$ is injective and $\Phi^{-1}$ is a $C^k$-coordinate chart.} $\Phi(u,v):\Omega\to M$ such that $\Phi(0)=p$, and for any $(u,v)\in\Omega$, the tangent vectors $\frac{\partial\Phi}{\partial u^1}(u,v),\dots,\frac{\partial\Phi}{\partial u^r}(u,v)\in T_{\Phi(u,v)}M$ form a basis of $\V_{\Phi(u,v)}$.
\end{prop}
See \cite{ShortFro} for a simple proof of Proposition \ref{Prop::Intro::SmoothFro}. It is well known that a similar theorem also holds for the $C^1$-subbundles with $C^1$-regular parameterization.
In this paper we extend the Proposition \ref{Prop::Intro::SmoothFro} to the log-Lipschitz setting.

\subsection{The main result}
Let $U\subseteq\R^n$ be an open subset, we say a continuous function $f:U\to\R^m$ is bounded \textbf{log-Lipschitz}, denoted by $f\in\Clog(U;\R^m)$, if
\begin{equation}\label{Eqn::Intro::ClogNorm}
    \|f\|_{\Clog(U;\R^m)}:=\sup_{x\in U}|f(x)|+\sup_{x,y\in U;0<|x-y|<1}|f(x)-f(y)|\Big(|x-y|\log\frac e{|x-y|}\Big)^{-1}<\infty.
\end{equation}

Note that every log-Lipschitz functions are Lipschitz, but not the converse.

Let $M$ be a smooth manifold, we say  a vector field $X$ on $M$ is  log-Lipschitz, if it has (locally) log-Lipschitz coefficients in every local coordinate systems.
\begin{defn}
We say a tangent subbundle $\V$ over $M$ is log-Lipschitz, if for any $p\in M$ there is a neighborhood $U\subseteq M$ of $p$ and log-Lipschitz vector fields $X_1,\dots,X_r$ on $U$ where $r=\rank \V$, such that the linear space $\V_q$ is spanned by $X_1(q),\dots,X_r(q)\in T_qM$ for every $q\in U$.
\end{defn}

For log-Lipschitz subbundles we introduce the notion of distributional\footnote{Throughout our paper, we use \textbf{distributions} as \textbf{generalized functions}, which are defined to be linear functionals on some test function spaces. For the corresponding differential geometry terminology, we always use \textbf{tangent subbundles}.} involutivity as follows. 

\begin{defn}\label{Def::Intro::DisInv}
Let $M$ be a $C^{1,1}$ manifold and let $\V\le TM$ be a log-Lipschitz tangent subbundle. We say $\V$ is \textbf{involutive in the sense of distributions}, if for any log-Lipschitz vector fields $X=\sum_{i=1}^nX^i\Coorvec{x^i}$, $Y=\sum_{j=1}^nY^j\Coorvec{x^j}$ on $\V$ and any log-Lipschitz 1-form $\theta=\sum_{k=1}^n\theta_kdx^k$ on $\V^\bot$, the generalized function
\begin{equation}\label{Eqn::Intro::DefInvEqn}
    \langle\theta,[X,Y]\rangle:=\sum_{i,j=1}^n\Big(\theta_iX^j\frac{\partial Y^i}{\partial x^j}-\theta_jY^i\frac{\partial X^j}{\partial x^i}\Big),
\end{equation}
is identically zero in the sense of distributions.
\end{defn}

Here $\V^\bot:=\{(p,l)\in T^*M:p\in M,\ l(v)=0,\forall v\in\V_p\}\le T^*M$ is the dual bundle of $\V$. 
Note that the right hand side of \eqref{Eqn::Intro::DefInvEqn} is not necessarily a $L^1_\loc$-function. It makes sense as a distribution, see Lemma \ref{Lem::Hold::MultLoc} and Corollary \ref{Cor::Hold::[X,Y]WellDef}. Also see Definition \ref{Def::DisInv::DefFunVF}. 

We see that distributional involutivity is the natural generalization to the classical involutivity. See \ref{Section::DisInv::CharInv} for a detailed discussion.

\medskip
Our main result is the following.
\begin{thm}\label{Thm::MainThm1}
Let $M$ be $n$-dimensional $C^{1,1}$ manifold, and let $\V\le TM$ be a rank-$r$ log-Lipschitz tangent subbundle, which is involutive in the sense of distributions.

Then for any $p\in M$ and any $\eps>0$, there is a neighborhood $\Omega\subseteq\R^r_u\times\R^{n-r}_v$ of $0$ and a map $\Phi(u,v):\Omega\to M$ (depending on $\eps$) satisfying:
\begin{enumerate}[parsep=-0.3ex,label=(\roman*)]
    \item\label{Item::MainThm1::Phi0} $\Phi(0)=p$.
    \item\label{Item::MainThm1::PhiReg} $\Phi\in C^{0,1-\eps}(\Omega;M)$ and $\frac{\partial\Phi}{\partial u^1},\dots,\frac{\partial\Phi}{\partial u^r}\in C^{0,1-\eps}(\Omega;TM)$.
    \item\label{Item::MainThm1::Phi-1} $\Phi$ is homeomorphic to its image, and $\Phi^{-1}\in C^{0,1-\eps}(\Phi(\Omega);\R^n)$.
    \item\label{Item::MainThm1::Span} For any $(u,v)\in \Omega$, $\frac{\partial\Phi}{\partial u^1}(u,v),\dots,\frac{\partial\Phi}{\partial u^r}(u,v)\in T_{\Phi(u,v)}M$ form a basis of the linear subspace $\V_{\Phi(u,v)}$.
\end{enumerate}
\end{thm}
\begin{remark}
Theorem \ref{Thm::MainThm1} \ref{Item::MainThm1::PhiReg} shows that both $\Phi$ and $\Phi^{-1}$ are $C^{0,1-\eps}$. In general $\Phi\in C^{0,1-\eps}$ does not imply $\Phi^{-1}\in C^{0,1-\eps}$ and vice versa, since we do not have the inverse function theorem for $C^{0,1-\eps}$-maps. 
\end{remark}

In Theorem \ref{Thm::MainThm1} the $C^{0,1-\eps}$-regularity is sharp, in the sense that there exists a log-Lipschitz subbundle such that either $\Phi\notin C^{0,1^-}$ or $\Phi^{-1}\notin C^{0,1^-}$, see Proposition \ref{Prop::Further::SharpProp}. Here $C^{0,1^-}=\bigcap_{0<\delta<1}C^{0,1-\delta}$.
In particular, we cannot pick a $\Phi$ such that $\Phi$ and $\Phi^{-1}$ are both log-Lipschitz.

\medskip
We are allowed to pick a  $\Phi$ such that both $\Phi$ and $\Phi^{-1}$ are $C^{0,1^-}$ when $\V$ is slightly regular:


Let $U\subseteq\R^n$ be an open subset, we say a continuous function $f:U\to\R^m$ is bounded \textbf{little log-Lipschitz}, denoted by $f\in\clog(U;\R^m)$, if $f$ is log-Lipschitz and 
$$\lim\limits_{r\to0}\sup\limits_{x,y\in U;0<|x- y|<r}|f(x)-f(y)|\Big(|x-y|\log\frac e{|x-y|}\Big)^{-1}=0.$$ 
Clearly $\clog(U;\R^m)\subsetneq\Clog(U;\R^m)$ is a closed subspace.

We define little log-Lipschitz subbundle in the identical way. Namely, a subbundle that is locally spanned by little log-Lipschitz vector fields.
\begin{defn}
We say a tangent subbundle $\V$ over $M$ is little log-Lipschitz, if for any $p\in M$ there is a neighborhood $U\subseteq M$ of $p$ and little log-Lipschitz vector fields $X_1,\dots,X_r$ on $U$ where $r=\rank \V$, such that the linear space $\V_q$ is spanned by $X_1(q),\dots,X_r(q)\in T_qM$ for every $q\in U$.
\end{defn}


\begin{thm}\label{Thm::MainThm2}
Let $M$ be $n$-dimensional $C^{1,1}$ manifold, and let $\V\le TM$ be a rank-$r$ little log-Lipschitz tangent subbundle, which is involutive in the sense of distributions.

Then for any $p\in M$, there is a neighborhood $\Omega\subseteq\R^r_u\times\R^{n-r}_v$ of $0$ and a map $\Phi(u,v):\Omega\to M$ satisfying:
\begin{enumerate}[parsep=-0.3ex,label=(\roman*)]
    \item\label{Item::MainThm2::Phi0} $\Phi(0)=p$.
    \item\label{Item::MainThm2::PhiReg} $\Phi\in  C^{0,1^-}(\Omega;M)$ and $\frac{\partial\Phi}{\partial u^1},\dots,\frac{\partial\Phi}{\partial u^r}\in  C^{0,1^-}(\Omega;TM)$.
    \item\label{Item::MainThm2::Phi-1} $\Phi$ is homeomorphic to its image, and $\Phi^{-1}\in C^{0,1^-}(\Phi(\Omega);\R^n)$.
    \item\label{Item::MainThm2::Span} For any $(u,v)\in \Omega$, $\frac{\partial\Phi}{\partial u^1}(u,v),\dots,\frac{\partial\Phi}{\partial u^r}(u,v)\in T_{\Phi(u,v)}M$ form a basis of the linear subspace $\V_{\Phi(u,v)}$.
\end{enumerate}
\end{thm}

Even in the setting of Theorem \ref{Thm::MainThm2}, we cannot make $\Phi$ and $\Phi^{-1}$ to be both log-Lipschitz in general, unless we take impose some stronger regularity assumptions on $\V$. See Remark \ref{Rmk::Further::PreciseCountGenThm}.

\subsection{Historical remarks for Frobenius type theorem in low regularities}

For vector fields $X$ and $Y$ which are not $C^1$, the commutator $[X,Y]$ may not be defined pointwise. Because of this, any Frobenius type theorem for subbundles of less than $C^1$-regularity requires a suitable notion of involutivity. This leads to the following:

\begin{ques}\label{Ques::Intro::Sec}
Let $\V$ be a tangent subbundle which is not $C^1$, and let vector fields $X,Y$ be sections of $\V$.
\begin{itemize}[parsep=-0.3ex]
    \item What is the good definition of the Lie bracket $[X,Y]$?
    \item With this definition how can we say whether $[X,Y]$ is a section of $\V$? 
\end{itemize}
\end{ques}

When $X$ and $Y$ are Lipschitz vector fields, the Lie bracket $[X,Y]$ makes sense as a vector field with $L^\infty$-coefficients. And we can talk about sections and involutivity almost everywhere. In \cite{Simic}, Simi\'c gives the following definition:
\begin{defn}[Almost everywhere characterization]\label{Def::Intro::AEsection}Let $M$ be a $C^{1,1}$-manifold, let  $\V$ be a Lipschitz subbundle of $TM$, and let $X,Y$ be two Lipschitz sections of $\V$.
We say $[X,Y]$ is an \textbf{almost everywhere section} of $\V$, if as a vector field with $L^\infty$-coefficients, $[X,Y](p)\in\V_p$ holds for almost every $p\in M$.

We say $\V$ is \textbf{almost everywhere involutive}, if for any Lipschitz sections $X,Y$ of $\V$, the Lie bracket $[X,Y]$ is an almost everywhere section of $\V$.
\end{defn}

 Later Rampazzo \cite{Rampazzo} gave another characterizations using set-valued Lie brackets:
\begin{defn}[Set valued characterization]\label{Def::Intro::SetInv}
Let $M$ be a $C^{1,1}$-manifold, let  $\V$ be a Lipschitz subbundle of $TM$, and let $X,Y$ be two Lipschitz sections of $\V$. The \textbf{set-valued Lie bracket} $[X,Y]_{set}\subseteq\V$ is defined as the convex hall $[X,Y]_{set}(p):=\operatorname{conv}([X,Y]_{set}^+(p))\subseteq T_pM$, $p\in M$, where
$$[X,Y]_{set}^+:=\overline{\{(q,[X,Y](q))\in TM:X,Y\text{ are differentiable at }q\}}\subseteq TM.$$
And we say $\V$ is \textbf{set-valuedly involutive}, if $[X,Y]_{set}\subseteq\V$ for every Lipschitz sections $X,Y$ of $\V$.
\end{defn}

Definitions \ref{Def::Intro::AEsection} and \ref{Def::Intro::SetInv} are equivalent, see \cite[Section 2]{Rampazzo} and \cite[Theorem 4.11 (I) and (III)]{Rampazzo}.

An advantage of set valued bracket is that it is defined pointwise everywhere.
Note that a Lipschitz function is differentiable almost everywhere, so $[X,Y]_{set}(p)\neq\varnothing$ for all $p\in M$.

\begin{prop}[Lipschitz Frobenius theorem, \cite{Simic,Rampazzo}]
Let $M$ be an $n$-dimensional $C^{1,1}$ manifold, and let $\V\le TM$ be a Lipschitz tangent subbundle of rank $r$. If $\V$ is either almost everywhere involutive or set-valuedly involutive, then for any $p\in M$, there is a neighborhood $\Omega\subseteq\R^r_u\times\R^{n-r}_v$ of $0$ and a bi-Lipschitz parameterization $\Phi(u,v):\Omega\to M$ satisfying the following:
\begin{enumerate}[nolistsep,label=(\roman*)]
    \item $\Phi(0)=p$.
    \item \label{LipFroConse2} $\frac{\partial\Phi}{\partial u^1},\dots,\frac{\partial\Phi}{\partial u^r}:\Omega\to TM$ are Lipschitz continuous maps.
    \item For any $(u,v)\in \Omega$, $\frac{\partial\Phi}{\partial u^1}(u,v),\dots,\frac{\partial\Phi}{\partial u^r}(u,v)\in T_{\Phi(u,v)}M$ form a basis of the linear subspace $\V_{\Phi(u,v)}$.
\end{enumerate}
\end{prop}

Both \cite{Simic} and \cite{Rampazzo} only deal with Lipschitz vector fields and Lipschitz subbundles. What happen if they are non-Lipschitz? Indeed, $[X,Y]$ is not necessarily $L^1_\loc$ so we cannot use an almost everywhere characterization. A non-Lipschitz section can be nowhere differentiable so we cannot use set-valued bracket either. 

To resolve this problem, we use generalized functions. When the vector fields $X,Y$ are log-Lipschitz, or more generally when $X,Y$ are H\"older $C^{0,\alpha}$ for some $\frac12<\alpha<1$, using techniques of paradifferential calculus we can define $[X,Y]$ as a vector field whose coefficients are distributions (see Corollary \ref{Cor::Hold::[X,Y]WellDef}). 
In fact the coefficients of $[X,Y]$ are in $C^{-1,\alpha}$: the space of distributions that can be written as the sum of derivatives of $C^{0,\alpha}$-functions, see Definition \ref{Def::Hold::HoldDef} \ref{Item::Hold::HoldDef::-1}.

By imitating Definition \ref{Def::Intro::DisInv} we can give a definition for the distributional sections of $\V$ as follows:
\begin{defn}[Distributional sections of subbundles]\label{Def::Intro::DisSec}
Let $0<\alpha<1$ and $1-\alpha<\beta<1$. Let $\V$ be a $C^{0,\alpha}$ tangent subbundle of a manifold $M$. We say a $C^{-1,\beta}$-vector field $Z=\sum_{i=1}^na^i\Coorvec{x^i}$ is a distributional section of $\V$, if for any log-Lipschitz section $\theta$ of $\V^\bot$, written as $\theta=\sum_{j=1}^n\theta_jdx^j$, the generalized function\footnote{Using Lemma \ref{Lem::Hold::MultLoc} \ref{Item::Hold::MultLoc::WellDef}, the right hand side of \eqref{Eqn::Intro::DisSec::DefEqn} makes sense as a distribution in $C^{-1,\beta}(M)$. Also see Definition \ref{Def::DisInv::DefFunVF}.}
\begin{equation}\label{Eqn::Intro::DisSec::DefEqn}
    \langle\theta,Z\rangle=\sum_{j=1}^na^j\theta_j\in\D'(M),
\end{equation}
is identically zero in the sense of distributions.

We denote by $C^{-1,\beta}_\loc(M;\V)$ the space of $C^{-1,\beta}$-sections of $\V$.
\end{defn}
When $\frac12<\beta<1$ we can define distributional involutivity of $\V$ by replacing every log-Lipschitz with H\"older $C^{0,\beta}$ in Definition \ref{Def::Intro::DisInv} (see Definition \ref{Def::DisInv::DefInv}). 

The characterization of distributional sections and distributional involutivity recover the above characterizations when $\V$ is Lipschitz. Indeed, for Lipschitz vector fields $X,Y$ that are sections of $\V$, $[X,Y]=0$ lays in $\V$ almost everywhere if and only if $\langle \theta,[X,Y]\rangle=0$, so Definitions \ref{Def::Intro::DisInv} and \ref{Def::Intro::AEsection} are equivalent. 

\medskip
Recently in \cite{ContFro}, a version of involutivity was introduced using approximations by smooth functions. \cite{ContFro} works on differential forms (see \cite[Definition 1.12]{ContFro}). An equivalent characterization using vector fields is as follows:
\begin{defn}\label{Def::Intro::AsyInv}
Let $\V$ be a continuous rank-$r$ tangent subbundle over $\R^n$.
We say $\V$ is \textbf{strongly asymptotically involutive}, if for any $p\in \R^n$ there is a neighborhood $U\subseteq \R^n$ of $p$, and a sequence of $C^1$-vector fields $\{X_1^\nu,\dots,X^\nu_r\}_{\nu=1}^\infty$ on $U$, such that
\begin{enumerate}[parsep=-0.3ex,label=(\roman*)]
    \item For each $j=1,\dots,r$, $\{X_j^\nu\}_{\nu=1}^\infty$ uniformly converges to a continuous vector field $X_j$, such that for each $q\in U$, $X_1|_q,\dots,X_r|_q\in T_qM$ form a linear basis of $\V_q$.
    \item There exist $\{c_{jk}^{l\nu}\}_{\nu=1}^\infty\subset C^0(\bar U)$ for $1\le j,k,l\le r$, such that 
    \begin{equation}\label{Eqn::Intro::AsyInv1}\exists t_0>0,\
        \lim\limits_{\nu\to\infty}\bigg(\max\limits_{1\le j,k\le r}\Big\|[X_j^\nu,X_k^\nu]-\sum_{l=1}^rc_{jk}^{l\nu}X_l^\nu\Big\|_{C^0(U;\R^n)}\bigg)\cdot\exp\big(t_0\cdot\max\limits_{1\le l\le r}\|\nabla X_l^\nu\|_{C^0(U;\R^{n\times n})}\big)=0.
    \end{equation}
\end{enumerate}
\end{defn}

In \cite{ContFro} it is shown that if $\V$ is strongly asymptotically involutive, then locally there exists a topological parameterization\footnote{By \textbf{topological parameterization} we mean a continuous map $\Phi:\Omega\subseteq\R^n\to M$ which is homeomorphic to its image.} $\Phi(u,v)\in C^1_u(C^0_v)$ such that $\Phi_*\Coorvec {u^1},\dots,\Phi_*\Coorvec{u^r}$ span $\V$. However, given a continuous subbundle, we do not know an algorithm to check whether such an approximation exists. In contrast, Definition \ref{Def::DisInv::DefInv} can be directly checked using coordinate systems. Nevertheless \eqref{Eqn::Intro::AsyInv1} is central to our proof of Theorems \ref{Thm::MainThm1} and \ref{Thm::MainThm2}. 

In the log-Lipschitz setting we prove that distributional involutivity implies strongly asymptotically involutivity. See Corollary \ref{Cor::Hold::CorAsyInv} and Remark \ref{Rmk::Hold::CorAsyInv}.

\subsection{Overview of the proof}\label{OverProof}

By Osgood uniqueness theorem \cite{Osgood,Markus},  a log-Lipschitz vector field $X$ always has ODE uniqueness, so we write $e^{tX}(p)$ for its \textbf{flow map}. That is, for any $p\in M$, the map $[t\mapsto e^{tX}(p)]:I\subseteq\R\to M$ is the maximal  $C^1$ solution $\gamma$ to the equation\footnote{In fact $e^{tX}$ can always be defined as long as $X$ is integrable and \ref{Eqn::Intro::ODEflowDef} has unique solution locally for all $p$ in the domain.}
\begin{equation}\label{Eqn::Intro::ODEflowDef}
    \dot\gamma(t)=X(\gamma(t))\in T_{\gamma(t)}M,\quad \gamma(0)=p.
\end{equation}
We denote $\exp_X(t,p):=e^{tX}(p)$.

For a given log-Lipschitz involutive subbundle $\V\le TM$, by linear algebra arguments we can find a log-Lipschitz local basis $(X_1,\dots,X_r)$ of $\V$ such that $[X_j,X_k]=0$, $1\le j,k\le r$ in the domain, see Lemma \ref{Lem::PfThm::CanGen}. To construct a (topological) parameterization $\Phi$, we can use the compositions of flows $e^{t^1X_1}\circ\dots\circ e^{t^rX_r}$.

The main difficulty to the proof of Theorems \ref{Thm::MainThm1} and \ref{Thm::MainThm2} is the flow commuting problem:
\begin{prob}[Flow commuting problem]\label{Prob::Intro::FlowComProb}
Let $X,Y$ be two vector fields such that their ODE flows $e^{tX}$, $e^{sY}$ are defined locally. 

Suppose the Lie bracket $[X,Y]$ is well-defined  and equals to 0, then do we have $e^{tX}\circ e^{sY}=e^{sY}\circ e^{tX}$ for small $t,s\in\R$?
\end{prob}

When $X$ and $Y$ are Lipschitz, it is known in \cite[Section 5]{RampazzoSussmanCommutators} for example.
In this paper we give positive answer for Problem \ref{Prob::Intro::FlowComProb} when $X,Y$ are log-Lipschitz with $[X,Y]=0$ holds in the sense of distributions (Proposition \ref{Prop::ODE::FlowComm}). The key is to show that all log-Lipschitz involutive subbundles are strongly asymptotically involutive as in Definition \ref{Def::Intro::AsyInv}, see Corollary \ref{Cor::Hold::CorAsyInv}. A simplified version of this is the following:  there are a smaller neighborhood $U$ and vector fields $\{X_\nu, Y_\nu\}_{\nu=1}^\infty\subset C^\infty(U;\R^n)$, such that $X_\nu\xrightarrow{ C^0(\bar U)} X$, $Y_\nu\xrightarrow{C^0(\bar U)} Y$, and 
 \begin{equation}\label{Eqn::Intro::AsyInv2}
     \exists t_0>0,\ \lim\limits_{\nu\to\infty}\|[X_\nu,Y_\nu]\|_{C^0(U;\R^n)}e^{t_0(\|\nabla X_\nu\|_{C^0(U;\R^{n\times n})}+\|\nabla Y_\nu\|_{C^0(U;\R^{n\times n})})}=0.
\end{equation}

To see such approximations $\{X_\nu, Y_\nu\}_{\nu=1}^\infty$ exist, we use the following: 

Let $(\Su_\nu)_{\nu=0}^\infty$ be a sequence of Littlewood-Paley summation operators (see Definition \ref{Def::Hold::LPchar}). We use $C^{-1,\beta}$ for the space of distributions that can be written as the formal sum of the derivatives of $C^{0,\beta}$-functions (see Definition \ref{Def::Hold::HoldDef} \ref{Item::Hold::HoldDef::-1}). Then
\begin{thm}[Approximating the vanished product]\label{Thm::Intro::ApproxThm}
Let $m,n\in\Z_+$ and let $0<\alpha,\beta<1$ satisfy $\alpha+\beta>0$. There is a $C_{n,m,\alpha,\beta}>0$ such that, if $\fb\in C^{0,\alpha}(\R^n;\R^m)$ and $\gb\in C^{-1,\beta}(\R^n;\R^m)$ satisfy $\fb\cdot \gb=0$ in the sense of distributions, then
\begin{equation}\label{Eqn::Intro::ApproxThmEqn}
    \|(\Su_\nu\fb)\cdot(\Su_\nu\gb)\|_{C^0(\R^n)}\le C_{n,m,\alpha,\beta}\cdot2^{-(\alpha+\beta-1)\nu}\|\fb\|_{C^{0,\alpha}(\R^n;\R^m)}\|\gb\|_{C^{-1,\beta}(\R^n;\R^m)},\qquad \forall\nu\in\Z_{\ge0}.
\end{equation}
\end{thm}

In application we take $m=2n$ and $\fb^j=(X^1,\dots,X^n,-Y^1,\dots,-Y^n)$, $\gb^j=(\frac{\partial Y^1}{\partial x^j},\dots,\frac{\partial Y^1}{\partial x^j},\frac{\partial X^1}{\partial x^j},\dots,\frac{\partial X^n}{\partial x^j})$, for each $j=1,\dots,n$, so that $\fb^j\cdot \gb^j=[X,Y]^j=0$.


\subsection*{Acknowledgement}
I would like to express my appreciation to my advisor Brian Street for his suggestions and encouragement, and  Andrew Lewis for informative discussions on modules and sheaves. I also want to thank Ruofan Jiang and Jiaqi Hou for suggestions on terminologies in algebra.

\section{Paraproduct Decompositions for Functions and Distributions}
In this section, using Bony's decompostion from Littlewood-Paley theory, 
we explain why $[X,Y]$ is defined when $X$ and $Y$ are log-Lipschitz vector fields and prove Theorem \ref{Thm::Intro::ApproxThm}.
\subsection{Preliminary tools}
First we introduce some tools for computing the products of H\"older functions and distributions.

Let $\{\phi_j\}_{j=0}^\infty\subset \Sc(\R^n;\R)$ be a sequence of real-valued Schwartz functions, whose Fourier transforms\footnote{We use $\hat f(\xi)=(2\pi)^{-\frac n2}\int_{\R^n}f(x)e^{- ix\xi}dx$ as the definition of Fourier transform. Note that $f$ is real-valued if and only if $\hat f$ is an even function.} satisfy the following:
\begin{itemize}[parsep=-0.3ex]
    \item $\hat\phi_0\in C_c^\infty(\R^n)$ is an even function and has support in $\{\xi\in\R^n:|\xi|<2\}$.
    \item $\hat\phi_j(\xi)=\hat\phi_0(2^{-j}\xi)-\hat\phi_0(2^{1-j}\xi)$ for $j\ge1$. Equivalently, $\phi_j(x)=2^{jn}\phi_0(2^jx)-2^{(j-1)n}\phi_0(2^{j-1}x)$ for $j\ge1$.
\end{itemize}

We call such $\{\phi_j\}_{j=0}^\infty$ a (Fourier) \textbf{dyadic resolution of unity}.

In this construction we have $\sum_{j=0}^\infty\hat\phi_j(\xi)=1$ for all $\xi\in\R^n$, and
\begin{equation}\label{Eqn::Hold::SuppofLPChar}
    \supp\hat\phi_0\subseteq\{|\xi|<2\};\quad\supp\hat\phi_j\subseteq\{2^{j-1}<|\xi|<2^{j+1}\}\text{ for }j\ge1.
\end{equation}

See \cite[Section 2.3.1]{TriebelTheoryOfFunctionSpacesI} for details. 

\begin{defn}\label{Def::Hold::LPchar}
Let $\{\phi_j\}_{j=0}^\infty$ be as above, we define the \textbf{Littlewood-Paley block operators} $\{\Su_j,\De_j\}_{j=0}^\infty$ that act on the space of tempered distributions $\Sc'(\R^n)$ as the following:
\begin{equation}\label{Eqn::Hold::DefBlockOP}
    \De_jf:=\phi_j\ast f,\qquad \Su_jf:=\sum_{k=0}^j\De_kf=2^{jn}\phi_0(2^j\cdot)\ast f,\quad j\ge0,\ f\in\Sc'(\R^n).
\end{equation}
\end{defn}

We require $\{\phi_j\}$ to be real so that for a real vector field $X$, the smoothing one $\Su_jX$ is still a real vector field.

From \eqref{Eqn::Hold::SuppofLPChar}, we have
\begin{equation}\label{Eqn::Hold::LPCharComp}
    \De_j\De_k=0,\quad\text{if }|j-k|\ge2;\qquad\De_j\Su_k=\begin{cases}\Su_k,&\text{if }j\le k-2,\\0,&\text{if }j\ge k+2.\end{cases}
\end{equation}


We use the following definitions for H\"older spaces:

\begin{defn}\label{Def::Hold::HoldDef}
Let $m\ge-1$ be an integer, let $0<\alpha\le1$, and let $\Omega\subseteq\R^n$ be an open set. The \textbf{bounded H\"older space} $C^{m,\alpha}(\Omega)$ is given as follows:
\begin{enumerate}[parsep=-0.3ex,label=(\roman*)]
    \item For $m\ge0$, we define $C^{m,\alpha}(\Omega):=\big\{f\in C^m(\bar\Omega):\sup\limits_{x,y\in \Omega,x\neq y}\frac{|\nabla^mf(x)-\nabla^mf(y)|}{|x-y|^\alpha}<\infty\big\}$ with norm $$\|f\|_{C^{m,\alpha}(\Omega)}:=\sum_{|\kappa|\le m}\sup\limits_{x\in\Omega}|\partial^\kappa f(x)|+\sum_{|\kappa|=m}\sup\limits_{x,y\in \Omega,x\neq y}|x-y|^{-\alpha}|\partial^\kappa f(x)-\partial^\kappa f(y)|.$$
    \item\label{Item::Hold::HoldDef::-1} For $m=-1$, we define $C^{-1,\alpha}(\Omega):=\{g+\divg G:G\in C^{0,\alpha}(\Omega;\R^n),g\in C^{0,\alpha}(\Omega)\}\subset\D'(\Omega)$ with norm
    $$\|f\|_{C^{-1,\alpha}(\Omega)}:=\inf\Big\{\|g\|_{C^{0,\alpha}(\Omega)}+\sum_{j=1}^n\|G_j\|_{C^{0,\alpha}(\Omega)}:f=g+\sum_{j=1}^n\frac{\partial G_j}{\partial x^j}\text{ in the sense of distributions}\Big\}.$$
\end{enumerate}
We use the following notations of \textbf{locally bounded H\"older spaces}: 
\begin{enumerate}[parsep=-0.3ex,label=(\roman*)]\setcounter{enumi}{2}
    \item\label{Item::Hold::HoldDef::Local} $C^{m,\alpha}_\loc(\Omega):=\{f\in\D'(\Omega'):f\in C^{m,\alpha}(\Omega'),\ \text{for every precompact open } \Omega'\Subset\Omega\}$, provided that $m\ge-1$.
    \item\label{Item::Hold::HoldDef::CAlp-} $C^{m,\alpha^-}_\loc(\Omega):=\bigcap_{0<\beta<\alpha}C^{m,\beta}_\loc(\Omega)$, provided that $m\ge-1$.
    \item\label{Item::Hold::HoldDef::CAlp+} $C^{m,\alpha^+}_\loc(\Omega):=\bigcup_{\alpha<\beta<1}C^{m,\beta}_\loc(\Omega)$, provided that $m\ge-1$ and $\alpha<1$.
\end{enumerate}
\end{defn}

We have an important characterization of H\"older spaces:

\begin{prop}\label{Prop::Hold::LPchar}Let $m\ge-1$ and $0<\alpha<1$. Then $f\in C^{m,\alpha}(\R^n)$ if and only if $f$ is a tempered distribution satisfying $\sup\limits_{j\in\Z_{\ge0}}2^{j(m+\alpha)}\|\De_jf\|_{L^\infty(\R^n)}<\infty$. Moreover there is $C_{m,\alpha}>0$ such that 
\begin{equation}\label{Eqn::Hold::LPChar}
    C_{m,\alpha}^{-1}\|f\|_{C^{m,\alpha}(\R^n)}\le\sup\limits_{j\in\Z_{\ge0}}2^{j(m+\alpha)}\|\De_jf\|_{L^\infty(\R^n)}\le C_{m,\alpha}\|f\|_{C^{m,\alpha}(\R^n)},\qquad\forall f\in\Sc'(\R^n).
\end{equation}
\end{prop}
\begin{proof}

Following the notations in \cite[Definition 2.3.1/2(i)]{TriebelTheoryOfFunctionSpacesI} we have the space 
\begin{equation*}
    B_{\infty\infty}^{m+\alpha}(\R^n):=\{f\in\Sc'(\R^n):\|f\|_{B_{\infty\infty}^{m+\alpha}} <\infty\},\  \|f\|_{B_{\infty\infty}^{m+\alpha}}:=\sup_{j\in\Z_{\ge0}}2^{j(m+\alpha)}\|\De_jf\|_{L^\infty(\R^n)},\quad m\ge-1,\ \alpha\in(0,1).
\end{equation*}

So \eqref{Eqn::Hold::LPChar} is equivalent as saying $\|f\|_{C^{m,\alpha}(\R^n)}\approx_{m,\alpha,\De}\|f\|_{B_{\infty\infty}^{m+\alpha}}$. 

The case $m\ge0$ can be found in \cite[Sections 2.2.2 and 2.5.7]{TriebelTheoryOfFunctionSpacesI} or \cite[Theorem 1.4.9]{GrafakosModern}. 
We explain the proof of $m=-1$ and $0<\alpha<1$ as the following.

Let $f\in C^{-1,\alpha}(\R^n)$, write $f=g_0+\sum_{j=1}^n\partial_jg_j$ where $g_0,\dots,g_n\in C^{0,\alpha}(\R^n)=B_{\infty\infty}^\alpha(\R^n)$. By \cite[Theorem 2.3.8(i)]{TriebelTheoryOfFunctionSpacesI} we have $\|\partial_jg_j\|_{B_{\infty\infty}^{\alpha-1}}\lesssim\|g_j\|_{B_{\infty\infty}^\alpha}$. Clearly $B^\alpha_{\infty\infty}(\R^n)\subset B^{\alpha-1}_{\infty\infty}(\R^n)$. Therefore
\begin{equation*}\textstyle
    \|f\|_{B_{\infty\infty}^{\alpha-1}}\le\|g_0\|_{B_{\infty\infty}^{\alpha-1}}+\sum_{j=1}^n\|\partial_jg_j\|_{B_{\infty\infty}^{\alpha-1}}\lesssim\sum_{j=0}^n\|g_j\|_{B_{\infty\infty}^{\alpha}}\approx\sum_{j=0}^n\|g_j\|_{C^{0,\alpha}(\R^n)}.
\end{equation*}
Taking the infimum of the decompositions $f=g_0+\sum_{j=1}^n\partial_jg_j$ we get $\|f\|_{B_{\infty\infty}^{\alpha-1}}\lesssim\|f\|_{C^{-1,\alpha}(\R^n)}$.

Conversely for $f\in B_{\infty\infty}^{\alpha-1}(\R^n)$ we consider $f=(I-\Delta)(I-\Delta)^{-1}f=(I-\Delta)^{-1}f+\sum_{j=1}^n\partial_j(-\partial_j(I-\Delta)^{-1}f)$. By \cite[Theorem 2.3.8(i)]{TriebelTheoryOfFunctionSpacesI} again (note that $I_s=(I-\Delta)^\frac s2$ in the reference) we have $\|(I-\Delta)^{-1}f\|_{B_{\infty\infty}^{\alpha+1}}+\sum_{j=1}^n\|\partial_j(I-\Delta)^{-1}f\|_{B_{\infty\infty}^\alpha}\lesssim\|f\|_{B_{\infty\infty}^{\alpha-1}}$. Since $B^{\alpha+1}_{\infty\infty}(\R^n)\subset B^{\alpha}_{\infty\infty}(\R^n)$,  we have 
\begin{equation*}
    \|(I-\Delta)^{-1}f\|_{C^{0,\alpha}(\R^n)}+\sum_{j=1}^n\|\partial_j(I-\Delta)^{-1}f\|_{C^{0,\alpha}(\R^n)}\approx\|(I-\Delta)^{-1}f\|_{B_{\infty\infty}^{\alpha}}+\sum_{j=1}^n\|\partial_j(I-\Delta)^{-1}f\|_{B_{\infty\infty}^\alpha}\lesssim\|f\|_{B_{\infty\infty}^{\alpha-1}}.
\end{equation*}
Take $g_0=(I-\Delta)^{-1}f$ and $g_j=-\partial_j(I-\Delta)^{-1}f$ we get $\|f\|_{C^{-1,\alpha}(\R^n)}\lesssim \|f\|_{B_{\infty\infty}^{\alpha-1}}$, completing the proof.
\end{proof}

Let $f\in C^{0,\alpha}(\R^n)$ and $g\in C^{-1,\beta}(\R^n)$. If $\alpha+\beta>1$, we are allowed to talk about the product $f\cdot g$, which makes sense as a distribution.
\begin{prop}\label{Prop::Hold::Product}
Let $0<\alpha,\beta<1$ satisfy $\alpha+\beta>1$. Then for $f\in C^{0,\alpha}(\R^n)$ and $g\in C^{-1,\beta}(\R^n)$, the product $fg$ makes sense as a distribution in $C^{-1,\beta}(\R^n)$. Moreover, there is a $C=C_{\alpha,\beta}>0$ such that
\begin{equation*}
    \|fg\|_{C^{-1,\beta}(\R^n)}\le C\|f\|_{C^{0,\alpha}(\R^n)}\|g\|_{C^{-1,\beta}(\R^n)}.
\end{equation*}
\end{prop}
This is the special case to \cite[Theorem 2.8.2(i)]{TriebelTheoryOfFunctionSpacesI} where we have $\rho=\alpha$ and $s=\beta-1$ in the reference.

In \cite{TriebelTheoryOfFunctionSpacesI} the operator $g\mapsto fg$ is called a ``pointwise muliplier'' associated with $f$. When $g\in C^0(\R^n)$ the operation is clearly pointwise, in the sense that $fg(x)$ depend only on $f(x)$ and $g(x)$. For general $g\in C^{-1,\beta}(\R^n)$, $fg$ is only a distribution which is not defined pointwisely. Nevertheless for an arbitrary open set $\Omega\subseteq\R^n$, the distribution $(fg)|_\Omega$ depends only on $f|_\Omega$ and $g|_\Omega$.

\begin{lem}\label{Lem::Hold::MultLoc}
Let $\Omega\subseteq\R^n$ be an open set, and let $0<\alpha,\beta<1$ satisfy $\alpha+\beta>1$. 
\begin{enumerate}[parsep=-0.3ex,label=(\roman*)]
    \item\label{Item::Hold::MultLoc::Loc} Suppose $\tilde f^1,\tilde f^2\in C^{0,\alpha}(\R^n)$ and $\tilde g^1,\tilde g^2\in C^{-1,\beta}(\R^n)$ satisfy $\tilde f^1|_\Omega=\tilde f^2|_\Omega$ and $\tilde g^1|_\Omega=\tilde g^2|_\Omega$. Then $(\tilde f^1\tilde g^1)|_\Omega=(\tilde f^2\tilde g^2)|_\Omega$ as distributions on $\Omega$.
    \item\label{Item::Hold::MultLoc::WellDef} In particular for $f\in C^{0,\alpha}_\loc(\Omega)$ and $g\in C^{-1,\beta}_\loc(\Omega)$, the product $fg$ is well-defined distribution in $C^{-1,\beta}_\loc(\Omega)$.
\end{enumerate}
\end{lem}
Also see \cite[Section 4]{BookProduct} for more details.
\begin{proof}

\ref{Item::Hold::MultLoc::Loc}: Let $i_1,i_2\in\{1,2\}$. Since $(\Su_\nu)_{\nu=0}^\infty$ is a Schwartz approximate identity, we have $\Su_\nu \tilde f^{i_1}\xrightarrow{C^{0,\alpha-\delta}}\tilde f^{i_1}$ and $\Su_\nu \tilde g^{i_2}\xrightarrow{C^{-1,\beta-\delta}}\tilde g^{i_2}$ for all $\delta>0$. Thus by Proposition \ref{Prop::Hold::Product} we have convergence $\lim_{j\to\infty}\Su_j\tilde f^{i_1}\Su_j\tilde g^{i_2}=\tilde f^{i_1}\tilde g^{i_2}$ in $C^{-1,\beta-\delta}(\R^n)$ for every $\delta>0$. By \cite[Lemma 4.2.2]{BookProduct} we see that $$(\tilde f^1\tilde g^1-\tilde f^2\tilde g^2)|_\Omega=((\tilde f^1-\tilde f^2)\tilde g^1+\tilde f^2(\tilde g^1-\tilde g^2))|_\Omega=\lim_{j\to\infty}\big(\Su_j(\tilde f^1-\tilde f^2)\cdot\Su_j\tilde g^1+\Su_j\tilde f^2\cdot\Su_j(\tilde g^1-\tilde g^2)\big)|_\Omega=0,\quad \text{in }\D'(\Omega).$$




\noindent\ref{Item::Hold::MultLoc::WellDef}: Let $f\in C^{0,\alpha}_\loc(\Omega)$ and $g\in C^{-1,\beta}_\loc(\Omega)$. By \ref{Item::Hold::MultLoc::Loc} it suffices to prove the well-definedness of $fg|_{\Omega'}\in C^{-1,\beta}(\Omega')$ on any precompact open $\Omega'\Subset\Omega$.

Let $\chi\in C_c^\infty(\Omega)$ be such that $\chi|_{\Omega'}\equiv1$. By taking zero extension outside $\Omega$ we have $\chi f\in C^{0,\alpha}(\R^n)$ and $\chi g\in C^{-1,\beta}(\R^n)$. By Proposition \ref{Prop::Hold::Product} we have $\chi f\cdot \chi g\in C^{-1,\beta}(\R^n)$. Since $\chi f|_{\Omega'}=f|_{\Omega'}$ and $\chi g|_{\Omega'}=g|_{\Omega'}$, by \ref{Item::Hold::MultLoc::Loc} we see that $fg|_{\Omega'}=\chi f\cdot \chi g|_{\Omega'}\in C^{-1,\beta}(\Omega')$ is defined, finishing the proof.
\end{proof}

In this way, we see that the Lie bracket $[X,Y]$ is defined as a distribution.
\begin{cor}\label{Cor::Hold::[X,Y]WellDef}
Let $\frac12<\alpha<1$, and let $\Omega\subseteq\R^n$ be an open set. If $X$ and $Y$ are $C^{0,\alpha}_\loc$ vector fields on $\Omega$, then $[X,Y]\in C^{-1,\alpha}_\loc(\Omega;\R^n)$ is a vector field with distribution coefficients.

In particular if $X$ and $Y$ are both log-Lipschitz, then $[X,Y]\in C^{-1,1^-}_\loc(\Omega;\R^n)$ (see Definition \ref{Def::Hold::HoldDef} \ref{Item::Hold::HoldDef::CAlp-}).
\end{cor}

Recall that for $X=\sum_{i=1}^nX^i\Coorvec{x^i}$ and $Y=\sum_{j=1}^nY^j\Coorvec{x^j}$ the Lie bracket $[X,Y]$ is given by
\begin{equation}\label{Eqn::Intro::LieBracketEqn}
    [X,Y]=\sum_{i,j=1}^n\Big(X^j\frac{\partial Y^i}{\partial x^j}-Y^j\frac{\partial X^i}{\partial x^j}\Big)\Coorvec{x^i}.
\end{equation}
\begin{proof} Write $X=:\sum_{j=1}^nX^j\Coorvec{x^j}$ and $Y=:\sum_{j=1}^nY^j\Coorvec{y^j}$, we have $X^j,Y^j\in C^{0,\alpha}_\loc(\Omega)$ and $\frac{\partial Y^i}{\partial x^j},\frac{\partial X^i}{\partial x^j}\in C^{-1,\alpha}_\loc(\Omega)$ for $1\le i,j\le n$. Applying Lemma \ref{Lem::Hold::MultLoc} \ref{Item::Hold::MultLoc::WellDef} since $\alpha+\alpha>1$, we get $X^j\frac{\partial Y^i}{\partial x^j}-Y^j\frac{\partial X^i}{\partial x^j}\in C^{-1,\alpha}_\loc(\Omega)$ for all $1\le i,j\le n$. By \eqref{Eqn::Intro::LieBracketEqn}, we conclude that $[X,Y]\in C^{-1,\alpha}_\loc(\Omega;\R^n)$.

Since $\Clog_\loc(\Omega)\subset C^{0,\alpha}_\loc(\Omega)$ for every $0<\alpha<1$, if $X$ and $Y$ are log-Lipschitz,  we have $[X,Y]\in C^{-1,\alpha}_\loc(\Omega;\R^n)$ for all $\frac12<\alpha<1$. Therefore $[X,Y]\in C^{-1,1^-}_\loc(\Omega;\R^n)$.
\end{proof}

Recall the Littlewood-Paley decomposition $f=\sum_{j=0}^\infty\De_jf$ and  $g=\sum_{j=0}^\infty\De_jg$. We can write $fg$ as
\begin{align}\notag
    fg&=\sum_{j,k=0}^\infty\De_jf\cdot\De_kg=\bigg(\sum_{j=0}^\infty\sum_{k=0}^{j-6}+\sum_{k=0}^\infty\sum_{j=0}^{j-6}+\sum_{j,k\ge0;|j-k|\le 6}\bigg)\De_jf\cdot\De_kg
    \\
    \label{Eqn::Hold::Bony}&= \sum_{j=0}^\infty\De_jf\cdot\Su_{j-6}g+\sum_{k=0}^\infty\De_kg\cdot\Su_{k-6}f+\sum_{j,k\ge0;|j-k|\le5}\De_jf\cdot\De_kg.
\end{align}
Here we use $\Su_l=0$ for $l\le-1$. The \eqref{Eqn::Hold::Bony} is called a \textbf{Bony decomposition} (also see \cite[Sections 2.6.1 and 2.8.1]{BahouriCheminDanchin}). It produces the following formula, which is  used in the proof of Proposition \ref{Prop::Hold::Product}:
\begin{lem}[Paraproduct formula]\label{Lem::Hold::ParaForm}Let $f$ and $g$ be tempered distributions such that the product $fg\in\Sc'(\R^n)$ is defined. Then we have the following decomposition:
\begin{equation}\label{Eqn::Hold::ParaDecomp}
\De_l(fg)=\De_l\sum_{j=l-2}^{l+2}(\De_jf)(\Su_{j-6}g)+\De_l\sum_{k=l-2}^{l+2}(\Su_{k-6}f)(\De_kg)+\De_l\sum_{\substack{|j-k|\le5\\j,k\ge l-2}}(\De_jf)(\De_kg),\quad l\ge0.
\end{equation}
Here we use $\De_j=\Su_j=0$ for $j\le-1$.
\end{lem}
One can find the proof in \cite[Section 2.8.1]{BahouriCheminDanchin} or \cite[Section 4.4.1]{GrafakosModern}. Note that we have an alternative decomposition $\De_l(fg)=\De_l\sum_{j=l-2}^{l+2}(\De_jf)(\Su_{j-3}g)+\De_l\sum_{k=l-2}^{l+2}(\Su_{k-3}f)(\De_kg)+\De_l\sum_{\substack{|j-k|\le2\\j,k\ge l-4}}(\De_jf)(\De_kg)$. It is the proof of Theorem \ref{Thm::Intro::ApproxThm} that requires us to use $\Su_{j-6}$ and $\Su_{k-6}$ rather than $\Su_{j-3}$ and $\Su_{k-3}$.
\begin{proof}[Proof of Lemma \ref{Lem::Hold::ParaForm}]
By \eqref{Eqn::Hold::DefBlockOP} we have $\widehat{\De_jf}=\hat\phi_j\hat f$, $\widehat{\De_kg}=\hat\phi_k\hat g$. Note that $(\De_jf\cdot\De_kg)^\wedge=\widehat{\De_jf}\ast\widehat{\De_kg}$, so $\supp(\De_jf\cdot\De_kg)^\wedge\subseteq\supp\widehat{\De_jf}+\supp\widehat{\De_kg}\subseteq\supp\hat\phi_j+\supp\hat\phi_k$. Using the support condition \eqref{Eqn::Hold::SuppofLPChar} we get
$$\supp\hat\phi_j+\supp\hat\phi_k\subseteq\big\{\max(2^{j-1}-2^{k+1},2^{k-1}-2^{j+1}) <|\xi|<2^{j+1}+2^{k+1}\},\quad\text{ for all }j,k\ge0.$$
Note that $\max(2^{j-1}-2^{k+1},2^{k-1}-2^{j+1})\ge 2^{\max(j,k)-2}>0$ when $|j-k|\ge3$. We see that when $|j-l|\ge 3$, $\supp(\De_jf\cdot\Su_{j-6}g)^\wedge\subseteq\{2^{j-2}<|\xi|<2^{j+2}\}$ has empty intersection with $\supp\hat\phi_l\subseteq\begin{cases}\{2^{l-1}<|\xi|<2^{l+1}\}&l\ge1,\\\{|\xi|<2\}&l=0.\end{cases}$
Therefore $\De_l\sum_{j=0}^\infty\sum_{k=0}^{j-6}\De_jf\cdot\Su_{j-6}g=\De_l\sum_{j=l-2}^{l+2}\De_jf\cdot\Su_{j-6}g$ for all $l\ge0$. By symmetry we also have $\De_l\sum_{k=0}^\infty\sum_{j=0}^{j-6}\De_kg\cdot\Su_{k-6}f=\De_l\sum_{k=l-2}^{l+2}\Su_{k-6}f\cdot\De_kg$.

For $|j-k|\le 6$ and $l\ge\max(j,k)+3$, we see that $\supp(\De_jf\cdot\De_kg)^\wedge\subseteq\{|\xi|<2^{\max{(j,k)}+2}\}$ is disjoint with $\supp\hat\phi_l\subseteq\{|\xi|>2^{l-1}\}$. Therefore $\De_l\sum_{j,k\ge0;|j-k|\le5}\De_jf\cdot\De_kg=\De_l\sum_{j,k\ge l-2;|j-k|\le5}\De_jf\cdot\De_kg$, finishing the proof.
\end{proof}

\subsection{Proof of Theorem \ref{Thm::Intro::ApproxThm} and corollaries on Lie brackets of vector fields}
Recall the Littlewood-Paley operators $\De_\nu=\phi_\nu\ast(-)$ and $\Su_\nu=\sum_{j=0}^\nu\phi_\nu\ast(-)=2^{n\nu}\phi_0(2^\nu\cdot)\ast(-)$ in Definition \ref{Def::Hold::HoldDef}, where $\{\phi_j\}_{j=0}^\infty $ is a dyadic resolution. Clearly $\|2^{n\nu}\phi_0(2^\nu\cdot)\|_{L^1}=\|\phi_0\|_{L^1}$ for all $\nu\ge0$, and $\|\phi_\nu\|_{L^1}=\|2^{n(\nu-1)}\phi_1(2^{\nu-1}\cdot)\|_{L^1}=\|\phi_1\|_{L^1}$ for all $\nu\ge1$. Therefore we have the following uniform bounds,
\begin{equation}\label{Eqn::Hold::PfApThm::PSC0Norm}
    \|\Su_\nu\|_{C^0\to C^0}\le\|2^{n\nu}\phi_0(2^\nu\cdot)\|_{L^1}=\|\phi_0\|_{L^1},\quad\|\De_\nu\|_{C^0\to C^0}\le\|\phi_\nu\|_{L^1}\le\max(\|\phi_0\|_{L^1},\|\phi_1\|_{L^1}),\quad\forall\nu\ge0.
\end{equation}

Fix $\alpha+\beta>1$,  $\fb\in C^{0,\alpha}(\R^n;\R^m)$ and $\gb\in C^{-1,\beta}(\R^n;\R^m)$ that satisfy $\fb\cdot\gb=0$ in the sense of distributions. By Proposition \ref{Prop::Hold::Product}, we know the equality $\fb\cdot\gb=0$ holds in $ C^{-1,\beta}(\R^n)$. Note that we already have the convergence $\lim\limits_{\nu\to\infty}(\Su_\nu\fb)\cdot(\Su_\nu\gb)=0$ in $ C^{-1,\beta-\delta}(\R^n)$, for every $\delta>0$. 
Theorem \ref{Thm::Intro::ApproxThm} says that this sequence  converges uniformly with quantitative speed, which is better than the convergence in some spaces of distributions.

\begin{proof}[Proof of Theorem \ref{Thm::Intro::ApproxThm}]We set $\|\fb\|_{ C^{0,\alpha}}=\|\gb\|_{ C^{-1,\beta}}=1$ without loss of generality.

By \eqref{Eqn::Hold::SuppofLPChar} $\Su_\nu f$ and $\Su_\nu g$ both have Fourier supports in $\{|\xi|<2^{\nu+1}\}$. Thus when $l\ge\nu+3$, $\supp(\Su_\nu f\cdot\Su_\nu g)^\wedge\subseteq\{|\xi|<2^{\nu+2}\}$  is disjoint from $\supp\hat\phi_l\subseteq\{|\xi|>2^{l-1}\}$, which means $\De_l(\Su_\nu f\cdot\Su_\nu g)=0$. Therefore we have
\begin{equation}\label{Eqn::Hold::PfApThm::Decomp}
    \Su_\nu \fb\cdot\Su_\nu \gb=\sum_{l=0}^\infty\De_l(\Su_\nu \fb\cdot\Su_\nu \gb)=\sum_{l=0}^{\nu+2}\De_l(\Su_\nu \fb\cdot\Su_\nu \gb)=\Su_{\nu-4}(\Su_\nu \fb\cdot\Su_\nu \gb)+\sum_{l=\nu-3}^{\nu+2}\De_l(\Su_\nu \fb\cdot\Su_\nu \gb).
\end{equation}

Our goal is to prove $\|\Su_{\nu-4}(\Su_\nu \fb\cdot\Su_\nu \gb)\|_{C^0(\R^n)}\lesssim2^{-(\alpha+\beta-1)\nu}$ and $\sum_{l=\nu-3}^{\nu+2}\|\De_l(\Su_\nu \fb\cdot\Su_\nu \gb)\|_{C^0(\R^n)}\lesssim2^{-(\alpha+\beta-1)\nu}$, where the implicit constants depend only on $m,n,\alpha,\beta$ and $\phi$.

\medskip
By \eqref{Eqn::Hold::ParaDecomp}, we know for any $l\in\Z_{\ge0}$,
\begin{equation}\label{Eqn::Hold::PfApThm::Para0}
\De_l(\fb\cdot\gb)=\De_l\sum_{j=l-2}^{l+2}(\De_j\fb)\cdot(\Su_{j-6}\gb)+\De_l\sum_{k=l-2}^{l+2}(\Su_{k-6}\fb)\cdot(\De_k\gb)+\De_l\sum_{j,k\ge l-2;|j-k|\le5}(\De_j\fb)\cdot(\De_k\gb)=0.
\end{equation}
So for $\nu,l\ge0$, we have decomposition
\begin{align}
    \label{Eqn::Hold::PfApThm::Para-1}
    \De_l(\Su_\nu\fb\cdot\Su_\nu\gb)=\De_l(\Su_\nu\fb\cdot\Su_\nu\gb)-\De_l(\fb\cdot\gb)=&\De_l\sum_{j=l-2}^{l+2}\big((\De_j\Su_\nu\fb)\cdot(\Su_{j-6}\Su_\nu\gb)-\De_j\fb\cdot\Su_{j-6}\gb\big)
    \\
    \label{Eqn::Hold::PfApThm::Para-2}
    &+\De_l\sum_{k=l-2}^{l+2}\big((\Su_{k-6}\Su_\nu\fb)\cdot(\De_k\Su_\nu\gb)-\Su_{k-6}\fb\cdot\De_k\gb\big)
    \\
    \label{Eqn::Hold::PfApThm::Para-3}
    &+\De_l\sum_{j,k\ge l-2;|j-k|\le5}\big((\De_j\Su_\nu\fb)\cdot(\De_k\Su_\nu\gb)-\De_j\fb\cdot\De_k\gb\big).
\end{align}

When $l\le \nu-4$ and $|j-l|\le 2$, we have $j\le\nu-2$. By \eqref{Eqn::Hold::LPCharComp} we have $\De_j\Su_\nu=\De_j$ and $\Su_{j-6}\Su_\nu=\Su_{j-6}$, so in this case \eqref{Eqn::Hold::PfApThm::Para-1} is zero. Similarly \eqref{Eqn::Hold::PfApThm::Para-2} vanishes as well. Thus for $\Su_{\nu-4}(\Su_\nu \fb\cdot\Su_\nu \gb)$ in \eqref{Eqn::Hold::PfApThm::Decomp} we have
\begin{align*}
    \Su_{\nu-4}(\Su_\nu \fb\cdot\Su_\nu \gb)
    =&\sum_{l=0}^{\nu-4}\De_l\sum_{j,k\ge l-2;|j-k|\le5}\big((\De_j\Su_\nu\fb\cdot \De_k\Su_\nu\gb)-\De_j\fb\cdot\De_k\gb\big)
    \\
    =&\sum_{l=0}^{\nu-4}\De_l\sum_{j,k\ge \nu-2;|j-k|\le5}\big((\De_j\Su_\nu\fb)\cdot(\De_k\Su_\nu\gb)-\De_j\fb\cdot\De_k\gb\big)\\
    =&\Su_{\nu-4}\sum_{j,k\ge \nu-2;|j-k|\le5}\big((\De_j\Su_\nu\fb)\cdot (\De_k\Su_\nu\gb)-\De_j\fb\cdot\De_k\gb\big).
\end{align*}

Using \eqref{Eqn::Hold::PfApThm::PSC0Norm} and Proposition \ref{Prop::Hold::LPchar} we have an estimate for \eqref{Eqn::Hold::PfApThm::Para-3}: for every $l\ge0$,
\begin{equation}\label{Eqn::Hold::PfApThm::Term3C0}
    \begin{aligned}
    &\sum_{\substack{j,k\ge l-2\\|j-k|\le5}}\big\|(\De_j\Su_\nu\fb)\cdot (\De_k\Su_\nu\gb)-\De_j\fb\cdot\De_k\gb\big\|_{C^0}\le \sum_{\substack{j,k\ge l-2\\|j-k|\le5}}\big(\|\Su_\nu\|_{C^0\to C^0}^2+1\big)\|\De_j\fb\|_{C^0}\|\De_k\gb\|_{C^0}
    \\
    \lesssim&_\phi\sum_{\substack{j,k\ge l-2\\|j-k|\le5}}\|\De_j\fb\|_{C^0}\|\De_k\gb\|_{C^0}\lesssim_{\phi,\alpha,\beta}\sum_{\substack{j,k\ge l-2\\|j-k|\le5}}2^{-j\alpha}2^{(1-\beta)k}\lesssim\sum_{q\ge l}2^{-q\alpha}2^{q(1-\beta)}\lesssim2^{-l(\alpha+\beta-1)}.
    \end{aligned}
\end{equation}

Taking $l=\nu$ in \eqref{Eqn::Hold::PfApThm::Term3C0} and using \eqref{Eqn::Hold::PfApThm::PSC0Norm} again we get
\begin{align}\label{Eqn::Hold::PfApThm::BddS}
    \|\Su_{\nu-4}(\Su_\nu \fb\cdot\Su_\nu \gb)\|_{C^0}\le\|\Su_{\nu-4}\|_{C^0\to C^0}\sum_{\substack{j,k\ge l-2\\|j-k|\le5}}\big\|(\De_j\Su_\nu\fb)\cdot (\De_k\Su_\nu\gb)-\De_j\fb\cdot\De_k\gb\big\|_{C^0}\lesssim2^{-\nu(\alpha+\beta-1)}.
\end{align}
This gives the first estimate $\|\Su_{\nu-4}(\Su_\nu \fb\cdot\Su_\nu \gb)\|_{C^0}\lesssim 2^{-\nu(\alpha+\beta-1)}$ for \eqref{Eqn::Hold::PfApThm::Decomp}.

\medskip Next we prove $\|\De_l(\Su_\nu \fb\cdot\Su_\nu \gb)\|_{C^0(\R^n)}\lesssim2^{-\nu(\alpha+\beta-1)}$ when $\nu-3\le l\le \nu+2$. In this case \eqref{Eqn::Hold::PfApThm::Para-1} and \eqref{Eqn::Hold::PfApThm::Para-2} do not vanish.

The term \eqref{Eqn::Hold::PfApThm::Para-1} is still good: by \eqref{Eqn::Hold::PfApThm::PSC0Norm} and Proposition \ref{Prop::Hold::LPchar} we have
\begin{equation}\label{Eqn::Hold::PfApThm::Term1C0}
    \begin{aligned}
    \Big\|\sum_{j=l-2}^{l+2}\big((\De_j\Su_\nu\fb)(\Su_{j-6}\Su_\nu\gb)-\De_j\fb\cdot\Su_{j-6}\gb\big)\Big\|_{C^0}&\lesssim\sum_{j=l-2}^{l+2}\big(\|\Su_\nu\|_{C^0\to C^0}^2+1\big)\|\De_j\fb\|_{C^0}\|\Su_{j-6}\gb\|_{C^0}
    \\
    &\lesssim\sum_{j=l-2}^{l+2}2^{-j\alpha}2^{j(1-\beta)}\approx2^{-l(\alpha+\beta-1)}.
\end{aligned}
\end{equation}

Combining \eqref{Eqn::Hold::PfApThm::Term1C0} and \eqref{Eqn::Hold::PfApThm::Term3C0} we have
\begin{equation}\label{Eqn::Hold::PfApThm::Term0C0}
\begin{aligned}
&\Big\|\De_l\sum_{j=l-2}^{l+2}(\De_j\fb)\cdot(\Su_{j-6}\gb)+\De_l\sum_{j,k\ge l-2;|j-k|\le5}(\De_j\fb)\cdot(\De_k\gb)\Big\|_{C^0}
\\
\lesssim&\|\De_l\|_{C^0\to C^0}\Big(\sum_{j=l-2}^{l+2}2^{-j\alpha}2^{j(1-\beta)}+\sum_{j,k\ge l-2;|j-k|\le5}2^{-j\alpha}2^{k(1-\beta)}\Big)\approx2^{-l(\alpha+\beta-1)}.
\end{aligned}
\end{equation}

By assumption $\nu-3\le l\le\nu+2$, so the right hand side above is bounded by $2^{-\nu(\alpha+\beta-1)}$.

It remains to estimate \eqref{Eqn::Hold::PfApThm::Para-2}. Note that for $|k-l|\le 2$ and $l\le\nu+2$ we have $k\le\nu+4$. By \eqref{Eqn::Hold::LPCharComp} we have $\Su_{k-6}\Su_\nu=\Su_{j-6}$. Therefore we can rewrite the summand in \eqref{Eqn::Hold::PfApThm::Para-2} as, 
\begin{equation}
    \textstyle(\Su_{k-6}\Su_\nu\fb)\cdot(\De_k\Su_\nu\gb)-\Su_{k-6}\fb\cdot\De_k\gb
    =\Su_{k-6}\fb\cdot(\De_k\Su_\nu\gb-\De_k\gb)=\Su_{k-6}\fb\cdot\big(\De_k\sum_{\mu=\nu+1}^{\nu+5}\De_\mu\gb\big).
\end{equation}
Here for the last  equality above, we use the fact $\De_k(\Su_\nu-\id)=\De_k\sum_{\mu=\nu+1}^\infty\De_\mu=\De_k\sum_{\mu=\nu+1}^{\nu+5}\De_\mu$. We define $$\textstyle\psi_\nu:=\sum_{\mu=\nu+1}^{\nu+5}\phi_\mu,\quad\nu\ge0.$$
Thus $\De_k\sum_{\mu=\nu+1}^{\nu+5}\De_\mu\gb=\De_k(\psi_\nu\ast\gb)=\psi_\nu\ast\De_k\gb$ (see Definition \ref{Def::Hold::LPchar}). 
Therefore for $x\in\R^n$,
\begin{equation}\label{Eqn::Hold::PfApThm::Para2Est}
    \begin{aligned}
    &\Su_{k-6}\fb\cdot(\psi_\nu\ast\De_k\gb )(x)=\int_{\R^n}(\Su_{k-6}\fb)(x)\cdot(\De_k\gb )(x-y)\psi_\nu(y)dy
    \\
    =&\big(\psi_\nu\ast(\Su_{k-6}\fb\cdot\De_k\gb)\big)(x)+\int_{\R^n}\psi_\nu(y)\big((\Su_{k-6}\fb)(x)-(\Su_{k-6}\fb)(x-y)\big)\cdot(\De_k\gb )(x-y)dy
    \\
    =&\big(\psi_\nu\ast(\Su_{k-6}\fb\cdot\De_k\gb)\big)(x)+\int_{\R^n}\Big(\int_0^1y\psi_\nu(y)\cdot(\nabla\Su_{k-6}\fb)(x-ty)dt\Big)\cdot(\De_k\gb )(x-y)dy.
\end{aligned}
\end{equation}

By scaling we have  $\|x\phi_\mu(x)\|_{L^1_x}=2^{1-\mu}\|x\phi_1(x)\|_{L^1_x}$ and $\|\phi_\mu\|_{L^1}=\|\phi_1\|_{L^1}$ for all $\mu\ge1$, thus
\begin{equation}\label{Eqn::Hold::PfApThm::PsiNorm}
    \|x\psi_\nu(x)\|_{L^1_x}=2^{-\nu}\|x\psi_0(x)\|_{L^1_x}\approx_\phi2^{-\nu},\quad\|\psi_\nu\|_{L^1}\approx_\phi1,\quad\forall \nu\ge0.
\end{equation}
Combining \eqref{Eqn::Hold::PfApThm::PsiNorm} with \eqref{Eqn::Hold::PfApThm::Para0} and \eqref{Eqn::Hold::PfApThm::Para2Est}, for $\nu-3\le l\le\nu+2$ we have
\begin{align*}
    &\Big\|\De_l\sum_{k=l-2}^{l+2}\big((\Su_{k-6}\Su_\nu\fb)\cdot(\De_k\Su_\nu\gb)-\Su_{k-6}\fb\cdot\De_k\gb\big)\Big\|_{C^0}
    \\
    \le&\Big\|\De_l\sum_{k=l-2}^{l+2}\psi_\nu\ast(\Su_{k-6}\fb\cdot\De_k\gb)\Big\|_{C^0}+\|\De_l\|_{C^0\to C^0}\sum_{k=l-2}^{l+2}\int_{\R^n}|y\psi_\nu(y)|\|\nabla\Su_{k-6}\fb\|_{C^0}\|\De_k\gb\|_{C^0}dy
    \\
    \le&\Big\|\psi_\nu\ast\De_l\Big(\sum_{j=l-2}^{l+2}\De_j\fb\cdot\Su_{j-6}\gb+\sum_{\substack{|j-k|\le5\\j,k\ge l-2}}\De_j\fb\cdot\De_k\gb\Big)\Big\|_{C^0}+\|\phi_l\|_{L^1}\sum_{k=l-2}^{l+2}\|y\psi_\nu\|_{L^1_y}\|\Su_{k-6}\nabla \fb\|_{C^0}\|\De_k\gb\|_{C^0}&\text{by }\eqref{Eqn::Hold::PfApThm::Para0}
    \\
    \lesssim&\|\psi_\nu\|_{L^1}2^{-\nu(\alpha+\beta-1)}+\|y\psi_\nu\|_{L^1_y}\sum_{k=l-2}^{l+2}\sum_{\mu=0}^{k-6}\|\De_\mu\nabla\fb\|_{C^0}\|\De_k\gb\|_{C^0}&\text{by }\eqref{Eqn::Hold::PfApThm::Term0C0}
    \\
    \lesssim&\|\psi_\nu\|_{L^1}2^{-\nu(\alpha+\beta-1)}+\|x\mapsto x\psi_\nu(x)\|_{L^1}\sum_{\mu=0}^\nu\|\nabla\fb\|_{C^{-1,\alpha}}2^{\mu(1-\alpha)}\|\gb\|_{C^{-1,\beta}}2^{l(1-\beta)}&\text{by }\eqref{Eqn::Hold::LPChar}
    \\\lesssim&2^{-\nu(\alpha+\beta-1)}+2^{-\nu}2^{\nu(1-\alpha)}2^{\nu(1-\beta)}\lesssim2^{-\nu(\alpha+\beta-1)}.&\text{by }\eqref{Eqn::Hold::PfApThm::PsiNorm}
\end{align*}
Here for the last inequality we use the fact $\|\nabla\fb\|_{C^{-1,\alpha}}\lesssim\|\fb\|_{C^{0,\alpha}}=1$.

Combining the above estimate with \eqref{Eqn::Hold::PfApThm::Term3C0} and \eqref{Eqn::Hold::PfApThm::Term1C0} we get $\|\De_l(\Su_\nu \fb\cdot\Su_\nu \gb)\|_{C^0(\R^n)}\lesssim2^{-\nu(\alpha+\beta-1)}$ when $\nu-3\le l\le \nu+2$, which is the second estimate for \eqref{Eqn::Hold::PfApThm::Decomp}.

We therefore have $\|\Su_\nu\fb\cdot\Su_\nu\gb\|_{C^0}\lesssim2^{-\nu(\alpha+\beta-1)}$ and finish the proof.
\end{proof}

We can now apply Theorem \ref{Thm::Intro::ApproxThm} to the vector fields with H\"older continuities.

\begin{prop}\label{Prop::Hold::InvVFApt}
Let $m\ge1$ and $0<\alpha,\beta<1$ be satisfy $\alpha>\max(\frac12,1-\beta)$. Let $\Omega\subseteq\R^n$ be an open set. Let $\Omega'\Subset\Omega$ be a precompact open subset and let $\chi\in C_c^\infty(\Omega)$ satisfies $\chi|_{\Omega'}\equiv1$. 

Let $X_1,\dots,X_m\in C^{0,\alpha}_\loc(\Omega;\R^n)$ be vector fields on $\Omega$. Assume there are distributions $(c_{ij}^k)_{i,j,k=1}^m\subset C^{-1,\beta}_\loc(\Omega)$ such that
\begin{equation}\label{Eqn::Hold::InvEqn1}
    [X_i,X_j]=\sum_{k=1}^mc_{ij}^kX_k,\quad 1\le i,j\le m.
\end{equation}

For $\nu\in\N$, we define 
\begin{equation*}
    X^\nu_i:=\Su_\nu(\chi X_i)=(\Su_\nu(\chi a_i^1),\dots,\Su_\nu(\chi a_i^n)),\quad c_{ij}^{k\nu}:=\Su_\nu(\chi c_{ij}^k),\quad 1\le i,j,k\le m.
\end{equation*}

Then there is a constant constant $C$ that does not depend on $\nu$ (but can depend on $X_i,c_{ij}^k,\chi$), such that
\begin{equation}\label{Eqn::Hold::InvVFApt::Main}
    \Big\|[X^\nu_i,X^\nu_j]-\sum_{k=1}^mc_{ij}^{k\nu}X_k^\nu\Big\|_{C^0(\Omega';\R^n)}\le C2^{-\nu\min(\alpha+\beta-1,2\alpha-1)},\quad\forall \nu\ge0.
\end{equation}
\end{prop}

Recall from Lemma \ref{Lem::Hold::MultLoc} \ref{Item::Hold::MultLoc::WellDef} and Corollary \ref{Cor::Hold::[X,Y]WellDef}, $c_{ij}^kX_k\in C^{-1,\beta}_\loc(\Omega;\R^n)$ and $[X_i,X_j]\in C^{-1,\alpha}_\loc(\Omega;\R^n)$ are defined as distributions.

\begin{proof}We can assume $\beta\le\alpha$, hence the right hand side of \eqref{Eqn::Hold::InvVFApt::Main} becomes $C2^{-\nu(\alpha+\beta-1)}$. By computation,
\begin{equation}\label{Eqn::Hold::InvEqn2}
    [\chi X_i,\chi X_j]-(X_i\chi)\cdot \chi X_j+(X_j\chi)\cdot\chi X_i-\sum_{k=1}^m\chi c_{ij}^k\cdot\chi X_k=0\quad\text{on }\R^n,\quad 1\le i,j\le m.
\end{equation}
By writing $X_i=:\sum_{l=1}^na^l\Coorvec{x^l}$ where $a^l\in C^{0,\alpha}_\loc(\Omega)$, \eqref{Eqn::Hold::InvEqn2} becomes
\begin{equation}\label{Eqn::Hold::InvEqn3}
    \sum_{q=1}^n\Big(\chi a_i^q\frac{\partial(\chi a_j^l)}{\partial x^q}-\chi a_j^q\frac{\partial(\chi a_i^l)}{\partial x^q}\Big)+\sum_{q=1}^n\Big(\chi a_i^l\cdot a_j^q\frac{\partial\chi}{\partial x^q}-\chi a_j^l\cdot a_i^q\frac{\partial\chi}{\partial x^q}\Big)-\sum_{k=1}^m\chi c_{ij}^k\cdot\chi a_k^l=0,\ 1\le i,j\le m,\ 1\le l\le n.
\end{equation}

We see that the above equality holds in $C^{-1,\min(\alpha,\beta)}(\R^n)=C^{-1,\beta}(\R^n)$.

In order to apply Theorem \ref{Thm::Intro::ApproxThm}, we define a function $\fb_{ij}^l\in C^{0,\alpha}(\R^n;\R^{2n+m+1})$ and a distribution  $\gb_{ij}^l\in C^{-1,\beta}(\R^n;\R^{2n+m+1})$ for $1\le i,j\le m$ and $1\le l\le n$, as
\begin{equation*}\arraycolsep=0.1pt
    \begin{array}{ccccccccccccccc}
     \fb_{ij}^l:=\big(&\chi a_i^1&,\dots,&\chi a_i^n&,&-\chi a_j^1&,\dots,&-\chi a_j^n&,&-\chi a_1^l&,\dots,&-\chi a_m^l&,&\textstyle\sum_{q=1}^n\big(\chi a_i^la_j^q\partial_q\chi-\chi a_j^l a_i^q\partial_q\chi\big)&\big),
     \\
     \gb_{ij}^l:=\big(&\partial_1(\chi a_j^l)&,\dots,&\partial_n(\chi a_j^l)&,&\partial_1(\chi a_i^l)&,\dots,&\partial_n(\chi a_i^l)&,&\chi c_{ij}^1&,\dots,&\chi c_{ij}^m&,&\centering1&\big).
\end{array}
\end{equation*}

The condition \eqref{Eqn::Hold::InvEqn3} is now the same as saying $\fb_{ij}^l\cdot\gb_{ij}^l=0\in C^{-1,\beta}(\R^n)$ for all $1\le i,j\le m$, $1\le l\le n$. 

By Theorem \ref{Thm::Intro::ApproxThm} we have $\|\Su_\nu\fb_{ij}^l\cdot\Su_\nu\gb_{ij}^l\|_{C^0}\lesssim 2^{-\nu(\alpha+\beta-1)}$ for all $i,j,l$ and $\nu$, in other words,
\begin{equation*}
    \bigg\|[X_i^\nu,X_j^\nu]+\Su_\nu(X_j\chi)\cdot X_i^\nu-\Su_\nu(X_i\chi)\cdot X_j^\nu-\sum_{k=1}^mc_{ij}^{k\nu}X_k^\nu\bigg\|_{C^0(\R^n)}\lesssim2^{-\nu(\alpha+\beta-1)},\quad 1\le i,j\le m,\quad\nu\ge0.
\end{equation*}

To prove \eqref{Eqn::Hold::InvVFApt::Main} it remains to show that $\|\Su_\nu(X_j\chi)\cdot X_i^\nu-\Su_\nu(X_i\chi)\cdot X_j^\nu\|_{C^0(\Omega')}\lesssim 2^{-\nu(\alpha+\beta-1)}$ for all $1\le i,j\le m$. Indeed, by assumption $\chi|_{\Omega'}\equiv1$, so $X_i\chi|_{\Omega'}\equiv0$. Hence for $1\le i\le m$,
\begin{equation}\label{Eqn::Hold::InvVFApt::Tmp}
    \begin{aligned}
    \|\Su_\nu(X_i\chi)\|_{C^0(\Omega')}=&\|\Su_\nu(X_i\chi)-X_i\chi\|_{C^0(\Omega')}\le \|\Su_\nu(X_i\chi)-X_i\chi\|_{C^0(\R^n)}\le\sum_{k=\nu+1}^\infty\|\De_\nu(X_i\chi)\|_{C^0(\R^n)}
    \\\lesssim&\|X_i\chi\|_{C^{0,\alpha}(\R^n)}2^{-\nu\alpha}\lesssim_{\alpha,X,\chi}2^{-\nu(\alpha+\beta-1)}.
\end{aligned}
\end{equation}

Therefore
\begin{align*}
    \sum_{i,j=1}^m\|\Su_\nu(X_j\chi)\cdot X_i^\nu-\Su_\nu(X_i\chi)\cdot X_j^\nu\|_{C^0(\Omega';\R^n)}&\le\sum_{i,j=1}^m\|\Su_\nu(X_j\chi)\|_{C^0(\Omega')}\|X_i^\nu\|_{C^0(\Omega;\R^n)}\\
    &\lesssim 2^{-\nu(\alpha+\beta-1)}\sup_{1\le i\le m}\|\chi X_i\|_{C^0(\R^n)}\lesssim2^{-\nu(\alpha+\beta-1)}.
\end{align*}
This completes the proof.
\end{proof}


Using Proposition \ref{Prop::Hold::InvVFApt} we can prove an estimate similar to \eqref{Eqn::Intro::AsyInv1}.

\begin{cor}\label{Cor::Hold::CorAsyInv}
Following from the assumptions of Proposition \ref{Prop::Hold::InvVFApt}, we have $\lim_{\nu\to\infty}\|X_i^\nu-X_i\|_{C^0(\Omega';\R^n)}=0$ for all $1\le i\le m$, and there is a $t_0>0$ such that
\begin{equation}\label{Eqn::Hold::CorAsyInv::Eqn1}
    \lim\limits_{\nu\to\infty}\max_{1\le i,j\le m}\Big\|[X^\nu_i,X^\nu_j]-\sum_{k=1}^mc_{ij}^{k\nu}X_k^\nu\Big\|_{C^0(\Omega';\R^n)}\exp\Big(t_0\cdot\max_{1\le l\le m}\|\nabla X_l^\nu\|_{C^0(\Omega';\R^{n\times n})}\Big)=0.
\end{equation}

In particular, for the case of two commutative vector fields ($m=2$ and $c_{ij}^k\equiv0$), namely $[X_1,X_2]=0$ in $\Omega$, we have
\begin{equation}\label{Eqn::Hold::CorAsyInv::Eqn2}
        \exists t_0>0,\ \lim\limits_{\nu\to\infty}\|[X_1^\nu,X_2^\nu]\|_{C^0(\Omega';\R^n)}e^{t_0(\|\nabla X_1^\nu\|_{C^0(\Omega';\R^{n\times n})}+\|\nabla X_2^\nu\|_{C^0(\Omega';\R^{n\times n})})}=0.
\end{equation}
\end{cor}
\begin{remark}\label{Rmk::Hold::CorAsyInv}For a distributional involutive log-Lipschitz subbundle $\V\le TM$, by Lemma \ref{Lem::PfThm::CanGen} locally we can find a log-Lipschitz basis $(X_1,\dots,X_r)$ for $\V$ such that $[X_i,X_j]=0$ for all $1\le i,j\le r$. So by \eqref{Eqn::Hold::CorAsyInv::Eqn1} with $c_{ij}^k=c_{ij}^{k\nu}\equiv0$, we see that the condition \eqref{Eqn::Intro::AsyInv1} is satisfied. Thus we conclude that a log-Lipschitz involutive subbundle is always strongly asymptotic involutive.
\end{remark}

In the proof we use the following blow-up estimate: 
\begin{lem}\label{Lem::ODE::LLisAdmissible}
There is a constant $C>0$ depending only on the choice of dyadic resolution $(\phi_j)_{j=0}^\infty$ for the operators $\Su_\nu$ in Definition \ref{Def::Hold::LPchar}, such that
\begin{equation*}
    \|\nabla(\Su_\nu f)\|_{C^0(\R^n;\R^n)}\le C\|f\|_{\Clog(\R^n)}\cdot\nu,\quad\forall \nu\in\Z_+,\quad f\in \Clog(\R^n).
\end{equation*}

\end{lem}
See \cite[Proposition 2.111 and Example p.119]{BahouriCheminDanchin} for a proof. This is based on the fact that the log-Lipschitz modulus of continuity $\mulog(r):=r(1-\log r)$ is an admissible modulus of continuity, see \cite[Definition 2.108]{BahouriCheminDanchin}.

\begin{proof}[Proof of Corollary \ref{Cor::Hold::CorAsyInv}]
By assumption $(X_i^\nu -X_i)|_{\Omega'}=(X_i^\nu -\chi X_i)|_{\Omega'}=(\Su_\nu(\chi X_i)-\chi X_i)|_{\Omega'}$ for all $1\le i\le m$. Thus similar to \eqref{Eqn::Hold::InvVFApt::Tmp},
\begin{equation*}
    \|X_i^\nu-X_i\|_{C^0(\Omega';\R^n)}\le\|\Su_\nu(\chi X_i)-\chi X_i\|_{C^0(\Omega';\R^n)}\le\sum_{\nu+1}^\infty\|\De_\nu(\chi X_i)\|_{C^0(\Omega';\R^n)}\lesssim_\alpha\|\chi X_i\|_{C^{0,\alpha}(\R^n;\R^n)}2^{-\nu\alpha}.
\end{equation*}
Let $\nu\to\infty$ we see that $\lim_{\nu\to\infty}\|X_i^\nu-X_i\|_{C^0(\Omega';\R^n)}=0$.

By Proposition \ref{Prop::Hold::InvVFApt} there is a $C_1>0$ and a $\delta>0$ (in fact $\delta=\min(\alpha+\beta-1,2\alpha-1)$) such that
\begin{equation*}
    \Big\|[X^\nu_i,X^\nu_j]-\sum_{k=1}^mc_{ij}^{k\nu}X_k^\nu\Big\|_{C^0(\Omega';\R^n)}\le C_12^{-\delta \nu},\quad\forall\nu\in\Z_+.
\end{equation*}
By Lemma \ref{Lem::ODE::LLisAdmissible} we can find a $C_2>0$ such that
\begin{equation*}
    \|\nabla X_i^\nu\|_{C^0(\Omega';\R^n)}\le\|\nabla X_i^\nu\|_{C^0(\R^m;\R^n)}\le C_2\nu,\quad\forall \nu\in\Z_+.
\end{equation*}
Therefore for $t_0>0$ we have
\begin{equation}\label{Eqn::Hold::CorAsyInv::Tmp}
    \Big\|[X^\nu_i,X^\nu_j]-\sum_{k=1}^mc_{ij}^{k\nu}X_k^\nu\Big\|_{C^0(\Omega';\R^n)}e^{t_0\cdot\max_{1\le l\le m}\|\nabla X_l^\nu\|_{C^0(\Omega';\R^{n\times n})}}\le C_12^{-\delta\nu}e^{t_0C_2\nu},\quad \nu\in\Z_+.
\end{equation}
    
Taking $t_0=(2C_2)^{-1}\delta\log2$, the right hand side of \eqref{Eqn::Hold::CorAsyInv::Tmp} is $C_12^{-\frac\delta2\nu}$, which goes to $0$ as $\nu\to\infty$. Therefore we get \eqref{Eqn::Hold::CorAsyInv::Eqn1} for such $t_0$.

Replacing $t_0$ by $\frac {t_0}2$ we get \eqref{Eqn::Hold::CorAsyInv::Eqn2}.
\end{proof}

\section{Commutativity of Flows of Log-Lipschitz Vector Fields}
In this section we revisit the ODE regularity estimate for log-Lipschitz vector fields. We use Theorem \ref{Thm::Intro::ApproxThm} to give a positive answer to the flow commuting problem in log-Lipschitz setting, which is the following:
\begin{prop}[The log-Lipschitz flowing commuting]\label{Prop::ODE::FlowComm}Let $ U\subseteq\R^n$ be an open set and let $X$ and $Y$ be two locally log-Lipschitz vector fields.

Suppose $[X,Y]=0$ in the sense of distributions on $ U$, then for any precompact open subset $ U_1\Subset U$ there is a small $\tau_0>0$, such that the flow maps $e^{tX},e^{sY},e^{tX}\circ e^{sY},e^{sY}\circ e^{tX}: U_1\to U$ are all defined for $-\tau_0<t,s<\tau_0$ and
\begin{equation}\label{Eqn::ODE::FlowCommEqn}
    e^{tX}\circ e^{sY}(p)=e^{sX}\circ e^{tY}(p)\in U,\quad\forall p\in U_1,\quad -\tau_0<t,s<\tau_0.
\end{equation}
\end{prop}
Recall that by Corollary \ref{Cor::Hold::[X,Y]WellDef} $[X,Y]$ is a vector field on $ U$ with distributional coefficients, hence it makes sense to say $[X,Y]=0$ as distributions.

In fact \eqref{Eqn::ODE::FlowCommEqn} holds for all $-T<t,s<T$ whenever $T>0$ is a (possibly) large number such that $e^{tX},e^{sY},e^{tX}\circ e^{sY},e^{sY}\circ e^{tX}: U_1\to U$ are all defined. In application the small time commutativity is enough.
\subsection{Estimates of log-Lipschitz ODE flows}

We recall some results in \cite[Chapters 2.11 and 3.1]{BahouriCheminDanchin}.

We define the log-Lipschitz modulus of continuity $\mulog:(0,1]\to(0,\infty)$ and a function $\M_\mulog:(0,1]\to[0,\infty)$ (also see \cite[(3.4)]{BahouriCheminDanchin}) as
\begin{equation}\label{Eqn::ODE::Mulog}
         \mulog(r):=r\log\tfrac er=r(1-\log r),\quad \M_\mulog(r):=\int_r^1\frac{ds}{\mulog(s)}=\log\log \frac er,\quad0< r\le 1.
\end{equation}
Thus \eqref{Eqn::Intro::ClogNorm} can be rewritten as
\begin{equation}\label{Eqn::ODE::ClogNorm}
    \|f\|_{\Clog( U;\R^m)}=\|f\|_{C^0( U;\R^m)}+\sup_{x,y\in  U;0<|x-y|<1}|f(x)-f(y)|\cdot\mulog(|x-y|)^{-1}.
\end{equation}

Clearly $\lim\limits_{r\to0}\M_\mulog(r)=+\infty$, so $\mulog$ satisfies the Osgood condition, namely
\begin{equation*}
    \int_0^1\frac{dr}{\mulog(r)}=+\infty.
\end{equation*}

For a log-Lipschitz vector field $X\in\Clog_\loc( U;\R^n)$ and a point $p\in U$, By Osgood's uniqueness result \cite{Osgood} (see \cite[Theorem 3.2]{BahouriCheminDanchin}),  the autonomous ODE $\dot\gamma(t)=X(\gamma(t))$, $\gamma(0)=p$ has a unique solution. Thus the ODE flow $e^{tX}$ is well-defined in local.

Suppose in addition that $X$ has compact support, by taking the zero extension outside $ U$ we can say $X\in\Clog_c(\R^n;\R^n)$. In this case $e^{tX}(p)\equiv p$ holds for all $t\in\R$ and $p\in U^c$, so by the maximal existence theorem (see \cite[Proposition 3.11]{BahouriCheminDanchin}) we see that $e^{tX}:\R^n\to\R^n$ is defined for all $t\in\R$. Later on we often assume the vector fields to have compact supports in $\R^n$ so that we do not need to worry about the domains of their ODE flows.

In order to give the regularity estimate of $e^{tX}$, we use the following Gronwall-Osgood inequality:
\begin{lem}[{\cite[Lemma 3.4 and Corollary 3.5]{BahouriCheminDanchin}}]\label{Lem::ODE::GroOsgInq}
Let $\mu:(0,1]\to\R_+$ be a Osgood modulus of continuity, that is, $\mu$ is a continuous increasing function satisfying $\lim_{r\to0^+}\mu(r)=0$ and $\int_0^1\frac{dr}{\mu(r)}=+\infty$. Let $\tau\in\R_+$, $\gamma\in L^1([0,\tau];\R_+)$ and $\rho:[0,\tau]\to[0,1]$  satisfy
\begin{equation}\label{Eqn::ODE::GroOsgAss}
    \begin{gathered}
    \rho(0)<1;\qquad\rho(t)\le\rho(0)+\int_0^t\gamma(s)\mu(\rho(s))ds,\quad 0\le s\le\tau;
    \\
    \int_0^\tau\gamma(t)dt\le\M_\mu(\rho(0)),\quad\text{where }\M_\mu(t):=\int_t^1\frac{ds}{\mu(s)}.
\end{gathered}
\end{equation}

Then
\begin{equation}\label{Eqn::ODE::GroOsgCon}
    \rho(t)\le\M_\mu^{-1}\Big(\M_\mu(\rho(0))-\int_0^t\gamma(s)ds\Big),\quad\forall 0\le t\le\tau.
\end{equation}
\end{lem}

For the special case where $\mu=\mulog$ we have the following:
\begin{lem}\label{Lem::ODE::Gronwall}
Let $0<a<1$, $b>0$ and $0<\tau\le\frac1b\log\log\frac ea$. Let $\rho:[0,\tau]\to\R$ be an integrable function such that $\rho(0)=a$ and
\begin{equation}\label{Eqn::ODE::Gronwall::Assumption}
    \rho(t)\le a+b\int_0^t\rho(s)(1-
    \log\rho(s))ds,\quad\forall 0\le t\le \tau.
\end{equation}
Then 
\begin{equation*}
    \rho(t)\le e\cdot a^{e^{-bt}},\quad\forall 0\le t\le\tau.
\end{equation*}
\end{lem}
Note that from the assumption of $\tau$ we have $e\cdot a^{e^{-bt}}\le1$.
\begin{proof}
Take $\gamma(t):\equiv b$ in Lemma \ref{Lem::ODE::GroOsgInq}. Note that by \eqref{Eqn::ODE::Mulog} $\int_0^\tau\gamma(s)ds=b\tau\le\log\log\frac ea=\M_\mulog(\rho(0))$, so the assumptions in \eqref{Eqn::ODE::GroOsgAss} are all satisfied. Thus we have the conclusion \eqref{Eqn::ODE::GroOsgCon}, since
\begin{equation*}
    \rho(t)\le \M_\mulog^{-1}(\M_\mulog(a)-bt)=e^{1-e^{\M_\mulog(a)-bt}}= e^{1-e^{\log\log\frac ea}\cdot e^{-bt}}=e^{1-e^{-bt}}\cdot a^{e^{-bt}}\le e\cdot a^{e^{-bt}},\quad\forall0\le t\le\tau.\qedhere
\end{equation*}
\end{proof}

Lemma \ref{Lem::ODE::Gronwall} establishes the H\"older estimate for log-Lipschitz ODE flows.

\begin{prop}\label{Prop::ODE::ODEReg}
Let $X$ be a log-Lipschitz vector field on $\R^n$ with compact support, and let $ U\subseteq\R^n$ be a bounded open set. We consider the map $\exp_X(t,x):=e^{tX}(x)$ for $t\in\R$ and $p\in\R^n$.
\begin{enumerate}[parsep=-0.3ex,label=(\roman*)]
    \item\label{Item::ODE::ODEReg::BigLL} For any $\eps>0$ there is a small $\tau_0>0$ such that $\exp_X\in C^{0,1-\eps}((-\tau_0,\tau_0)\times  U;\R^n)$. \textnormal{(Also see \cite[Theorem 3.7]{BahouriCheminDanchin})}
    \item\label{Item::ODE::ODEReg::LittleLL} If $X$ is little log-Lipschitz, then   $\exp_X\in C^{0,1^-}((-T,T)\times  U;\R^n)$ for all $T>0$. 
\end{enumerate}
\end{prop}

\begin{proof}
Since $X$ has compact support, we know $\exp_X$ is globally defined. Clearly $\exp_X$ is uniformly Lipschitz in $t$ since $|\exp_X(t_1,x)-\exp_X(t_2,x)|=|e^{(t_1-t_2)X}(e^{t_2X}(x))|\le|t_1-t_2|\|X\|_{C^0(\R^n;\R^n)}$. Thus to prove \ref{Item::ODE::ODEReg::BigLL} it suffices to show the following:
\begin{equation}\label{Eqn::ODE::ODEReg::BigLLClaim}
    \forall\eps>0,\ \exists\tau_0,C_1>0\quad\text{such that }|e^{tX}(x_1)-e^{tX}(x_2)|\le C_1|x_1-x_2|^{1-\eps},\quad\text{for } |t|<\tau_0\text{ and } x_1,x_2\in U.
\end{equation}
And if in addition $X$ is little log-Lipschitz, then to prove \ref{Item::ODE::ODEReg::LittleLL} it suffices to prove the following:
\begin{equation}\label{Eqn::ODE::ODEReg::LittleLLClaim}
    \forall\eps, T>0,\ \exists C_2>0\quad\text{such that }|e^{tX}(x_1)-e^{tX}(x_2)|\le C_2|x_1-x_2|^{1-\eps},\quad\text{for } |t|<T\text{ and } x_1,x_2\in U.
\end{equation}

We now prove \eqref{Eqn::ODE::ODEReg::BigLLClaim} for \ref{Item::ODE::ODEReg::BigLL} and \eqref{Eqn::ODE::ODEReg::LittleLLClaim} for \ref{Item::ODE::ODEReg::LittleLL}.

\medskip
\noindent\ref{Item::ODE::ODEReg::BigLL}: Since $ U$ is a bounded set, by enlarging the constant $C_1$ it suffices to prove \eqref{Eqn::ODE::ODEReg::BigLLClaim} with $|x_1-x_2|<\frac1e$. And by replacing $X$ by $-X$ we only need to prove the case $t\ge0$.

Note that by \eqref{Eqn::ODE::ClogNorm} and \eqref{Eqn::ODE::Mulog}, we have for every $x_1,x_2\in\R^n$ and $t\ge0$,
\begin{equation}
    |e^{tX}x_1-e^{tX}x_2|\le\Big|x_1-x_2+\int_0^t\big(X(e^{sX}x_1)-X(e^{sX}x_2)\big)ds\Big|\le|x_1-x_2|+\|X\|_{\Clog}\int_0^{t}\mulog\big(|e^{sX}x_1-e^{sX}x_2|\big)ds.
\end{equation}
For $x_1,x_2$ that satisfy $|x_1-x_2|<\frac1e$, we take 
\begin{equation*}
    a:=|x_1-x_2|,\quad b:=\|X\|_{\Clog},\quad \rho(t):=|e^{tX}(x_1)-e^{tX}(x_2)|,\quad \tau:=b^{-1}\cdot\log 2\le\tfrac1b\log\log\tfrac ea.
\end{equation*}

Now \eqref{Eqn::ODE::Gronwall::Assumption} is satisfied, and by Lemma \ref{Lem::ODE::Gronwall} we get
$|e^{tX}(x_1)-e^{tX}(x_2)|\le e|x_1-x_2|^{e^{-bt}}$ for all $0\le t\le\tau$ and $|x_1-x_2|<\frac1e$.

Therefore for $\eps>0$, taking $\tau_0=\min(\tau,\frac1b\ln\frac1{1-\eps})$ we get 
\begin{equation*}
    |e^{tX}(x_1)-e^{tX}(x_2)|\le e|x_1-x_2|^{e^{-bt}}\le e|x_1-x_2|^{e^{-b\tau_0}}\le e|x_1-x_2|^{1-\eps},\quad\forall0\le t\le\tau_0\text{ and }|x_1-x_2|<\tfrac1e.
\end{equation*}

Since $U$ is a bounded domain, let $L=\delta^{-1}\operatorname{diam}U$. We then have $ |e^{tX}(x_1)-e^{tX}(x_2)|\le eL|x_1-x_2|^{1-\eps}$ for all $x_1,x_2\in U$. We obtain \eqref{Eqn::ODE::ODEReg::BigLLClaim} and complete the proof of \ref{Item::ODE::ODEReg::BigLL}.

\medskip
\noindent\ref{Item::ODE::ODEReg::LittleLL}:
Fix a small $\eps>0$ and a large $T>0$ from the assumption, the key is to find a $\delta=\delta(\eps, T)>0$ such that $|e^{tX}(x_1)-e^{tX}(x_2)|\le e|x_1-x_2|^{1-\eps}$ for all $|x_1-x_2|<\delta$ and $-T<t<T$. Replacing $X$ by $-X$ it suffices to estimate the case $0\le t<T$.

Let $b:=\frac 1T\log\frac1{1-\eps}>0$. By assumption $X$ is little log-Lipschitz, so there is a $\delta_1>0$ such that
\begin{equation*}
    |X(x_1)-X(x_2)|\le (b/e)^\frac1{1-\eps}|x_1-x_2|(1-\log|x_1-x_2|),\quad\forall x_1,x_2\in\R^n\text{ satisfying }|x_1-x_2|<\delta_1.
\end{equation*}

Take $\delta=\min(e^{1-e^{bT}},\delta_1)$, we have $T\le\frac1b\log\log\frac e\delta$. Now for any $x_1,x_2$ such that $|x_1-x_2|<\delta$, we take 
\begin{equation*}
    a:=|x_1-x_2|,\quad \rho(t):=|e^{tX}(x_1)-e^{tX}(x_2)|,\quad \tau:=T\le\tfrac1b\log\log\tfrac ea.
\end{equation*}
So \eqref{Eqn::ODE::Gronwall::Assumption} is satisfied, and by Lemma \ref{Lem::ODE::Gronwall} we get
$|e^{tX}(x_1)-e^{tX}(x_2)|\le e|x_1-x_2|^{e^{-bt}}\le e|x_1-x_2|^{1-\eps}$ for all $0\le t<T$ and $|x_1-x_2|<\delta$. 

Since $U$ is a bounded domain, let $L=\delta^{-1}\operatorname{diam}U$, we then have $|e^{tX}(x_1)-e^{tX}(x_2)|\le eL |x_1-x_2|^{1-\eps}$ for all $x_1,x_2\in U$. This gives \eqref{Eqn::ODE::ODEReg::LittleLLClaim} and finishes the proof of \ref{Item::ODE::ODEReg::LittleLL}.
\end{proof}
\begin{remark}
In Proposition \ref{Prop::ODE::ODEReg} \ref{Item::ODE::ODEReg::BigLL}, even though $\exp_X\in C^{0,1-\eps}$ near $t=0$ holds for arbitrary $\eps\in(0,1)$, it is possible that $\exp_X$ cannot be $C^{0,1^-}$ near $t=0$. See Proposition \ref{Prop::Further::SharpProp}.
\end{remark}

For a log-Lipschitz vector field $X\in \Clog_c(\R^n;\R^n)$ with compact support, we can consider the Schwartz approximation $X^\nu:=\Su_\nu X$ for $\nu\in\Z_{\ge0}$, where $\Su_\nu$ is the Littlewood-Paley summation operator given in Definition \ref{Def::Hold::LPchar}. The following result shows the uniform convergence $e^{tX^\nu}\to e^{tX}$ when $\nu\to\infty$.

\begin{lem}\label{Lem::ODE::FlowConv}
Let $X$ be a log-Lipschitz vector field with compact support. For $\nu\in\Z_{\ge 0}$, define $X^\nu:=\Su_\nu X$. Then for every $T>0$, we have the following uniform convergence:
\begin{equation}\label{Eqn::ODE::FlowConv}
    \lim\limits_{\nu\to\infty}\sup_{|t|<T,x\in\R^n}|e^{tX^\nu}(x)-e^{tX}(x)|=0.
\end{equation}
\end{lem}
Note that for each $\nu$, the flow $e^{tX^\nu}$ is still globally defined because $X^\nu$ is bounded (even though it does not have compact support).
\begin{proof}
When $x_1,x_2\in\R^n$ satisfy $|x_1-x_2|<1$, by Definition \ref{Def::Hold::LPchar} we have for every $\nu\ge0$,
\begin{equation}\label{Eqn::ODE::FlowConv::Tmp}
    \begin{aligned}
    |X^\nu(x_1)-X^\nu(x_2)|=&|\Su_\nu X(x_1)-\Su_\nu(x_2)|\le\int_{\R^n}2^{n\nu}|\phi_0(2^\nu y)||X(x_1-y)-X(x_2-y)|dy
    \\
    \le&\|2^{n\nu}\phi_0(2^\nu\cdot)\|_{L^1}\|X\|_{\Clog}\cdot\mulog(|x_1-x_2|)=\|\phi_0\|_{L^1}\cdot\|X\|_{\Clog}\cdot\mulog(|x_1-x_2|).
\end{aligned}
\end{equation}
Thus $\|X^\nu\|_{\Clog}\le \|\phi_0\|_{L^1}\|X\|_{\Clog}$ is uniformly bounded for all $\nu\in\Z_{\ge0}$. Since we have uniform convergence $X^\nu\to X$ as $\nu\to\infty$, applying \cite[Proposition 3.9]{BahouriCheminDanchin} we get \eqref{Eqn::ODE::FlowConv}.
\end{proof}

\subsection{Commutativity of flows: Proof of Proposition \ref{Prop::ODE::FlowComm}}

Let $X,Y:U\subseteq\R^n\to\R^n$ be two Lipschitz vector fields such that $[X,Y]=0$ as distributions. We are going to prove $e^{-tX}e^{-sY}e^{tX}e^{sY}=\id$ for small $|t|$ and $|s|$.

\begin{proof}[Proof of Proposition \ref{Prop::ODE::FlowComm}]
Let $ U\subseteq\R^n$ be an open set and let $X,Y: U\to\R^n$ be two log-Lipschitz vector fields such that $[X,Y]=0$ in the sense of distributions. By shrinking $ U$ if necessary, we can assume that $\|X\|_{\Clog( U;\R^n)}$ and $\|Y\|_{\Clog( U;\R^n)}$ are both finite.

Let $ U'\Subset U$ be a precompact subset such that $ U_1\Subset U'$, where $ U_1$ is the precompact subset in the assumption of Proposition \ref{Prop::ODE::FlowComm}. Let $\chi\in C_c^\infty( U)$ be such that $\chi|_{ U'}\equiv1$. Define a small $\tau_1>0$ as $$\tau_1:=\big(4\max(\|X\|_{C^0},\|Y\|_{C^0})\big)^{-1}\cdot\dist( U_1,\partial U').$$
Clearly $\chi X(x)=X(x)$ and $\chi Y(x)=Y(x)$ hold for all $x\in U'\subseteq\{p\in U:\chi(p)=1\}$, so
\begin{equation}\label{Eqn::ODE::PfFlowComm::Flow=}
    e^{-tX}e^{-sY}e^{t'X}e^{s'Y}(x)=e^{-t(\chi X)}e^{-s(\chi Y)}e^{t'(\chi X)}e^{s'(\chi Y)}(x)\in U',\quad\forall x\in U_1,\quad |t|,|s|,|t'|,|s'|<\tau_1.
\end{equation}

Let $X^\nu:=\Su_\nu(\chi X)$ and $Y^\nu:=\Su_\nu(\chi Y)$ be the Schwartz approximations for $\nu\in\Z_+$. We have the commutator formula from \cite[Lemma 4.1]{RampazzoSussmanCommutators}: for $x\in\R^n$ and $t,s\in\R$,
\begin{equation}\label{Eqn::ODE::PfFlowComm::CommFormula}
    e^{-tX^\nu}e^{-sY^\nu}e^{tX^\nu}e^{sY^\nu}(x)-x=\int_0^t\int_0^s\left([Y^\nu,X^\nu](e^{-tX^\nu}e^{-s' Y^\nu}e^{-t' X^\nu})\right)\left(e^{(t-t')X^\nu}e^{s' Y^\nu}(x)\right)dt' ds'.
\end{equation}

We claim that there is a $\tau\in(0,\tau_1]$ such that \eqref{Eqn::ODE::PfFlowComm::CommFormula} converges to $0$ uniformly for $x\in U_1$ and $|t|,|s|<\tau$, that is
\begin{equation}\label{Eqn::ODE::PfFlowComm::Claim}
    \lim_{\nu\to\infty}\sup_{-\tau<t,s<\tau}\|e^{-tX^\nu}e^{-sY^\nu}e^{tX^\nu}e^{sY^\nu}-\id\|_{C^0( U_1;\R^n)}.
\end{equation}

Once this is done, by \eqref{Eqn::ODE::PfFlowComm::Flow=} and Lemma \ref{Lem::ODE::FlowConv} we get for $x\in U_1$ and $|t|,|s|<\tau$,
$$ e^{-tX}e^{-sY}e^{tX}e^{sY}(x)=e^{-t(\chi X)}e^{-s(\chi Y)}e^{t(\chi X)}e^{s(\chi Y)}(x)=\lim\limits_{\nu\to\infty}e^{-tX^\nu}e^{-sY^\nu}e^{tX^\nu}e^{sY^\nu}(x)=x.$$This would complete the proof of the proposition.

\medskip
To prove \eqref{Eqn::ODE::PfFlowComm::Claim}, taking absolute value on both sides of \eqref{Eqn::ODE::PfFlowComm::CommFormula} we have
\begin{align*}
    &\|e^{-tX^\nu}e^{-sY^\nu}e^{tX^\nu}e^{sY^\nu}-\id\|_{C^0( U_1;\R^n)}
    \\\le&\bigg|\int_0^t\int_0^s\|[Y^\nu,X^\nu]\|_{C^0( U';\R^n)}\left\|\nabla\big(e^{-tX^\nu}e^{-s'Y^\nu}e^{-t' X^\nu}\big)\right\|_{C^0( U';\R^{n\times n})}dt'ds'\bigg|
    \\
    \le&\|[X^\nu,Y^\nu]\|_{C^0( U';\R^n)}\int_{-|t|}^{|t|}\int_{-|s|}^{|s|}\|\nabla (e^{-tX^\nu})\|_{C^0(\R^n;\R^{n\times n})}\|\nabla( e^{-s'Y^\nu})\|_{C^0(\R^n;\R^{n\times n})}\|\nabla (e^{-t'X^\nu})\|_{C^0(\R^n;\R^{n\times n})}dt'ds'
    \\
    \le&\|[X^\nu,Y^\nu]\|_{C^0( U';\R^n)}\int_{-|t|}^{|t|}\int_{-|s|}^{|s|}e^{(|t|+|t'|)\|\nabla X^\nu\|_{C^0(\R^n;\R^{n\times n})}}e^{|s'|\|\nabla Y^\nu\|_{C^0(\R^n;\R^{n\times n})}}dt'ds'.
\end{align*}
Here the last inequality is obtained by the Gronwall's inequality for Lipschitz vector fields, see \cite[Proposition 3.10]{BahouriCheminDanchin}.

By Proposition \ref{Prop::Hold::InvVFApt} with $X_1=X$, $X_2=Y$, $c_{ij}^k\equiv0$ and $\alpha=\beta=\frac23$, we can find a $C_1>0$ such that
\begin{equation}\label{Eqn::ODE::PfFlowComm::EstLie}
    \|[X^\nu,Y^\nu]\|_{C^0( U';\R^n)}\le C_12^{-\frac23\nu},\quad\forall \nu\in\Z_+.
\end{equation}

By Lemma \ref{Lem::ODE::LLisAdmissible} we can find a $C_2>0$ that does not depend on $\nu$, such that
\begin{equation}\label{Eqn::ODE::PfFlowComm::EstExp}
    \|\nabla X^\nu\|_{C^0(\R^n;\R^{n\times n})}\le C_2\nu,\quad \|\nabla Y^\nu\|_{C^0(\R^n;\R^{n\times n})}\le C_2\nu,\quad\forall\nu\in\Z_+.
\end{equation}

Combining \eqref{Eqn::ODE::PfFlowComm::EstLie} and \eqref{Eqn::ODE::PfFlowComm::EstExp}, taking $\tau=\min(1,\tau_1,\frac{\log2}{9C_2})$, we see that for $|t|,|s|<\tau$,
\begin{equation*}
    \|e^{-tX^\nu}e^{-sY^\nu}e^{tX^\nu}e^{sY^\nu}-\id\|_{C^0( U_1;\R^n)}\le C_12^{-\frac23\nu}\int_{-1}^1\int_{-1}^1e^{2\tau C_2\nu}e^{\tau C_2\nu}dt'ds'=4C_1 2^{(\frac{3C_2\tau}{\log 2}-\frac23)\nu}=4C_12^{-\frac13\nu}.
\end{equation*}
Let $\nu\to\infty$, the right hand side above goes to 0, finishing the proof of \eqref{Eqn::ODE::PfFlowComm::Claim} and hence the proposition.
\end{proof}

Using Proposition \ref{Prop::ODE::ODEReg} and taking induction on the number of commutative vector fields, we obtain the following:
\begin{cor}\label{Cor::ODE::MultFlow}
Let $r\in\Z_+$ and let $ U\subseteq\R^n$ be an open set. Let $X_1,\dots,X_r$ be (locally) log-Lipschitz vector fields on $ U$ that are pairwise commutative in the sense of distributions, i.e. $[X_i,X_j]=0$ as distributions for $1\le i,j\le r$. Then for any precompact open set $ U_1\Subset U$, there is a $\tau_0>0$ such that the map
\begin{equation}\label{Eqn::ODE::MultFlow::Psi}
    \Psi^X(t,x)=\Psi^X(t^1,\dots,t^r,x):=e^{t^1X_1}\circ\dots\circ e^{t^rX_r}(x),
\end{equation}
is defined on $x\in U_1$ and $|t|<\tau_0$ and satisfies
\begin{equation}\label{Eqn::ODE::MultFlow::Eqnddt}
    \frac{\partial\Psi^X}{\partial t^j}(t,x)=X_j(\Psi^X(t,x)),\quad x\in U_1,\quad|t|<\tau_0,\quad j=1,\dots,r.
\end{equation}
Moreover
\begin{enumerate}[parsep=-0.3ex,label=(\roman*)]
    \item\label{Item::ODE::MultFlow::LLReg} For any $\eps>0$ there is a $0<\tau\le\tau_0$, such that $\Psi^X:B^r(\tau)\times U_1\to U$ is a $C^{0,1-\eps}$-map, and $\Coorvec{t^j}\Psi^X: B^r(\tau)\times U_1\to\R^n$ are $C^{0,1-\eps}$-maps for $1\le j\le r$.
    \item\label{Item::ODE::MultFlow::LittleLLReg} If $X_1,\dots,X_r$ are little log-Lipschitz, then $\Psi^X:B^r(\tau_0)\times U_1\to U$ is a $C^{0,1^-}$-map, and $\Coorvec{t^j}\Psi^X: B^r(\tau_0)\times U_1\to\R^n$ are $C^{0,1^-}$-maps for $1\le j\le r$.
\end{enumerate}
\end{cor}
\begin{proof}
We see that there is a small $\tilde \tau_0>0$ such that the composited ODE flows $e^{t^1X_{j_1}}\circ \dots\circ e^{t^rX_{j_r}}$ are all defined for $j_1,\dots,j_r\in\{1,\dots, r\}$ and $t^1,\dots,t^r\in(-\tilde \tau_0,\tilde\tau_0)$. Indeed by shrinking $ U$ we can assume $X_1,\dots,X_r$ are all bounded vector fields on $ U$, then taking $\tilde \tau_0:=(\|X_1\|_{C^0}+\dots+\|X_r\|_{C^0})^{-1}\cdot\dist( U_1,\partial U)$ will do.

Now $\Psi^X(t,x)$ is defined on $x\in U_1$ and $|t|<\tilde \tau_0$. By Proposition \ref{Prop::ODE::FlowComm} and by choosing $0<\tau_0<\tilde \tau_0$ small enough we have
\begin{align*}
    &\frac{\partial\Psi^X}{\partial t^j}(t,x)=\Coorvec{t^j}e^{t^1X_1}\circ\dots\circ e^{t^rX_r}(x)=\Coorvec{t^j}e^{t^jX_j}\circ e^{t^1X_1}\circ\dots\circ e^{t^{j-1}X_{j-1}}\circ e^{t^{j+1}X_{j+1}}\circ\dots\circ e^{t^rX_r}(x)
    \\
    =&X_j\circ e^{t^jX_j}\circ e^{t^1X_1}\circ\dots\circ e^{t^{j-1}X_{j-1}}\circ e^{t^{j+1}X_{j+1}}\circ\dots\circ e^{t^rX_r}(x)=X_j\circ e^{t^1X_1}\circ\dots\circ e^{t^rX_r}(x)=X_j(\Psi^X(t,x)),
\end{align*}
for all $x\in U_1$ and $|t|<\tau_0$. This completes the proof of \eqref{Eqn::ODE::MultFlow::Eqnddt}.

\medskip\noindent\ref{Item::ODE::MultFlow::LLReg}: By possibly shrinking $\tau_0$ we can assume that $\Psi^X(B^r(0,\tau_0)\times U_1)\Subset U$. Denote $\tilde U_1:=\Psi^X(B^r(0,\tau_0)\times U_1)$. We see that $\Psi^X(t,x)$ is the composition map
\begin{equation}\label{Eqn::ODE::MultFlow::CompMap}
    \Psi^X(t,x)=\exp_{X_1}(t^1,(\exp_{X_2}(t^2,(\dots(\exp_{X_r}(t^r,x))\dots)))),\quad\text{where }\exp_{X_j}: B^r(0,\tau_0)\times\tilde U_1\to U.
\end{equation}

For $0<\eps<1$, take a $\tilde\eps\in(0,\eps]$ such that $(1-\tilde\eps)^{r+1}\ge1-\eps$. By Proposition \ref{Prop::ODE::ODEReg} \ref{Item::ODE::ODEReg::BigLL}, there is a $\tau\in(0,\tau_0]$ such that $\exp_{X_j}\in C^{0,1-\tilde\eps}((-\tau,\tau)\times\tilde U_1;U)$. Therefore $\Psi^X\in C^{0,(1-\tilde\eps)^r}((-\tau,\tau)^r\times U_1;U)$ and thus $\Psi^X\in C^{0,(1-\tilde\eps)^n}(B^r(0,\tau)\times U_1;U)\subseteq C^{0,1-\eps}(B^r(0,\tau)\times U_1;U)$. 

And since $X_j\in\Clog\subset C^{0,1-\tilde\eps}$ for each $1\le j\le r$, we get $X_j\circ\Psi^X\in C^{0,(1-\tilde\eps)^{r+1}}((-\tau,\tau)^r\times U_1;U)\subseteq C^{0,1-\eps}((-\tau,\tau)^r\times U_1;U)$, thus by \eqref{Eqn::ODE::MultFlow::Eqnddt} we have $\frac{\partial\Psi^X}{\partial t^j}=X_j\circ\Psi^X\in C^{0,1-\eps}(B^r(0,\tau)\times U_1;U)$.

\medskip\noindent\ref{Item::ODE::MultFlow::LittleLLReg}: When $X_1,\dots,X_r\in\clog$, by Proposition \ref{Prop::ODE::ODEReg} \ref{Item::ODE::ODEReg::LittleLL} we have $\exp_{X_j}\in C^{0,1^-}((-\tau_0,\tau_0)\times\tilde U_1;U_1)$ for $1\le j\le r$. Taking compositions and using \eqref{Eqn::ODE::MultFlow::CompMap} we get $\Psi^X\in C^{0,1^-}((-\tau_0,\tau_0)\times U_1;U)$, thus $\Psi^X\in C^{0,1^-}(B^r(0,\tau)\times U_1;U)$. 

Since $X_j\in\clog\subset C^{0,1^-}$, by \eqref{Eqn::ODE::MultFlow::Eqnddt} we get $\frac{\partial\Psi^X}{\partial t^j}=X_j\circ\Psi^X\in C^{0,1^-}(B^r(0,\tau_0)\times U_1;U)$ for $1\le j\le r$.
\end{proof}
\section{Proofs of Main Results (Theorems \ref{Thm::MainThm1} and \ref{Thm::MainThm2})}\label{Section::PfThm}
We begin the proof by taking a good local basis following the idea of \cite{ShortFro}.

\begin{lem}[Canonical local basis]\label{Lem::PfThm::CanGen} Let $\tilde U\subseteq\R^n$ be an open set and let $\V$ be a continuous tangent subbundle over $\tilde U$ with rank $r$.
Then for any $p\in\tilde  U$ there are a linear coordinate system $(x,y)=(x^1,\dots,x^r,y^1,\dots,y^{n-r})$ for $\R^n$ and a neighborhood $ U\subseteq\tilde U$ of $p$, such that $\V|_ U$ has continuous generators $X_1,\dots,X_r$ of the form
\begin{equation}\label{Eqn::PfThm::CanGen::Main}
    X_j=\Coorvec{x^j}+\sum_{k=1}^{n-r}b_j^k\Coorvec{y^k},\qquad j=1,\dots,r,
\end{equation}
where $b_j^k\in C^0_\loc( U)$ satisfies $b_j^k(p)=0$ for $1\le j,k\le n$. 
Moreover
\begin{enumerate}[parsep=-0.3ex,label=(\roman*)]
    \item\label{Item::PfThm::CanGen::Reg} If $\V$ is $\Clog$, $\clog$ or $C^{0,\alpha}$ ($0<\alpha<1$), then $X_1,\dots,X_r$ are also $\Clog$, $\clog$ or $C^{0,\alpha}$ ($0<\alpha<1$) respectively.
    \item\label{Item::PfThm::CanGen::Inv} If $\V$ is log-Lipschitz and distributional involutive, then $X_1,\dots,X_r$ are pairwise commutative as distributions.
\end{enumerate}
\end{lem}


\begin{proof}The existence of $X_1,\dots,X_r$ and the result \ref{Item::PfThm::CanGen::Reg} can be obtained using the same argument of \cite[Lemma 1]{ShortFro}. For this method we omit the details to reader.

Alternatively the subbundle $\V$ can be viewed as the map $\V:\tilde U\to\Gr(r,\R^n)$ where $\Gr(r,\R^n)=\mathbf{O}(\R^n)/(\mathbf{O}(\R^r)\times\mathbf{O}(\R^{n-r}))$ is the Grassmannian space containing all rank $r$ linear subspaces. Let $(\rho^1,\dots,\rho^n)$ be the standard coordinates for $\R^n$. We know $\Gr(r,\R^n)$ is covered by the standard affine coordinate charts $\{\phi_I:\mathbf{U}_I\subset\Gr(r,\R^n)\to\R^r\times\R^{n-r}\}_{I\subseteq\{1,\dots,n\};\#I=r}$, where
\begin{equation*}
     \phi_I^{-1}(u):=\Span\Big(\Coorvec{\rho^{i_a}}+\sum_{b=1}^{n-r}u_a^b\Coorvec{\rho^{i_{r+b}}}\Big)_{a=1}^r,\quad \text{for }u=(u_a^b)_{r\times(n-r)}\in\R^r\times\R^{n-r};\quad \mathbf U_I:=\phi_I^{-1}(\R^r\times\R^{n-r}).
\end{equation*}
Here we use $1\le i_1<\dots<i_r\le n$ as indices in $I$, and $1\le i_{r+1}<\dots<i_n\le n$ as indices not in $I$.

Say $\V(p)\in\mathbf U_I$ for some $|I|=r$, we can take $x^j=\rho^{i_j}$ for $1\le j\le r$ and $y^j=\rho^{i_{r+j}}$ for $1\le j\le n-r$. By continuity $\V(q)\in\mathbf U_I$ holds for $q\in\tilde U$ closed to $p$. Thus we can find a such $U\subseteq\tilde U$ with $(b_j^k)_{r\times (n-r)}:=\phi_I\circ\V|_U$. 

Set $\Xs\in\{\Clog,\clog,C^{0,\alpha}\}$. If $\V\in\Xs$, we see that $\phi_I\circ\V|_U:U\to\R^r\times\R^{n-r}$ is in $\Xs$. Thus $(b_j^k)_{r\times (n-r)}:U\to \R^r\times\R^{n-r}$ have the same regularity, which gives \ref{Item::PfThm::CanGen::Reg}.

\medskip\noindent\ref{Item::PfThm::CanGen::Inv}: The proof is slightly different from \cite[Corollary]{ShortFro}. Since $X_1,\dots,X_r$ form a local basis of $\V$ on $ U$, $\V^\bot|_ U$ has sections
\begin{equation}\label{Eqn::PfThm::CanGen::DualForm}
    \theta^l:=dy^l-\sum_{j=1}^rb_j^ldx^j,\quad 1\le l\le n-r.
\end{equation}

On the other hand by direct computation 
\begin{equation}\label{Eqn::PfThm::CanGen::LieBraEqn}
    [X_j,X_k]=\sum\limits_{l=1}^{n-r}\Big(\frac{\partial b_k^l}{\partial x^j}-\frac{\partial b_j^l}{\partial x^k}+\sum\limits_{m=1}^{n-r}\Big(b_j^m\frac{\partial b_k^l}{\partial y^m}-b_k^m\frac{\partial b_j^l}{\partial y^m}\Big)\Big)\Coorvec{y^l}=\sum_{l=1}^{n-r}\langle dy^l,[X_j,X_k]\rangle\Coorvec{y^l}=\sum_{l=1}^{n-r}\langle \theta^l,[X_j,X_k]\rangle\Coorvec{y^l}.
\end{equation}

By assumption (see Definition \ref{Def::Intro::DisInv}) $\langle\theta^l,[X_j,X_k]\rangle=0$ as distributions, so $[X_j,X_k]=0$ as distributions for all $1\le j,k\le r$.
\end{proof}

We now begin the proof of Theorems \ref{Thm::MainThm1} and \ref{Thm::MainThm2}. It is enough to assume $M$ to be an open subset of $\R^n$, see Remark \ref{Rmk::DisInv::RmkMfldObj} \ref{Item::DisInv::RmkMfldObj::CheckCover}. We postpone the discussion of distributions on $C^{1,1}$-manifolds to Section \ref{Section::DisInv::DisDef}.

\begin{proof}[Proof of Theorem \ref{Thm::MainThm1}]
By passing to a local coordinate chart of the manifold $M$, we can assume that $M\subseteq\R^n$ is an open subset with the base point $p=0\in\R^n$.

By Lemma \ref{Lem::PfThm::CanGen} we can find a neighborhood $U\subseteq M$ of $p=0$, a linear coordinate system $(x,y)=(x^1,\dots,x^r,y^1,\dots,y^{n-r})$ for $\R^n$, and commutative log-Lipschitz vector fields $X_1,\dots,X_r$ defined on $U$ that have the form \eqref{Eqn::PfThm::CanGen::Main} and span $\V$ at every point in $U$.

Let $\Omega''\subseteq\R^{n-r}$ be an open neighborhood of $0\in\R^{n-r}$ such that $\{(0,v)\in\R^r\times\R^{n-r}:v\in\Omega''\}\Subset U$. We can take a small enough number $\tau_0>0$, so the following map $\Phi(u,v):B^r(0,\tau_0)\times\Omega''\subset\R^r\times\R^{n-r}\to U$ is defined:
\begin{equation}\label{Eqn::PfThm::DefofPhi}
    \Phi(u,v):=\Phi(u^1,\dots,u^r,v):=e^{u^1X_1}\circ\dots\circ e^{u^rX_r}(0,v),\quad |u|<\tau_0,\quad v\in\Omega''.
\end{equation}
Clearly $\Phi(0,0)=0$, this gives \ref{Item::MainThm1::Phi0}.

\medskip
By Corollary \ref{Cor::ODE::MultFlow} and by choosing $\tau_0>0$ small enough, $\Phi$ is well-defined, differentiable in $u$, and satisfies
\begin{equation*}
    \Coorvec{u^j}\Phi(u,v)=X_j(\Phi(u,v)),\quad1\le j\le r,\quad|u|<\tau_0,\quad v\in\Omega''.
\end{equation*}

Since $X_1,\dots,X_r$ are linearly independent and span $\V$ in the domain, we see that $\frac{\partial\Phi}{\partial u^1}(u,v),\dots,\frac{\partial\Phi}{\partial u^r}(u,v)$ spans $\V_{\Phi(u,v)}$ for $(u,v)$ in the domain $B^r(0,\tau_0)\times\Omega''$. This gives \ref{Item::MainThm1::Span}.

\medskip
Note that $\Phi(u,v)=\Psi^X(u,(0,v))$ where $\Psi^X$ is the map in \eqref{Eqn::ODE::MultFlow::Psi}. For a given $\eps>0$, by Corollary \ref{Cor::ODE::MultFlow} \ref{Item::ODE::MultFlow::LLReg} we can find a $\tau\in(0,\tau_0]$ such that $\Phi\in C^{0,1-\eps}(B^r(0,\tau)\times\Omega'';U)$ and $\frac{\partial\Phi}{\partial u^j}\in C^{0,1-\eps}(B^r(0,\tau)\times\Omega'';\R^n)$ for $1\le j\le r$.

Therefore we can take $\Omega:=B^r(0,\tau)\times\Omega''$. Since we have identified $M$ as a subset of $\R^n$, we see that $\Phi\in C^{0,1-\eps}(\Omega;M)$ and $\frac{\partial\Phi}{\partial u^j}\in C^{0,1-\eps}(\Omega;TM)$ for $1\le j\le r$. This proves \ref{Item::MainThm1::PhiReg}. 

\medskip
It remains to show \ref{Item::MainThm1::Phi-1} that $\Phi:\Omega\to\Phi(\Omega)$ is homeomorphism and $\Phi^{-1}\in C^{0,1-\eps}(\Phi(\Omega);\R^n)$. 

We can assume $(x,y)$ to be the standard coordinate system for $M\subseteq\R^n$. Thus we have the natural coordinate maps $ x^j(q^1,\dots,q^n)=q^j$ for $1\le j\le r$ and $ y^k(q^1,\dots,q^n)=q^{r+k}$ for $1\le k\le n-r$ and $q\in\R^n$. Define $\lambda=(\lambda^1,\dots\lambda^{n-r}):\Phi(\Omega)\subset\R^n_q\to\R^{n-r}$ as 
\begin{equation}\label{Eqn::PfThm::lambda}
    \lambda^j(q):=y^j(e^{-q^rX_r}\dots e^{-q^1X_1}(q)),\quad 1\le j\le n-r.
\end{equation}
By Corollary \ref{Cor::ODE::MultFlow} \ref{Item::ODE::MultFlow::LLReg} with possibly shrinking $\Omega$, we have that $\lambda^1,\dots,\lambda^{n-r}\in C^{1-\eps}(\Phi(\Omega))$. So $(x,\lambda)$ is a $C^{0,1-\eps}$  map defined on $\Phi(\Omega)$.

We are going to show $\Phi\circ(x,\lambda)=( x,\lambda)\circ\Phi=\id$. This implies that  $\Phi\in C^{0,1-\eps}$ is homeomorphic onto its image such that $\Phi^{-1}\in C^{0,1-\eps}$ as well. In particular $\Phi(\Omega)$ is an open set.

Indeed, since $X_j$ are of the form $X_j=\Coorvec{x^j}+\sum_{k=1}^{n-r}b_j^k\Coorvec{y^k}$, we see that $\Phi$ can be written as
$$\Phi(u,v)=(u,\phi(u,v)),\quad\text{where }\phi\in C^{0,1-\eps}(\Omega;\R^{n-r}).$$
Also we have
\begin{equation}\label{Eqn::PfThm::Tmp}
    e^{-q^rX_r}\dots e^{-q^1X_1}(q)=(0, y(e^{-q^rX_r}\dots e^{-q^1X_1}(q))),\quad  x(e^{-q^rX_r}\dots e^{-q^1X_1}(q))=0,\quad q\in\Phi(\Omega).
\end{equation}

Therefore for $(u,v)\in\Omega$, using $\Phi(u,v)=(u,\phi(u,v))$ we have
\begin{equation}\label{Eqn::PfThm::PfInv1}
    \begin{aligned}
    ( x,\lambda)(\Phi(u,v))=&(x(u,\phi(u,v)),y(e^{-u^rX_r}\dots e^{-u^1X_1}(\Phi(u,v))))
    \\
    =&(u, y(e^{-u^rX_r}\dots e^{-u^1X_1}e^{u^1X_1}\dots e^{u^rX_r}(0,v)))
    \\
    =&(u, y(0,v))=(u,v).
\end{aligned}
\end{equation}
For $(x,y)\in\Phi(\Omega)$, by \eqref{Eqn::PfThm::Tmp} we have
\begin{equation}\label{Eqn::PfThm::PfInv2}
    \begin{aligned}
    \Phi((x(q),\lambda(q))=&\Phi(x,y(e^{-q^rX_r}\dots e^{-q^1X_1}(q)))
    \\
    =&e^{q^1X_1}\dots e^{q^rX_r}(0, y(e^{-q^rX_r}\dots e^{-q^1X_1}(q)))
    \\
    =&e^{q^1X_1}\dots e^{q^rX_r}(e^{-q^rX_r}\dots e^{-q^1X_1}(q))=q.
\end{aligned}
\end{equation}
We get $\Phi\circ(x,\lambda)=\id$ and $(x,\lambda)\circ\Phi=\id$, finishing the whole proof.
\end{proof}

To prove Theorem \ref{Thm::MainThm2}, we follows the same construction of $\Phi$ in \eqref{Eqn::PfThm::DefofPhi}.
\begin{proof}[Proof of Theorem \ref{Thm::MainThm2}]
We can still assume $M$ is an open subset of $\R^n$. Let $X_1,\dots,X_r$ and $B^r(0,\tau_0)\times\Omega''$ be as in the proof of Theorem \ref{Thm::MainThm1}.
Set $\Omega:=B^r(0,\tau_0)\times\Omega''$, define $\Phi:\Omega\to M$ to be the same as \eqref{Eqn::PfThm::DefofPhi} and $\lambda^1,\dots,\lambda^r:\Phi(\Omega)\to\R$ to be as \eqref{Eqn::PfThm::lambda}. Following the proof of Theorem \ref{Thm::MainThm1} we get \ref{Item::MainThm2::Phi0}, \ref{Item::MainThm2::PhiReg} and \ref{Item::MainThm2::Span} without any change. Also we know $\Phi^{-1}=( x,\lambda):\Phi(\Omega)\to\R^n$ by \eqref{Eqn::PfThm::PfInv1} and \eqref{Eqn::PfThm::PfInv2}.

The only thing we need to prove is \ref{Item::MainThm1::Phi-1} that $\Phi\in C^{0,1^-}(\Omega;\R^n)$ and $\lambda\in C^{0,1^-}(\Phi(\Omega);\R^{n-r})$. This result follows immediately from Corollary \ref{Cor::ODE::MultFlow} \ref{Item::ODE::MultFlow::LittleLLReg} since we have $\Phi(u,v)=\Psi^X(u,(0,v))$ and $(0,\lambda(x,y))=\Psi^X(-x,(x,y))$, where $\Psi^X$ is in \eqref{Eqn::ODE::MultFlow::Psi}.
\end{proof}
\section{Distributional Characterizations of Sections and Involutivities}\label{Section::DisInv}
In this section we characterize sections and involutivity in the sense of distributions.

\begin{conv}\label{Conv::DisInv::-1log}
We use $C^{0,\mathrm{LogL}}=\Clog$ and we use the formal ordering $0<1-\mathrm{LogL}<0^+<\alpha<1^-<\mathrm{LogL}<1$ for all $0<\alpha<1$. In this way for $\beta,\gamma\in[0,1)\cup\{0^+,1^-,\mathrm{LogL} \}$, $\beta\le\gamma$ if and only if $C^{0,\beta}\supseteq C^{0,\gamma}$.

For an open subset $\Omega\subseteq\R^n$ we define
\begin{equation*}
    \textstyle C^{-1,\mathrm{LogL}}_\loc(\Omega):=\big\{\sum_{j=1}^N\sum_{i=1}^nf_{ij}\cdot\partial_ig_j:N\in\Z_+,\ f_{ij}\in C^{0,0^+}_\loc(\Omega),\ g_j\in \Clog_\loc(\Omega)\big\}.
\end{equation*}
\end{conv}
\begin{remark}
The definition of $C^{-1,\mathrm{LogL}}$ guarantees the product map $[(f,g)\mapsto fg]:C^{0,\alpha}_\loc(\Omega)\times C^{-1,\alpha}_\loc(\Omega)\to C^{-1,\alpha}_\loc(\Omega)$ is defined for $\alpha=\mathrm{LogL}$. Recall from Lemma \ref{Lem::Hold::MultLoc} \ref{Item::Hold::MultLoc::WellDef} this is true for $\frac12<\alpha<1$.

\end{remark}
\begin{remark}
    In \cite{LidingThesis} a different space $\mathscr C^{\mathrm{LogL}-1}$ was used. We have $C^{-1,\mathrm{LogL}}\subseteq \mathscr C^{\mathrm{LogL}-1}$, and the inclusion is probably strict.
    
    In \cite[Definition 2.1.21]{LidingThesis} we define $\|f\|_{\mathscr  C^{\mathrm{LogL}-1}(\R^n)}:=\sup_{j\ge1}j^{-1}\|\Su_jf\|_{L^\infty(\R^n)}$ where $\Su_j$ are in Definition \ref{Def::Hold::LPchar}. For an open subset $\Omega\subseteq\R^n$ we define $\|f\|_{\mathscr 
 C^{\mathrm{LogL}-1}(\Omega)}=\inf\{\|\tilde f\|_{\mathscr  C^{\mathrm{LogL}-1}(\R^n)}:\tilde f|_\Omega=f\}$ and $\mathscr 
 C^{\mathrm{LogL}-1}_\loc(\Omega)=\{f\in\D'(\Omega):f|_{\Omega'}\in \mathscr 
 C^{\mathrm{LogL}-1}(\Omega')\text{ for every open }\Omega'\Subset\Omega\}$.

    Using this definition, Lemma \ref{Lem::DisInv::PushForwardFuncSpaces} below also holds if one replaces $C^{-1,\beta}$ by $\mathscr 
 C^{\mathrm{LogL}-1}$, see \cite[Lemmas 2.1.22 (ii) and 2.1.36]{LidingThesis}. In particular $[(f,g)\mapsto fg]:\Clog_\loc(\Omega)\times \mathscr 
 C^{\mathrm{LogL}-1}_\loc(\Omega)\to \mathscr 
 C^{\mathrm{LogL}-1}_\loc(\Omega)$ holds. One can see that a $\mathscr 
 C^{\mathrm{LogL}-1}$-vector field can be defined on a $C^{1,1}$-manifold, see also \cite[Remark 3.1.2 (a)]{LidingThesis}.
\end{remark}


\subsection{Distributional sections of subbundles}\label{Section::DisInv::DisDef}
We first explain how (some) distributions can be defined on non-smooth manifolds. 

For $0<\beta\le1$, the $C^{-1,\beta}$-functions can be defined on $C^{1,1}$-manifolds due to the following:
\begin{lem}\label{Lem::DisInv::PushForwardFuncSpaces}
Let $U,V\subseteq\R^n$ be two open sets. Let $\varphi:U\to V$ be a $C^{1,1}$-diffeomorphism (i.e. $\varphi,\varphi^{-1}$ are both $C^{1,1}_\loc$-maps). Let $\alpha\in(0,1]\cup\{\mathrm{LogL}\}$, $\beta\in(1-\alpha,1)\cup\{\mathrm{LogL}\}$ and let $\rho\in C^{0,\alpha}_\loc(U)$.

If $f\in C^{-1,\beta}_\loc(V)$ then $\rho(f\circ\varphi)\in C^{-1,\beta}_\loc(U)$ as well.
\end{lem}
\begin{proof}
Write $f=g_0+\sum_{j=1}^n\partial_jg_j$ where $g_0,\dots,g_n\in C^{0,\beta}_\loc(V)$. Clearly $g_j\circ\varphi\in C^{0,\beta}_\loc(U)$ for $j=0,\dots,n$.

Denote $\psi=(\psi^1,\dots,\psi^n):=\varphi^{-1}:V\to U$, which is a $C^{1,1}$-map by assumption. Using chain rule and product rule,
\begin{equation*}
    \rho((\partial_jg_j)\circ\varphi)=\sum_{k=1}^n(\partial_k(g_j\circ\varphi))\cdot\rho((\partial_j\psi^k)\circ\varphi).
\end{equation*}
We have $g_j\circ\varphi\in C^{0,\beta}_\loc(U)$ and $(\partial_j\psi^k)\circ\varphi\in C^{0,1}_\loc(U)$, so $\partial_k(g_j\circ\varphi)\in C^{-1,\beta}_\loc(U)$, $\rho((\partial_j\psi^k)\circ\varphi)\in C^{0,\alpha}_\loc(U)$.
Thus by Lemma \ref{Lem::Hold::MultLoc} \ref{Item::Hold::MultLoc::WellDef} for $\beta\in (1-\alpha,1)$ and Convention \ref{Conv::DisInv::-1log} for $\beta=\mathrm{LogL}$, we have the product $(\partial_k(g_j\circ\varphi))\cdot\rho((\partial_j\psi^k)\circ\varphi)\in C^{-1,\beta}_\loc(U)$. Taking sum over $k=1,\dots,n$ we get $\rho((\partial_jg_j)\circ\varphi)\in C^{-1,\beta}_\loc(U)$ for $j=1,\dots,n$.


Clearly $\rho(g_0\circ\varphi)\in C^{0,\min(\alpha,\beta)}_\loc(U)\subset C^{-1,\beta}_\loc(U)$. Taking sum over $j=0,\dots,n$ we get the result.
\end{proof}

Using Lemma \ref{Lem::DisInv::PushForwardFuncSpaces}, we can define objects on manifolds via patchings, which are  analogous to \cite[Definition 6.3.3]{Hormander}.

\begin{defn}\label{Def::DisInv::DefFunVF}
 Let  $M$ be a $n$ dimensional $C^{1,1}$-manifold with maximal $C^{1,1}$-atlas $\As=\{\psi:U_\psi\subseteq M\to\R^n\}_\psi$. Namely, each $\psi\in\As$ is a homeomorphism $\psi=(\psi^1,\dots,\psi^n):U_\psi\xrightarrow{\sim}\psi(U_\psi)\subseteq\R^n$; we have $\psi\circ\phi^{-1}\in C^{1,1}_\loc(\phi(U_\phi\cap U_\psi);\R^n)$ whenever $\psi,\phi\in\As$ satisfy $U_\phi\cap U_\psi\neq \varnothing$; and $\As$ is maximal with these properties.
 
 Let $\beta\in(0,1)$ and $m\in\{0,-1\}$.
\begin{enumerate}[parsep=-0.3ex,label=(\roman*)]
    \item A (locally) $ C^{m,\beta}$-function $f$ is a collection $\{f_\psi\in C^{m,\beta}_\loc(\psi(U_\psi)):\psi\in\As\}$ such that
    \begin{equation}\label{Eqn::DisInv::DefFunVF::TransFun}
        f_\psi=f_\phi\circ(\phi\circ\psi^{-1}),\quad\text{on }\psi(U_\psi\cap U_\phi)\subseteq\R^n\text{ whenever }U_\psi\cap U_\phi\neq\varnothing.
    \end{equation}
    
    We denote by $ C^{m,\beta}_\loc(M)$ the space of all $ C^{m,\beta}$-functions on $M$.
    \item\label{Item::DisInv::DefFunVF::VF} A (locally) $ C^{m,\beta}$-vector field $X$ is a collection $\{X_\psi=(X_\psi^1,\dots,X_\psi^n)\in C^{m,\beta}_\loc(\psi(U_\psi);\R^n):\psi\in\As\}$ such that
    \begin{equation}\label{Eqn::DisInv::DefFunVF::TransVF}
        X_\psi^i=\sum_{j=1}^n\big(X_\phi^j\circ (\phi\circ \psi^{-1})\big)\cdot\Big(\frac{\partial \psi^i }{\partial \phi^j}\circ \psi^{-1}\Big),\quad\text{on }\psi(U_\psi\cap U_\phi)\subseteq\R^n\text{ whenever }U_\psi\cap U_\phi\neq\varnothing.
    \end{equation}
    We denote by $ C^{m,\beta}_\loc(M;\X)$ the space of all $ C^{m,\beta}$-vector fields on $M$.
\end{enumerate}
\end{defn}

\begin{remark}\label{Rmk::DisInv::RmkMfldObj}
\begin{enumerate}[parsep=-0.3ex,label=(\alph*)]
    \item Here \eqref{Eqn::DisInv::DefFunVF::TransFun} and \eqref{Eqn::DisInv::DefFunVF::TransVF} are given by the pullback equations $f_\psi=(\phi\circ\psi^{-1})^*f_\phi$ and $X_\psi=(\phi\circ\psi^{-1})^*X_\phi$. By Lemma \ref{Lem::DisInv::PushForwardFuncSpaces} the right hand sides of \eqref{Eqn::DisInv::DefFunVF::TransFun} and \eqref{Eqn::DisInv::DefFunVF::TransVF} are well-defined $C^{-1,\beta}_\loc$ functions. 
    \item\label{Item::DisInv::RmkMfldObj::CheckCover} In particular, to show a function $f$ (or a vector fields $X$) is $C^{-1,\beta}$, it suffices to show $f_{\psi_j}\in C^{-1,\beta}_\loc(\psi_j(U_{\psi_j}))$ (or $X_{\psi_j}\in C^{-1,\beta}_\loc(\psi_j(U_{\psi_j});\R^n)$) on a coordinate cover $\{\psi_j\}_j$. Namely a subcollection of atlas $\{\psi_j:U_{\psi_j}\to\R^n\}_j$ such that $\bigcup_jU_{\psi_j}=M$. In this way to study the local problems we can use one coordinate chart near a fixed point and work on subsets in Euclidean spaces.
    \item We do not talk about ``bounded $C^{m,\beta}$ functions or vector fields'' on a general manifold, since a manifold may not have structures like metric to control the size of functions or vector fields.
\end{enumerate}
\end{remark}
By investigating on the transition maps we can also define $C^{-1,\alpha}$-differential forms on $M$ in the same way. We omit the details to readers.

\begin{note}
Let $M$ be a $n$-dimensional $C^{1,1}$ manifold. For $\alpha\in(0,1]\cup\{\mathrm{LogL}\}$ and $m\in\{0,-1\}$, $0\le k\le n$, we denote by $C^{m,\alpha}_\loc(M;\Lambda^k)$ the space of all (locally) $C^{m,\alpha}$ differential $k$-forms on $M$. And we denote $C^{m,\alpha}_\loc(M;\Lambda^\bullet):=\bigoplus_{k=0}^nC^{m,\alpha}_\loc(M;\Lambda^k)$.
\end{note}

\begin{remark}\label{Rmk::DisInv::Functional}
A $C^{m,\beta}$-function on manifold can be identified as a linear functional by the way of \cite[Chapter 3.1]{GeomDist}. We sketch the construction as the following: 

For a (smooth or non-smooth) $n$-dimensional manifold $M$ with atlas $\As=\{\psi:U_\psi\subseteq M\to\R^n\}$, we define the \textbf{volume bundle} $\Vol_M$ to be a real rank 1 vector bundle on $M$ as
\begin{gather}\label{Eqn::DisInv::DefVolBundle}
    \Vol_M:=\Big(\coprod_{\psi\in\As}\psi(U_\psi)\times\R\Big)\Big/\left\{(x,t)\sim(\phi\circ\psi^{-1}(x),g_{\phi\psi}(x)t):x\in \psi(U_\phi\cap U_\psi),t\in\R\right\},
    \\
    \text{where }g_{\phi\psi}(x):=|\det D(\phi\circ\psi^{-1})(x)|,\quad x\in \psi(U_\phi\cap U_\psi).\notag
\end{gather}

When $(M,\As)$ is a $C^{1,1}$-manifold, then $\Vol_M$ is a Lipschitz vector bundle. In this case take a $C^{1,1}$ coordinate cover $\{\psi_i:U_{\psi_i}\subset M\to\R^n\}_{i\in I}$ of $M$ and take a $C^{1,1}$ partition of unity $\{\chi_i\in C_\loc^{1,1}(U_{\psi_i})\}_{i\in I}$, we have a pairing: for $f=(f_\psi)_{\psi\in\As}\in C^0_\loc(M)$,
\begin{equation*}
    \langle f,s\rangle:=\sum_{i\in I}\int_{\psi_i(U_{\psi_i})}\chi_i(\psi_i(x))f_{\psi_i}(x)(\psi_i^*s)(x)dx,\quad s\in C^{0,1}_c(M;\Vol_M).
\end{equation*}
The pairing is still defined for $f\in C^{-1,\beta}_\loc(M)$, thus we have a map $C^{-1,\beta}_\loc(M)\hookrightarrow C_c^{0,1}(M;\Vol_M)'$ for $0<\beta<1$. The pairing is canonical in the sense that it does not depend on the choice of coordinates and partition of unity. Therefore we can say $C^{-1,\beta}_\loc(M)\subset C_c^{0,1}(M;\Vol_M)'$ for all $0<\beta<1$. 

Similarly we have the embedding of the space of vector fields $C^{-1,\beta}_\loc(M;\X)\subset C_c^{0,1}(M;\Lambda^1\otimes\Vol_M)'$ and the space of differential forms $C^{-1,\beta}_\loc(M;\Lambda^k)\subset C_c^{0,1}\big(M;(\bigwedge^k TM)\otimes\Vol_M\big)'$ for $1\le k\le n$. 

If $M$ is a oriented manifold, then $\Vol_M$ is canonically isomorphic to $\bigwedge^nT^*M$, the top degree alternating tensor of cotangent bundle. In this case a distribution is indeed an $n$-current. We leave the details to the reader.
\end{remark}

Next we define distributional sections on vector bundles.
\begin{lem}\label{Lem::DisInv::DisSecVB}
    Let $\alpha\in(0,1]\cup\{\mathrm{LogL}\}$, let $M$ be a $C^{1,1}$-manifold, and let $\E$ be a $C^{0,\alpha}$-vector bundle over $M$. Let $\beta,\gamma\in (0,1]\cup\{\mathrm{LogL}\}$ satisfy $\gamma\le\alpha$ and $\beta>1-\gamma$.
    \begin{enumerate}[parsep=-0.3ex,label=(\roman*)]
        \item\label{Item::DisInv::DisSecVB::FiniteGen} The space of sections $C^{0,\gamma}_\loc(M;\E)$ is a finitely generated module over the ring of function $C^{0,\gamma}_\loc(M)$. In particular if $s_1,\dots,s_N\in C^{0,\gamma}_\loc(M;\E)$ generate the module, then $s_1(p),\dots,s_N(p)$ span the fibre $\E_p$ for every $p\in M$.
        \item\label{Item::DisInv::DisSecVB::Ten} We have equality of tensor products,
        \begin{equation}\label{Eqn::DisInv::DisSecVB::Ten}
            C^{-1,\beta}_\loc(M)\otimes_{C^{0,\alpha}_\loc(M)}C^{0,\alpha}_\loc(M;\E)=C^{-1,\beta}_\loc(M)\otimes_{C^{0,\gamma}_\loc(M)}C^{0,\gamma}_\loc(M;\E).
        \end{equation}
    \end{enumerate}
\end{lem}
\begin{remark}\label{Rmk::DisInv::RmkDisSecVB}
\begin{enumerate}[parsep=-0.3ex,label=(\alph*)]
    \item Algebraically tensor products can only be equal in the sense of a unique isomorphism. Here we can make them equal set theoretically by  embedding both sides of \eqref{Eqn::DisInv::DisSecVB::Ten} to a larger ambient spaces. For example, we can view them as linear functionals on $C^{0,\gamma}_c(M;\E^*\otimes\Vol_M)$ (see Remark \ref{Rmk::DisInv::Functional}).
    \item\label{Item::DisInv::RmkDisSecVB::NoDepReg} We can use either side of \eqref{Eqn::DisInv::DisSecVB::Ten} to define the space of $C^{-1,\beta}$-sections of $\E$. We see that the definition does not depend on the regularity assumption $\E\in C^{0,\alpha}$: the space is the same if we view $\E$ as a $C^{0,\gamma}$-vector bundle. In fact the space $C^{-1,\beta}_\loc(M;\E)$ is defined as long as $\E\in C^{0,(1-\beta)^+}$.
    \item When $\E\le TM$ is a tangent subbundle, \eqref{Eqn::DisInv::DisSecVB::Ten} coincides with $C^{-1,\beta}_\loc(M;\E)$ as given in Definition \ref{Def::Intro::DisSec}. See Proposition \ref{Prop::DisInv::CharSec}.
\end{enumerate}
\end{remark}
\begin{proof}[Proof of Lemma \ref{Lem::DisInv::DisSecVB}]
\ref{Item::DisInv::DisSecVB::FiniteGen}: This is the non-smooth version of the Swan's Theorem. One can find the proof in \cite[Theorem 5.3]{Lewis}.

\medskip
\noindent\ref{Item::DisInv::DisSecVB::Ten}: Taking finitely many $C^{0,\alpha}$-sections $s_1,\dots,s_N$ that generates $C^{0,\alpha}_\loc(M;\E)$. Clearly a $C^{0,\gamma}$-section of $\E$ can be written as some $C^{0,\gamma}$-linear combinations of $s_1,\dots,s_N$, which means
\begin{equation*}
    C^{0,\gamma}_\loc(M;\E)=C^{0,\gamma}_\loc(M)\otimes_{C^{0,\alpha}_\loc(M)}C^{0,\alpha}_\loc(M;\E).
\end{equation*}

On the other hand, if $R$ is a commutative ring, $S$ is an $R$-commutative algebra, $P$ is a $S$-module and $Q$ is a $R$-module, then we have
\begin{equation*}
    P\otimes_RQ=P\otimes_S(S\otimes_RQ).
\end{equation*}
Taking $R=C^{0,\alpha}_\loc(M)$, $S=C^{0,\gamma}_\loc(M)$, $P=C^{-1,\beta}_\loc(M)$ and $Q=C^{0,\alpha}_\loc(M;\E)$ we get  \eqref{Eqn::DisInv::DisSecVB::Ten}.
\end{proof}

Now we can give some properties of distributional sections of a tangent subbundle. 

\begin{prop}[Characterization of distributional sections]\label{Prop::DisInv::CharSec}
Let $\alpha\in(0,1]\cup\{\mathrm{LogL}\}$, $\beta\in(1-\alpha,1)\cup\{\mathrm{LogL}\}$, let $M$ be an $n$-dimensional $C^{1,1}$-manifold and let $\V\le TM$ be  $C^{0,\alpha}$-subbundle of rank $r$.
Let $X\in C^{-1,\beta}_\loc(M;\X)$ and let $1-\beta<\gamma\le\alpha$, $1-\alpha<\delta\le\beta$. The following are equivalent:
\begin{enumerate}[parsep=-0.3ex,label=(S.\arabic*)]
    \item\label{Item::DisInv::CharSec::Dual} $\langle\theta,X\rangle=0$ as a generalized function on $M$ for all $C^{0,\alpha}$ 1-form $\theta$ on $\V^\bot$. This is the restatement of $X\in C^{-1,\beta}_\loc(M;\V)$ from Definition \ref{Def::Intro::DisInv}.
    \item\label{Item::DisInv::CharSec::Pullback} $X\in C^{-1,\delta}_\loc(M;\V)$, or equivalently $X\in C^{-1,\delta}_\loc(M;\V)\cap C^{-1,\beta}_\loc(M;\X)$.
    \item\label{Item::DisInv::CharSec::Ten1} $X\in C^{-1,\beta}_\loc(M)\otimes_{C^{0,\alpha}_\loc(M)}C^{0,\alpha}_\loc(M;\V)$. That is, there are finitely many $Y_1,\dots,Y_N\in C^{0,\alpha}_\loc(M;\V)$ and $f^1,\dots, f^N\in C^{-1,\beta}_\loc(M)$ such that $X=\sum_{i=1}^Nf^iY_i$ in the sense of distributions. 
    \item\label{Item::DisInv::CharSec::Ten2}  $X\in C^{-1,\beta}_\loc(M)\otimes_{C^{0,\gamma}_\loc(M)}C^{0,\gamma}_\loc(M;\V)$.
\end{enumerate}
\end{prop}
\begin{remark}
By the equivalence of \ref{Item::DisInv::CharSec::Ten1} and \ref{Item::DisInv::CharSec::Ten2} (since they are all equivalent to \ref{Item::DisInv::CharSec::Dual}), we see that definition of $C^{-1,\beta}$-sections does not depend on the regularity assumption that $\V$ is $C^{0,\alpha}$. Also see Remark \ref{Rmk::DisInv::RmkDisSecVB} \ref{Item::DisInv::RmkDisSecVB::NoDepReg}.
\end{remark}
\begin{proof}[Proof of Proposition \ref{Prop::DisInv::CharSec}]
Clearly \ref{Item::DisInv::CharSec::Dual} implies \ref{Item::DisInv::CharSec::Pullback}. The direction \ref{Item::DisInv::CharSec::Ten1} $\Rightarrow$ \ref{Item::DisInv::CharSec::Dual} follows from the fact that
\begin{equation*}
    \langle \theta,X\rangle=\sum_{i=1}^N\langle\theta,f^iY_i\rangle=\sum_{i=1}^Nf^i\langle\theta,Y_i\rangle=0,\quad\text{for every }C^{0,\alpha}\text{ 1-form }\theta\text{ on }\V^\bot.
\end{equation*}
The equivalence of \ref{Item::DisInv::CharSec::Ten1} and \ref{Item::DisInv::CharSec::Ten2} follows from Lemma \ref{Lem::DisInv::DisSecVB} \ref{Item::DisInv::DisSecVB::Ten} by considering $\E=\V$.

It remains to show  \ref{Item::DisInv::CharSec::Pullback} $\Rightarrow$ \ref{Item::DisInv::CharSec::Ten1}. 

\medskip
By Lemma \ref{Lem::DisInv::DisSecVB} \ref{Item::DisInv::DisSecVB::FiniteGen} we can find finitely many $C^{0,\alpha}$-vector fields $Y_1,\dots,Y_N$ that generates the module  $C^{0,\alpha}_\loc(M;\V)$ over the ring $C^{0,\alpha}_\loc(M)$. In particular $Y_1(p),\dots,Y_N(p)$ span $\V_p\le T_pM$ for all $p\in M$.

We can find an open cover $\{U_\kappa\}_\kappa$, associated with $r$-element subsets $I_\kappa\subseteq\{1,\dots,N\}$ and vector fields $Z^\kappa_1,\dots,Z^\kappa_{n-r}\in C^{0,\alpha}_\loc(U_\kappa)$ for each $\kappa$, such that $\{Y_i\}_{i\in I}\cup\{Z^\kappa_1,\dots,Z^\kappa_{n-r}\}$ form a local basis for $TM|_{U_\kappa}$. So for each $\kappa$, $(Y_i)_{i\in I_\kappa}$ are linearly independent and span $\V|_{U_\kappa}$.

Let $(\mu_\kappa^1,\dots,\mu_\kappa^r,\theta_\kappa^1,\dots,\theta_\kappa^{n-r})$ be  $C^{0,\alpha}$ 1-forms defined on $U_\kappa$ that form the dual basis to vector fields $(Y_{\sigma_\kappa(1)},\dots,Y_{\sigma_\kappa(r)},Z^\kappa_1,\dots,Z^\kappa_{n-r})$ where we take $\sigma_\kappa:\{1,\dots,r\}\to I_\kappa$ such that $\sigma_\kappa(1)<\dots<\sigma_\kappa(r)$. Therefore every $\tilde X\in C^{-1,\beta}_\loc(U_\kappa;\X)$ can be written as
\begin{equation*}
    \tilde X=\sum_{j=1}^r\langle\mu_\kappa^j,\tilde X\rangle Y_{\sigma_\kappa(j)}+\sum_{k=1}^{n-r}\langle\theta_\kappa^k,\tilde X\rangle Z^\kappa_k,\quad\text{on }U_\kappa\text{ for each }\kappa.
\end{equation*}

Let $\{\chi_\kappa\in C_c^{1,1}(U_\kappa)\}_\kappa$ be a $C^{1,1}$ partition of unity to the open cover $\{U_\kappa\}$, thus we have
\begin{equation*}
    X=\sum_{\kappa}\Big(\sum_{j=1}^r\langle\chi_\kappa\mu_\kappa^j, X\rangle\cdot Y_{\sigma_\kappa(j)}+\sum_{k=1}^{n-r}\langle\chi_\kappa\theta_\kappa^k,X\rangle \cdot Z^\kappa_k\Big).
\end{equation*}
Clearly $\chi_\kappa\theta_\kappa^1,\dots,\chi_\kappa\theta_\kappa^{n-r}\in C^{0,\alpha}(M;\V^\bot)$ are globally defined, since $\theta_\kappa^1,\dots,\theta_\kappa^{n-r}$ are $C^{0,\alpha}$-sections of $V^\bot|_{U_\kappa}$.

By assumption $X\in C^{-1,\beta}_\loc(M;\X)\cap C^{-1,\delta}_\loc(M;\V)$, we have $\langle\chi_\kappa\theta_\kappa^k,X\rangle=0 \in C^{-1,\delta}_\loc(M)$, so 
\begin{equation}\label{Eqn::DisInv::CharSecPf::Tmp}
    X=\sum_{\kappa}\sum_{j=1}^r\langle\chi_\kappa\mu_\kappa^j, X\rangle\cdot Y_{\sigma_\kappa(j)})=\sum_{i=1}^N\Big(\sum_\kappa\langle\chi_\kappa\mu_\kappa^{\sigma_\kappa^{-1}(j)}, X\rangle\Big)\cdot Y_j.
\end{equation}
Here we let $\mu^{\sigma_\kappa^{-1}(j)}_\kappa=0$ if $j\notin I_\kappa(=\sigma_\kappa\{1,\dots,r\})$.

Since $\langle\chi_\kappa\mu_\kappa^j, X\rangle\in C^{-1,\beta}_\loc(M)$ and the sums in \eqref{Eqn::DisInv::CharSecPf::Tmp} are locally finite, we know $X$ is the linear combinations of $Y_1,\dots,Y_N\in C^{0,\alpha}_\loc(M;\V)$ with coefficients in $C^{-1,\beta}_\loc(M)$. Therefore $X\in C^{-1,\beta}_\loc(M)\otimes_{C^{0,\alpha}_\loc(M)}C^{0,\alpha}_\loc(M;\V)$, which is the condition \ref{Item::DisInv::CharSec::Ten1}. 
\end{proof}

\subsection{Characterizations of involutivity}\label{Section::DisInv::CharInv}
In this part we provide some equivalent characterizations of distributional involutivity for a tangent subbundle. For completeness we restate the definition of distributional involutivity including the case $\frac12<\alpha<1$.
\begin{defn}\label{Def::DisInv::DefInv}
    Let $\alpha\in(\frac12,1)\cup\{\mathrm{LogL}\}$ and let $M$ be a $C^{1,1}$-manifold. We say a $C^{0,\alpha}$-subbundle $\V$ of $TM$ is distributional involutive, if for every vector fields $X,Y\in C^{0,\alpha}_\loc(M;\V)$ and 1-form $\theta\in C^{0,\alpha}_\loc(M;\V^\bot)$, the canonical pairing $\langle\theta,[X,Y]\rangle\in C^{-1,\alpha}_\loc(M)$ is identically zero in the sense of distributions.
\end{defn}
    

Recall that in the smooth setting, when $\V$ is a smooth tangent subbundle over a smooth manifold $M$, we have the following equivalent characterizations of involutivity:
\begin{enumerate}[parsep=-0.3ex,label=(\ref{Section::DisInv}.\Alph*)]
    \item\label{Item::DisInv::SmoothChar::LieAlg} $C^\infty_\loc(M;\V)$ is a Lie subalgebra of $C^\infty_\loc(M;\X)$, i.e. if $X$ and $Y$ are smooth sections of $\V$ then $[X,Y]$ is also a smooth section of $\V$.
    \item\label{Item::DisInv::SmoothChar::LieDer}  If $X\in C^\infty_\loc(M;\V)$ and $\theta\in C^\infty_\loc(M;\V^\bot )$, then the Lie derivative $\Lie X\theta\in C^\infty_\loc(M;\V^\bot)$.
    \item\label{Item::DisInv::SmoothChar::DiffIdeal}  On the exterior algebra $C^\infty_\loc(M;\Lambda^\bullet)=\bigoplus_{k=0}^{\dim M}C^\infty_\loc(M;\Lambda^k)$, the (two-sided) ideal $C^\infty_\loc(M;\V^\bot)\wedge C^\infty_\loc(M;\Lambda^\bullet)$ is closed under differential.
    \item\label{Item::DisInv::SmoothChar::2Form} On 2-forms we have $d(C^\infty_\loc(M;\V^\bot))\subseteq C^\infty_\loc(M;\V^\bot)\wedge C^\infty_\loc(M;\Lambda^1)$.
\end{enumerate}

These characterizations can all be generalized to non-smooth subbundles.

We first describe some of the characterizations involving Lie brackets. 

\begin{prop}[Characterizations of distributional involutivity I]\label{Prop::DisInv::CharInv1}
Let $\alpha\in(\frac12,1)\cup\{\mathrm{LogL}\}$, let $M$ be a $n$-dimensional $C^{1,1}$-manifold and let $\V$ be a rank $r$ $C^{0,\alpha}$-subbundle of $TM$. 

Let $\frac12<\beta<\alpha$. The following conditions for $\V$ are equivalent:
\begin{enumerate}[parsep=-0.3ex,label=(I.\arabic*)]
    \item\label{Item::DisInv::CharInv1::Sec1} For every sections $X,Y\in C^{0,\alpha}_\loc(M;\V)$, we have $[X,Y]\in C^{-1,\alpha}_\loc(M;\V)$.
    \item\label{Item::DisInv::CharInv1::Sec2} For every sections $X,Y\in C^{0,\beta}_\loc(M;\V)$, we have $[X,Y]\in C^{-1,\beta}_\loc(M;\V)$.
    \item\label{Item::DisInv::CharInv1::Pair1} For every $X,Y\in C^{0,\alpha}_\loc(M;\V)$ and $\theta\in C^{0,\alpha}_\loc(M;\V^\bot)$, we have $\langle\theta,[X,Y]\rangle=0\in C^{-1,\alpha}_\loc(M)$. In other words $\V$ is distributional involutive (by Definition \ref{Def::DisInv::DefInv}).
    \item\label{Item::DisInv::CharInv1::Pair2} For every $X,Y\in C^{0,\alpha}_\loc(M;\V)$ and $\theta\in C^{0,\alpha}_\loc(M;\V^\bot)$, we have $d\theta(X,Y)=0\in C^{-1,\alpha}_\loc(M)$.
    \item\label{Item::DisInv::CharInv1::Pair3} For every $X\in C^{0,\alpha}_\loc(M;\V) $ and $\theta\in C^{0,\alpha}_\loc(M;\V^\bot)$, we have $\Lie X\theta\in C^{-1,\alpha}_\loc(M;\V^\bot)$.
    \item\label{Item::DisInv::CharInv1::Gen1} For every $p\in M$ there is a neighborhood $U\subseteq M$ of $p$ and a $C^{0,\alpha}$-local basis $(X_1,\dots,X_r)$ for $\V|_U$ such that $X_1,\dots,X_r$ are pairwise commutative in the sense of distributions.
    \item\label{Item::DisInv::CharInv1::Gen2} If $X_1,\dots,X_{N}\in C^{0,\alpha}_\loc(M;\V)$ are global generators for $\V$ (that is $X_1,\dots,X_N$ span $\V$ at every point in $M$), then there are $c_{ij}^k\in C^{-1,\alpha}_\loc(M)$ for $1\le i,j,k\le N$ such that 
    \begin{equation}\label{Eqn::DisInv::CharInv1::Gen2cijk}
        [X_i,X_j]=\sum_{k=1}^Nc_{ij}^kX_k\quad\text{on }M,\quad \forall\ 1\le i,j\le N.
    \end{equation}
\end{enumerate}
\end{prop}
\begin{remark}
\begin{enumerate}[parsep=-0.3ex,label=(\alph*)]
    \item In the statements, \ref{Item::DisInv::CharInv1::Sec1} and \ref{Item::DisInv::CharInv1::Pair3} are the direct generalizations to \ref{Item::DisInv::SmoothChar::LieAlg} and \ref{Item::DisInv::SmoothChar::LieDer} respectively.
    \item By the equivalence of \ref{Item::DisInv::CharInv1::Sec1} and \ref{Item::DisInv::CharInv1::Sec2}, we see that the definition of distributional involutivity does not depend on the regularity assumption $\V\in C^{0,\alpha}$. In fact it is enough to assume $\V\in C^{0,{\frac12}^+}$ so that the Lie bracket is well-defined as a distributional vector field. Also see Section \ref{Section::DisInv::Border}.
    
    \smallskip
    In this way the characterizations are still equivalent if we replace every $\alpha$ by $\beta$ in \ref{Item::DisInv::CharInv1::Pair1} - \ref{Item::DisInv::CharInv1::Gen2}.

    \item The above characterizations  do not use the information of $\dim M$ and $\rank \V$. Thus it is possible to extend the definition of involutive condition to the case where $M$ is not necessary a manifold or $\V$ is not a subbundle.
    
    \smallskip
    For example, let $\alpha\in(\frac12,1)\cup\{\mathrm{LogL}\}$ and $\beta\in(1-\alpha,1)$. Let $\mathscr M\subseteq C^{0,\alpha}_\loc(M;\X)$ be a submodule over the ring $C^{0,\alpha}_\loc(M;\X)$. We can say $\mathscr M$ is $C^{-1,\beta}$-involutive, if $[\mathscr M,\mathscr M]\subseteq \mathscr M\otimes_{C^{0,\alpha}_\loc(M)}C^{-1,\beta}_\loc(M)$. See \cite[Definition 3.0.2, Examples 3.4.9 and 3.4.11]{LidingThesis} for more discussions.
    \item In particular from \ref{Item::DisInv::CharInv1::Gen2} we can focus on the specific choice of $X_1,\dots,X_N$ rather than on the ambient space $\V$. In this setting we only consider the module spanned by $X_1,\dots,X_N$, where $\Span(X_1(p),\dots,X_N(p))\le T_pM$ have constant rank among all $p\in M$. 
    
\end{enumerate}
\end{remark}

\begin{proof}[Proof of Proposition \ref{Prop::DisInv::CharInv1}]We consider the directions \ref{Item::DisInv::CharInv1::Sec2} $\Rightarrow$ \ref{Item::DisInv::CharInv1::Sec1} $\Leftrightarrow$ \ref{Item::DisInv::CharInv1::Pair1} $\Leftrightarrow$ \ref{Item::DisInv::CharInv1::Pair2} $\Leftrightarrow$ \ref{Item::DisInv::CharInv1::Pair3}; \ref{Item::DisInv::CharInv1::Pair1} $\Rightarrow$ \ref{Item::DisInv::CharInv1::Gen1} $\Rightarrow$ \ref{Item::DisInv::CharInv1::Sec2}, and \ref{Item::DisInv::CharInv1::Sec1} $\Leftrightarrow$ \ref{Item::DisInv::CharInv1::Gen2}.

When $X,Y\in C^{0,\alpha}_\loc(M;\V)$ and $\theta\in C^{0,\alpha}_\loc(M;\V^\bot)$, by Lemma \ref{Lem::Hold::MultLoc} \ref{Item::Hold::MultLoc::WellDef} and Convention \ref{Conv::DisInv::-1log} we have $[X,Y]\in C^{-1,\alpha}_\loc(M;\X)$, $\langle\theta,[X,Y]\rangle,d\theta(X,Y)\in C^{-1,\alpha}_\loc(M)$ and $\Lie X\theta\in C^{-1,\alpha}_\loc(M;\Lambda^1)$.
 
Clearly \ref{Item::DisInv::CharInv1::Sec2} $\Rightarrow$ \ref{Item::DisInv::CharInv1::Sec1}. 

The direction \ref{Item::DisInv::CharInv1::Sec1} $\Leftrightarrow$ \ref{Item::DisInv::CharInv1::Pair1} follows from Definition \ref{Def::Intro::DisSec}, where $[X,Y]$ is a section of $\V$ if and only if $\langle\theta,[X,Y]\rangle=0$ for all $\theta\in C^{0,\alpha}_\loc(M;\V^\bot)$.

\medskip
\noindent\ref{Item::DisInv::CharInv1::Pair1} $\Leftrightarrow$ \ref{Item::DisInv::CharInv1::Pair2} $\Leftrightarrow$ \ref{Item::DisInv::CharInv1::Pair3}: For $X,Y\in C^{0,\alpha}_\loc(M;\V)$ and $\theta\in C^{0,\alpha}_\loc(M;\V^\bot)$, in the sense of distributions we have $d\theta(X,Y)=X\langle \theta, Y\rangle-Y\langle\theta,X\rangle-\langle\theta,[X,Y]\rangle=-\langle\theta,[X,Y]\rangle$ and $\langle\Lie X\theta,Y\rangle=X\langle\theta,Y\rangle-\langle\theta,[X,Y]\rangle$. Thus
\begin{equation}\label{Eqn::DisInv::CharInv1::Tmp}
    \langle\theta,[X,Y]\rangle=-d\theta(X,Y)=\langle\Lie X\theta,Y\rangle.
\end{equation}
If one term in \eqref{Eqn::DisInv::CharInv1::Tmp} is zero then the other two vanish as well. Thus we get \ref{Item::DisInv::CharInv1::Pair1} $\Leftrightarrow$ \ref{Item::DisInv::CharInv1::Pair2}.
Note that \ref{Item::DisInv::CharInv1::Pair3} is equivalent as saying $\langle\Lie X\theta,Y\rangle=0$ for all $X,Y\in C^{0,\alpha}_\loc(M;\V)$ and $\theta\in C^{0,\alpha}_\loc(M;\V^\bot)$, thus \ref{Item::DisInv::CharInv1::Pair3} is equivalent to  \ref{Item::DisInv::CharInv1::Pair1} and \ref{Item::DisInv::CharInv1::Pair2}.
    
    \medskip
\noindent \ref{Item::DisInv::CharInv1::Pair1} $\Rightarrow$ \ref{Item::DisInv::CharInv1::Gen1}: This is done by Lemma \ref{Lem::PfThm::CanGen}, where the result Lemma \ref{Lem::PfThm::CanGen} \ref{Item::PfThm::CanGen::Inv} still holds if we assume $\V$ is $C^{0,\alpha}$.

\medskip
\noindent \ref{Item::DisInv::CharInv1::Gen1} $\Rightarrow$ \ref{Item::DisInv::CharInv1::Sec2}: Let $Y_1,Y_2\in C^{0,\beta}_\loc(M;\V)$, by locality of distributions it suffices to show that for any $p\in M$ there is a neighborhood $U\subseteq M$ of $p$ such that $[Y_1,Y_2]\in C^{-1,\beta}_\loc(U;\V|_U)$. Fix $p\in M$, and let $U\subseteq M$ and  $X_1,\dots,X_r\in C^{0,\alpha}_\loc(U;\V|_U)$ be as in \ref{Item::DisInv::CharInv1::Gen1}. By Proposition \ref{Prop::DisInv::CharSec} \ref{Item::DisInv::CharSec::Dual} $\Rightarrow$ \ref{Item::DisInv::CharSec::Ten1} we can write $Y_i=\sum_{j=1}^rf_i^jX_j$ for some $f_i^j\in C^{-1,\beta}_\loc(U)$, $i=1,2$ and $1\le j\le r$. Therefore
\begin{equation*}
    [Y_1,Y_2]=\sum_{i,j=1}^r[f_1^iX_i,f_2^jX_j]=\sum_{i,j=1}^rf_1^i(X_if_2^j)X_j-f_2^j(X_jf_1^i)X_i+f_1^if_2^j[X_i,X_j]=\sum_{i,j=1}^r\big(f_1^j(X_jf_2^i)-f_2^j(X_jf_1^i)\big)X_i.
\end{equation*}

We see that $[Y_1,Y_2]$ is the $C^{-1,\beta}_\loc$ linear combinations of $X_1,\dots,X_r$, thus $[Y_1,Y_2]\in C^{-1,\beta}_\loc(U)\otimes_{C^{0,\alpha}_\loc(U)}C^{0,\alpha}(U;\V|_U)$. By Proposition \ref{Prop::DisInv::CharSec} \ref{Item::DisInv::CharSec::Ten1}$\Rightarrow$\ref{Item::DisInv::CharSec::Dual} $[Y_1,Y_2]\in C^{-1,\beta}_\loc(U;\V|_U)$. Since $p\in M$ is arbitrary, we get $[Y_1,Y_2]\in C^{-1,\beta}_\loc(M;\V)$, proving \ref{Item::DisInv::CharInv1::Sec2}.

\medskip The proof of \ref{Item::DisInv::CharInv1::Gen2} $\Rightarrow$ \ref{Item::DisInv::CharInv1::Sec1} is identical to the proof of \ref{Item::DisInv::CharInv1::Gen1} $\Rightarrow$ \ref{Item::DisInv::CharInv1::Sec2} except we do not need to localize the objects on an open subset. We omit the details to the readers.

\medskip
\noindent \ref{Item::DisInv::CharInv1::Sec1} $\Rightarrow$ \ref{Item::DisInv::CharInv1::Gen2}: Note that by Lemma \ref{Lem::DisInv::DisSecVB} \ref{Item::DisInv::DisSecVB::FiniteGen} we can find finitely many $X_1,\dots,X_N\in C^{0,\alpha}_\loc(M;\V)$ that span $\V$ at every point.

By assumption \ref{Item::DisInv::CharInv1::Sec1}, $[X_i,X_j]\in C^{-1,\alpha}(M;\V)$ for all $1\le i,j\le N$.
So by Proposition \ref{Prop::DisInv::CharSec} \ref{Item::DisInv::CharSec::Dual}$\Rightarrow$\ref{Item::DisInv::CharSec::Ten1}, $[X_i,X_j]\in C^{-1,\alpha}_\loc(M)\otimes_{C^{0,\alpha}_\loc(M)}C^{0,\alpha}(M;\V)$.
Since $X_1,\dots,X_N$ generates $C^{0,\alpha}(M;\V)$ as a module over $C^{0,\alpha}_\loc(M)$, we can write $[X_i,X_j]=\sum_{k=1}^Nc_{ij}^kX_k$ for some $c_{ij}^k\in C^{-1,\alpha}_\loc(M)$. This gives the condition \ref{Item::DisInv::CharInv1::Gen2}.
\end{proof}

In Proposition \ref{Prop::DisInv::CharInv1} we mostly focus on the side of vector fields. We can consider the characterizations using differential forms which are related to the sections of $\V^\bot$.

\begin{prop}[Characterizations of distributional involutivity II]\label{Prop::DisInv::CharInv2}
Let $\alpha\in(\frac12,1)\cup\{\mathrm{LogL}\}$, let $M$ be a $n$-dimensional $C^{1,1}$-manifold and let $\V$ be a rank $r$ $C^{0,\alpha}$-subbundle of $TM$. 

Then $\V$ is distributional involutive if and only if  either of the following condition holds:
\begin{enumerate}[parsep=-0.3ex,label=(I.\arabic*)]\setcounter{enumi}{7}
    
        \item\label{Item::DisInv::CharInv2::Ann} Let $1\le k\le n$ and $\omega\in C^{0,\alpha}_\loc(M;\Lambda^k)$ satisfy $\omega(Y_1,\dots,Y_k)=0$ for every $Y_1,\dots,Y_k\in C^{0,\alpha}_\loc(M;\V)$. Then we have $d\omega(Y_1,\dots,Y_{k+1})=0$ as distribution for every $Y_1,\dots,Y_{k+1}\in C^{0,\alpha}_\loc(M;\V)$.   
    
    \item\label{Item::DisInv::CharInv2::Diff} $d(C^{0,\alpha}_\loc(M;\V^\bot))\subseteq (C^{0,\alpha}_\loc(M;\V^\bot)\wedge C^{0,\alpha}_\loc(M;\Lambda^1))\otimes_{C^{0,\alpha}_\loc(M)}C^{-1,\alpha}_\loc(M)$. That is, for every $\lambda\in C^{0,\alpha}_\loc(M;\V^\bot)$, $d\lambda\in C^{-1,\alpha}_\loc(M;\Lambda^2)$ can be written as a finite sum
    \begin{equation*}
        d\lambda=\sum_{i=1}^Nf_i\cdot\theta^i\wedge\mu^i,\quad\text{for some }N\ge1,\quad f_i\in C^{-1,\alpha}_\loc(M),\quad \theta^i\in C^{0,\alpha}_\loc(M;\V^\bot),\quad \mu^i\in C^{0,\alpha}_\loc(M;\Lambda^1).
    \end{equation*}

    \item\label{Item::DisInv::CharInv2::coSec} For every $\lambda^1,\dots,\lambda^{n-r+1}\in C^{0,\alpha}_\loc(M;\V^\bot)$ we have 
    \begin{equation}\label{Eqn::DisInv::CharInv2::coSecEqn}
        \lambda^1\wedge\dots\wedge\lambda^{n-r}\wedge d\lambda^{n-r+1}=0\in C^{-1,\alpha}_\loc(M;\Lambda^{n-r+2}).
    \end{equation}
\end{enumerate}
\end{prop}

Here \ref{Item::DisInv::CharInv2::Diff} is the direct generalization to \ref{Item::DisInv::SmoothChar::2Form}.

We shall see that \ref{Item::DisInv::CharInv2::Ann} generalizes \ref{Item::DisInv::SmoothChar::DiffIdeal}, because the space of all differential forms that annihilate sections of $\V$ is indeed an ideal of the exterior algebra $C^{0,\alpha}_\loc(M;\Lambda^\bullet)$.

To prove the equivalence of \ref{Item::DisInv::CharInv2::Diff} and \ref{Item::DisInv::CharInv2::Ann} we need the following:

\begin{lem}\label{Lem::DisInv::IdealChar}
    Let $\alpha$ and $\V\le TM$ satisfies the assumptions of Proposition \ref{Prop::DisInv::CharInv2}. Then for $1\le k\le n$,
    \begin{equation}\label{Eqn::DisInv::IdealChar::Eqn}
        C^{0,\alpha}_\loc(M;\V^\bot)\wedge C^{0,\alpha}_\loc(M;\Lambda^{k-1})=\{\omega\in C^{0,\alpha}_\loc(M;\Lambda^k):\omega(X_1,\dots,X_k)=0\ \forall X_1,\dots,X_k\in C^{0,\alpha}_\loc(M;\V)\}.
    \end{equation}
\end{lem}
\begin{proof}

Let $W$ be a $n$-dimensional vector space and $V\le W$ is a rank $r$ subspace. By linear algebra, as the subspaces of alternating tensors we have
\begin{equation}\label{Eqn::DisInv::IdealChar::PointwiseEqn}
    V^\bot\wedge \textstyle\bigwedge^{k-1}W^*=\big\{\phi\in\textstyle\bigwedge^{k}W^*:\phi(v_1,\dots,v_k)=0\text{ for all }v_1,\dots,v_k\in V\big\}.
\end{equation}
For each $p\in M$, by taking $W=T_pM$ and $V=\V_p$ in \eqref{Eqn::DisInv::IdealChar::PointwiseEqn}, we get \eqref{Eqn::DisInv::IdealChar::Eqn} and complete the proof.
\end{proof}

\begin{proof}[Proof of Proposition \ref{Prop::DisInv::CharInv2}]

By Lemma \ref{Lem::DisInv::IdealChar} we see that 
\ref{Item::DisInv::CharInv2::Diff} is the special case to \ref{Item::DisInv::CharInv2::Ann} when $k=1$. So \ref{Item::DisInv::CharInv2::Ann} $\Rightarrow$ \ref{Item::DisInv::CharInv2::Diff}.


\medskip
We are going to prove \ref{Item::DisInv::CharInv1::Gen2} $\Rightarrow$  \ref{Item::DisInv::CharInv2::Ann} and \ref{Item::DisInv::CharInv2::Diff} $\Rightarrow$ \ref{Item::DisInv::CharInv2::coSec} $\Rightarrow$ \ref{Item::DisInv::CharInv1::Gen1}. This would complete the proof since by Proposition \ref{Prop::DisInv::CharInv1}, \ref{Item::DisInv::CharInv1::Gen2} and \ref{Item::DisInv::CharInv1::Gen1} are both equivalent characterizations of involutivity.

\medskip
\noindent\ref{Item::DisInv::CharInv1::Gen2} $\Rightarrow$ \ref{Item::DisInv::CharInv2::Ann}: Let $X_1,\dots,X_N$ satisfy the condition \ref{Item::DisInv::CharInv1::Gen2}, and let $\omega$ be a $k$-form annihilating all sections of $\V$. In particular $\omega(X_{i_1},\dots,X_{i_k})=0$ for all $\{i_1,\dots,i_k\}\subseteq\{1,\dots,N\}$. By passing to $C^{0,\alpha}$-linear combinations it suffices to show that $\omega(Y_1,\dots,Y_{k+1})=0$ for $\{Y_1,\dots,Y_{k+1}\}\subseteq\{X_1,\dots,X_N\}$. 

Indeed by direct computation along with \eqref{Eqn::DisInv::CharInv1::Gen2cijk} we have
\begin{align*}
    d\omega(X_{i_0},\dots,X_{i_{k}})=&\sum_{0\le u<v\le k}(-1)^{u+v}\omega([X_{i_u},X_{i_v}],X_{i_0},\dots , X_{i_{u-1}},X_{i_{u+1}}\dots ,X_{i_{v-1}},X_{i_{v+1}},\dots ,X_{i_k})
    \\
    =&\sum_{0\le u<v\le k}(-1)^{u+v}\sum_{l=1}^Nc_{i_ui_v}^l\cdot\omega(X_l,X_{i_0},\dots , X_{i_{u-1}},X_{i_{u+1}}\dots ,X_{i_{v-1}},X_{i_{v+1}},\dots ,X_{i_k})=0.
\end{align*}
This gives \ref{Item::DisInv::CharInv2::Ann}.

\medskip
\noindent\ref{Item::DisInv::CharInv2::Diff} $\Rightarrow$ \ref{Item::DisInv::CharInv2::coSec}:
Let $\lambda^1,\dots,\lambda^{n-r+1}\in C^{0,\alpha}_\loc(M;\V^\bot)$. By assumption \ref{Item::DisInv::CharInv2::Diff} we can write $d\lambda^{n-r+1}=\sum_{i=1}^Nf_i\theta^i\wedge\mu^i$ for some $\theta^i\in C^{0,\alpha}_\loc(M;\V^\bot)$. So
\begin{equation*}
    \lambda^1\wedge\dots\wedge\lambda^{n-r}\wedge d\lambda^{n-r+1}=\sum_{i=1}^N \lambda^1\wedge\dots\wedge\lambda^{n-r}\wedge \theta^i\wedge (f_i\mu^i).
\end{equation*}

On the other hand $\V^\bot $ is a cotangent subbundle with rank $(n-r)$, we see that $\lambda^1,\dots,\lambda^{n-r},\theta^i$ are linear dependent at every point in $M$, thus $\lambda^1\wedge\dots\wedge\lambda^{n-r}\wedge \theta^i\equiv0$. We conclude that $\lambda^1\wedge\dots\wedge\lambda^{n-r}\wedge d\lambda^{n-r+1}=0$ as a distributional differential form.

\medskip
\noindent\ref{Item::DisInv::CharInv2::coSec} $\Rightarrow$ \ref{Item::DisInv::CharInv1::Gen1}: Fixing a point $p\in M$, by Lemma \ref{Lem::PfThm::CanGen} we can find a $C^{1,1}$-chart $(x^1,\dots,x^r,y^1,\dots,y^{n-r}):U\subseteq M\to\R^n$ near $p$ and $C^{0,\alpha}$-vector fields $X_1,\dots,X_r$ on $U$ that span $\V|_U$ and have the form $X_j=\Coorvec{x^j}+\sum_{k=1}^{n-r}b_j^k\Coorvec{y^k}$. What we need is to show that $[X_j,X_k]=0$ for all $1\le j<k\le r$.

Clearly $(X_1,\dots,X_r,\Coorvec{y^1},\dots,\Coorvec{y^{n-r}})$ span the tangent space at every point in $U$. The dual basis has the form $(dx^1,\dots,dx^r,\theta^1,\dots,\theta^{n-r})$ where $\theta^l=dy^l-\sum_{j=1}^rb_j^ldx^j$, $l=1,\dots,n-r$ (see \eqref{Eqn::PfThm::CanGen::DualForm}). Thus for $1\le l\le n-r$
\begin{align*}
    d\theta^l&=\sum_{i,j=1}^r\frac{\partial b_i^l}{\partial x^j}dx^i\wedge dx^j+\sum_{i=1}^r\sum_{k=1}^{n-r}\frac{\partial b_i^l}{\partial y^k}dx^i\wedge dy^k
    \\
    &=\sum_{1\le i<j\le r}\Big(\frac{\partial b_i^l}{\partial x^j}-\frac{\partial b_j^l}{\partial x^i}+\sum_{k=1}^{n-r}\big(b_j^k\frac{\partial b_i^l}{\partial y^k}-b_i^k\frac{\partial b_i^k}{\partial y^l}\big)\Big)dx^i\wedge dx^j+\sum_{i=1}^r\sum_{k=1}^{n-r}\frac{\partial b_j^l}{\partial y^k}dx^i\wedge\theta^k
    \\
    &\equiv\sum_{1\le i<j\le r}\langle\theta^l,[X_i,X_j]\rangle dx^i\wedge dx^j\pmod{\theta^1,\dots,\theta^r}.
\end{align*}
Thus $$\theta^1\wedge\dots\wedge\theta^{n-r}\wedge d\theta^l=\sum_{1\le i<j\le r}\langle\theta^l,[X_i,X_j]\rangle \theta^1\wedge\dots\wedge\theta^{n-r}\wedge dx^i\wedge dx^j,\quad 1\le l\le n-r.$$
By shrinking $U$ if necessary we can find global sections $\lambda^1,\dots,\lambda^r\in C^{0,\alpha}_\loc(M;\V^\bot)$ such that $\lambda^l|_U=\theta^l$. By assumption \ref{Item::DisInv::CharInv2::coSec} we see that $\langle\theta^l,[X_i,X_j]\rangle=0\in C^{-1,\alpha}_\loc(U)$. Since $[X_i,X_j]=\sum_{k=1}^{n-r}\langle\theta^l,[X_i,X_j]\rangle\Coorvec{y^l}$ we know $[X_i,X_j]=0$ are commutative, proving the condition \ref{Item::DisInv::CharInv1::Gen1}.
\end{proof}

\subsection{The borderline of definiability}\label{Section::DisInv::Border}

In Section \ref{Section::DisInv::CharInv} we give equivalent characterizations of involutivity on a subbundle $\V\le TM$ which is $C^{0,\frac12+\eps}$ for some $\eps>0$. We can push the definability to its limit by merely assuming $\V$ is $L^\infty\cap H^\frac12$ and the manifold is 1 derivative better than $L^\infty\cap H^\frac12$.

Here $H^\frac12$ is the fractional $L^2$-Sobolev space: for an open set $\Omega\subseteq\R^n$, we use
\begin{equation}\label{Eqn::Borderline::H12Norm}
    \|f\|_{H^\frac12(\Omega)}^2:=\|f\|_{L^2(\Omega)}^2+\int_{\Omega\times\Omega}\frac{|f(x)-f(y)|^2}{|x-y|^{n+1}}dxdy.
\end{equation}
Similarly we define $H^{\frac32}(\Omega)=\{f\in H^\frac12(\Omega):\nabla f\in H^\frac12(\Omega;\R^n)\}$ and $H^{-\frac12}(\Omega):=\{g_0+\sum_{i=1}^n\partial_ig_i:g_0,\dots,g_n\in H^\frac12(\Omega)\}$. For $s=-\frac12,\frac12,\frac32$ we define $H^s_\loc(\Omega)$ in the usual way.
One can check that we have equality with the standard Besov spaces $H^s(\Omega)=B_{22}^s(\Omega)$ when $\Omega$ is $\R^n$ or a bounded smooth domain. See \cite[Proposition 1.59 and Definition 2.68]{BahouriCheminDanchin} and \cite[Section 3.4.2]{TriebelTheoryOfFunctionSpacesI}.

\begin{lem}\label{Lem::Borderline::H12lem}
    Let $U,V\subseteq\R^n$ be two open sets and let $G\in C^{0,1}_\loc\cap H^\frac32_\loc(U;V)$ be a bijection such that $G^{-1}\in C^{0,1}_\loc\cap H^\frac32_\loc(V;U)$. We have the following property of multiplications and compositions:
    \begin{enumerate}[parsep=-0.3ex,label=(P.\arabic*)]
        \item\label{Item::Borderline::Alg} $L^\infty_\loc\cap H^\frac12_\loc(U)$ is an algebra: if $f,g\in L^\infty_\loc\cap H^\frac12_\loc(U)$, so is $fg$.
        \item\label{Item::Borderline::Comp} If $f\in L^\infty_\loc\cap H^\frac12_\loc(V)$ and $F\in C^{0,1}_\loc\cap H^\frac32_\loc(V)$, then $f\circ G\in L^\infty_\loc\cap H^\frac12_\loc(U)$ and $F\circ G\in C^{0,1}_\loc\cap H^\frac32_\loc(U)$.
        \item\label{Item::Borderline::Mult} If $f\in L^\infty_\loc\cap H^\frac12_\loc(U)$ and $g\in  H^{-\frac12}_\loc+L^1_\loc(U)$ then $fg\in H^{-\frac12}_\loc(U)+L^1_\loc(U)$.
        \item\label{Item::Borderline::DistComp} If $f\in H^{-\frac12}_\loc+L^1_\loc(V)$  then $f\circ G\in H^{-\frac12}_\loc+L^1_\loc(U)$.
    \end{enumerate}
\end{lem}

\begin{proof}[Sketch]\ref{Item::Borderline::Alg} follows from the fact that $L^\infty\cap H^\frac12(\R^n)$ is an algebra. See \cite[Corollary 2.86]{BahouriCheminDanchin}.

For \ref{Item::Borderline::Comp}, clearly $f\circ G\in L^\infty_\loc$. By the norm characterization \eqref{Eqn::Borderline::H12Norm} we know $H^\frac12$ is preserved under bi-Lipschitz maps, so $f\circ G\in H^\frac12_\loc(U)$.

The result $F\circ G\in C^{0,1}_\loc\cap H^\frac32_\loc$ follows from taking $f$ as the components of $\nabla F$ and the chain rule  $\nabla (F\circ G)=((\nabla F)\circ G)\cdot\nabla G$. 

For \ref{Item::Borderline::Mult}, we use \eqref{Eqn::Hold::Bony} to write the three sums decomposition $fg=\Pi(f,g)+\Pi(g,f)+R(f,g)$. By \cite[Theorem 2.82]{BahouriCheminDanchin} $\Pi:L^\infty(\R^n)\times H^{-\frac12}(\R^n)\to H^{-\frac12}(\R^n)$ and $\Pi:H^{-\frac12}(\R^n)\times L^\infty(\R^n)\to H^{-\frac12}(\R^n)$ are both bounded. One can check that $R:H^\frac12(\R^n)\times H^{-\frac12}(\R^n)\to L^1(\R^n)$ and clearly $[(f,g)\mapsto fg ]:L^\infty(\R^n)\times L^1(\R^n)\to L^1(\R^n)$. Thus by passing to local we have the desired result.

For \ref{Item::Borderline::DistComp} clearly $[f\mapsto f\circ G]:L^1_\loc(V)\to L^1_\loc(U)$, it suffices $f\in H^{-\frac12}_\loc(U)$. We can write $f=g_0+\sum_{i=1}^n\partial_ig_i$ for $g_0,\dots,g_n\in H^\frac12_\loc(U)$. By \eqref{Eqn::Borderline::H12Norm} and that $G$ is locally bi-Lipschitz we know $g_i\circ G\in H^\frac12_\loc(U)$. Write $\Psi=(\Psi^1,\dots,\Psi^n):=G^{-1}$, we have
\begin{equation*}
    (\partial_ig_i)\circ G=\sum_{j=1}^n(\partial_j(g_i\circ G))\cdot((\partial_i\Psi^j)\circ G).
\end{equation*}
By \ref{Item::Borderline::Mult} we know $(\partial_ig_i)\circ G\in H^{-\frac12}_\loc+L^1_\loc$, concluding the result.
\end{proof}

We can now  define $C^{0,1}\cap H^\frac32$ manifolds in a natural way: the $C^{0,1}\cap H^\frac32$ atlas $\As$ on an $n$-dimensional topological manifold $M$ is the maximal collection of homeomorphisms $\psi:U_\psi\subseteq M\to\R^n $ such that the transition maps $\phi\circ\psi^{-1}:\psi(U_\phi\cap U_\psi)\to\R^n$ are $C^{0,1}_\loc\cap H^\frac32_\loc$ whenever the domain is nonempty. By Lemma \ref{Lem::Borderline::H12lem} \ref{Item::Borderline::Alg}, a $C^{0,1}\cap H^\frac32$ coordinate cover uniquely determines a maximal $C^{0,1}\cap H^\frac32$ atlas.

Similar to Definition \ref{Def::DisInv::DefFunVF} \ref{Item::DisInv::DefFunVF::VF}, let $\Xs=\{L^\infty\cap H^\frac12,H^{-\frac12}+L^1\}$, a $\Xs$-vector field on $M$ can be defined as a collection $X=\{X_\psi\in \Xs_\loc(U_\psi;\R^n)\}_{\psi\in \As}$ satisfying \eqref{Eqn::DisInv::DefFunVF::TransVF}. Using \ref{Item::Borderline::Alg} and \ref{Item::Borderline::Comp} for $\Xs=L^\infty\cap H^\frac12$, and using \ref{Item::Borderline::Mult} and \ref{Item::Borderline::DistComp} for $\Xs=H^{-\frac12}+L^1$, we see that to check $X\in\Xs$ it suffices to check $X_{\psi_i}\in\Xs_\loc(U_{\psi_i};\R^n)$ where $\{\psi_i\}$ is a $C^{0,1}\cap H^\frac32$ coordinates cover.

If $M$ is in addition $C^1$ then we can define a $L^\infty\cap H^\frac12$ tangent subbundle in a natural way, since $TM$ is a topological manifold. However we cannot define a ``continuous'' tangent bundle on Lipschitz manifolds. In this case we can define a rank-$r$ $L^\infty\cap H^\frac12$-subbundle  $\V$ as a collection $\{\V_\psi\in H^\frac12_\loc(U_\psi;\mathbf{Gr}(r,\R^n))\}_{\psi\in\As}$ of Sobolev mappings, such that
\begin{align*}
    \label{Eqn::Borderline::TransBdd}&(\nabla(\phi\circ\psi^{-1})(x))\bullet_{\mathrm{GL}}\V_\psi(x)=\V_\phi\circ\phi\circ\psi^{-1}(x)\text{ almost every }x\in U_\phi\cap U_\psi,\text{ provided }U_\phi\cap U_\psi\neq\varnothing.
\end{align*}
Here $\mathbf{Gr}(r,\R^n)$ is the Grassmannian space endowed with a standard group action
\begin{equation*}
    \operatorname{GL}(\R^n)\curvearrowright\mathbf{Gr}(r,\R^n),\quad A\bullet_{\mathrm{GL}}\Span(v_1,\dots,v_r):=\Span(Av_1,\dots,Av_r).
\end{equation*}
Note that  the compactness of $\mathbf{Gr}(r,\R^n)$ ensures every maps to it is bounded, thus there is no need to say $\V_\psi\in L^\infty_\loc(U_\psi;\mathbf{Gr}(r,\R^n))$.

We can define the dual bundle $\V^\bot$ as the collection $\{\V^\bot_\psi \}_{\psi\in\As}$, where $\V^\bot_\psi$ is given by the composition map $U_\psi\xrightarrow{\V_\psi}\mathbf{Gr}(r,\R^n)\xrightarrow{\cong}\mathbf{Gr}(n-r,\R^n)$. Clearly $\V^\bot_\psi\in H^\frac12_\loc(U_\psi;\mathbf{Gr}(n-r,\R^n))$ for each $\psi$.

We can define sections and involutivity as the following: still we assume $M\in C^{0,1}\cap H^\frac32$ and $\V\in L^\infty\cap H^\frac12$: 
\begin{itemize}[parsep=-0.3ex]
    \item We say a $L^\infty$-vector field $X=\{X_\psi\}_\psi$ is a section of $\V=\{\V_\psi\}$ if $X_\psi(x)\in \V_\psi(x)(\subseteq\R^n)$ for almost every $x\in U_\psi$ and for every $\psi\in\As$. The section of $\V^\bot$ is defined similarly.
    \item We say a $(H^{-\frac12}+L^1)$ vector field $X=\{X_\psi\}_\psi$ is a distributional section of $\V$, if $\langle \theta,X_\psi\rangle=0\in (H^{-\frac12}+L^1)(U_\psi)$ holds for all $\psi\in\As$ and 1-form $\theta\in(L^\infty\cap H^\frac12)(U_\psi;\V_\psi^\bot)$.
    \item We say $\V$ is distributional involutive, if for every $(L^\infty\cap H^\frac12)$-sections $X,Y$, the $(H^{-\frac12}+L^1)$ vector field $[X,Y]$ is a distributional section of $\V^\bot$.
\end{itemize}

We can ask whether Propositions \ref{Prop::DisInv::CharSec}, \ref{Prop::DisInv::CharInv1} and \ref{Prop::DisInv::CharInv2} are still true if we replace $(C^{0,\alpha},C^{-1,\alpha})$ with $(L^\infty\cap H^\frac12,H^{-\frac12}+L^1)$. For example
\begin{ques}\label{Ques::DisInv::CharLim}
\begin{enumerate}[parsep=-0.3ex,label=(\roman*)]
    \item Let $X\in (H^{-\frac12}_\loc+L^1_\loc)(M;\X)$ be a distributional section of $\V$. Do we have $X\in (L^\infty_\loc\cap H^\frac12_\loc)(M;\V)\otimes_{L^\infty_\loc\cap H^\frac12_\loc(M)}(H^{-\frac12}_\loc+L^1_\loc)(M)$?
    \item Assume that $\V\in L^\infty\cap H^\frac12$ has rank $r$. Suppose for every $\lambda^1,\dots,\lambda^{n-r+1}\in (L^\infty_\loc\cap H^\frac12_\loc)(M;\V^\bot)$ we have $\lambda^1\wedge\dots\wedge\lambda^{n-r}\wedge d\lambda^{n-r+1}=0$. Is $\V$ involutive in the sense of distributions?
\end{enumerate}
\end{ques}


It is not known to the author whether the above questions are true for general $L^\infty\cap H^\frac12$ tangent subbundles. 
\section{Further Remarks}

\subsection{A PDE counterpart}
We can interpret Theorem \ref{Thm::MainThm1} in terms of the first-order PDE system following from \cite[Section 5]{Rampazzo}:

\begin{thm}\label{Thm::Further::PDEThm}Let $M$ be a $n$-dimensional $C^{1,1}$-manifold. Let $L_1,\dots,L_r$ be some first order differential operators on $M$ with log-Lipschitz coefficients, such that
\begin{equation}\label{Eqn::Further::PDEAssumption1}
    L_jL_k-L_kL_j=\sum_{l=1}^rc_{jk}^l\cdot L_l,
\end{equation}
for some $c_{jk}^l\in C^{-1,1^-}_\loc(M)$, $1\le j,k,l\le r$, in the sense of distributions.

Suppose there is a $(n-r)$-dimensional $C^1$-submanifold  $S\subset \R^n$ such that 
\begin{equation}\label{Eqn::Further::PDEAssumption2}
    \Span (L_1|_q\dots,L_r|_q)\oplus T_qS=T_q M,\qquad\forall q\in S.
\end{equation}

Let $0<\beta\le1$ and $h\in C^{0,\beta}(S)$. Then for any $p\in S$ and any $0<\eps<\beta$, there is a neighborhood $U\subseteq M$ of $p$, such that there exists a unique solution $f\in C^{0,\beta-\eps}(U)$ to the Cauchy problem
\begin{equation}\label{Eqn::Further::CountPDE}
    \begin{cases}L_jf=0,&1\le j\le r,\\f|_{S\cap U}=h.\end{cases}
\end{equation}
\end{thm}Here $L_jf=0$ holds in the sense of distributions.

Applying Lemma \ref{Lem::Hold::MultLoc} \ref{Item::Hold::MultLoc::WellDef} with $L_j\in C^{0,1^-}_\loc(U;\X)$ and $\nabla f\in C^{-1,\beta-\eps}(U;\R^n)$, we know that $L_jf\in C^{-1,\beta-\eps}_\loc(U)$ is defined as a distribution.

If in addition $L_1,\dots,L_r$ have little log-Lipschitz coefficients, we can show that $f$ given \eqref{Eqn::Further::CountPDE} is indeed $C^{0,\beta^-}$. We omit the proof to reader.

The proof is based on the construction of $\Phi$ in \eqref{Eqn::PfThm::DefofPhi}, along with the following proposition on inverse functions. 
\begin{prop}\label{Prop::Further::Mu}
Let $(x,y)=(x^1,\dots,x^r,y^1,\dots,y^{n-r})$ be the standard coordinate system of $\R^r\times\R^{n-r}$. Let $U\subseteq\R^n$ be an open neighborhood of $0$. Let $X_1,\dots,X_r$ be some commutative log-Lipschitz vector fields on $U$ that have the form $X_j=\Coorvec{x^j}+\sum_{k=1}^{n-r}b_j^k\Coorvec{y^k}$ for $1\le j\le r$.

Let $\Gamma=(\Gamma^1,\dots,\Gamma^{n-r}):\R^{n-r}\to\R^r$ be a $C^1$-map such that $\Gamma(0)=0$ and $\nabla \Gamma(0)=0$. We define a map $\mu$ as
\begin{equation}\label{Eqn::Further::Mu}
    \mu(s):= y(e^{-\Gamma(s)\cdot X}(\Gamma(s),s))=y(e^{-\Gamma^1(s)X_1}\dots e^{-\Gamma^r(s)X_r}(\Gamma(s),s)),
\end{equation}
provided the ODE flows are defined.

Then for any $\eps>0$ there is an open neighborhood $U''=U''_\eps\subseteq \R^{n-r}$ of $0\in\R^{n-r}$ such that $\mu:U''\to\R^{n-r}$ is a $C^{0,1-\eps}$-map which is homeomorphic to its image, and $\mu^{-1}:\mu(U'')\to\R^{n-r}$ is also $C^{0,1-\eps}$.
\end{prop}
\begin{proof}The local definedness of $\mu$ at $s=0$ is clear, since $\Gamma(0)=0\in\R^r$ and the ODE flows are all defined near $0\in\R^n$. Note that  by construction we have $$(0,\mu(s))=e^{-\Gamma(s)\cdot X}(\Gamma(s),s),\quad\text{for }s\text{ in the domain}.$$

For a given $\eps>0$, the $C^{0,1-\eps}$-regularity of $\mu$ near $0$ follows from Corollary \ref{Cor::ODE::MultFlow} \ref{Item::ODE::MultFlow::LLReg} and the fact that $(0,\mu(s))=\Psi^X(-\Gamma(s),(\Gamma(s),s))$ where $\Psi^X$ is the map in \eqref{Eqn::ODE::MultFlow::Psi}. It remains to show that there are a $c=c_\eps>0$ and a neighborhood $U''=U''_\eps\subseteq\R^{n-r}$ of $0$, such that
\begin{equation}\label{Eqn::Further::MuInj}
    |e^{-\Gamma(s_1)\cdot X}(\Gamma(s_1),s_1)-e^{-\Gamma(s_2)\cdot X}(\Gamma(s_2),s_2)|\ge c|s_1-s_2|^{\frac1{1-\eps}},\quad\forall s_1,s_2\in U''.
\end{equation}

Once \eqref{Eqn::Further::MuInj} is done, we know $s\mapsto(0,\mu(s))$ is injective near $0$, so is $\mu$. Then by invariance of domain $\mu:U''\to\mu(U'')\subseteq\R^{n-r}$ is homeomorphism. Therefore $\mu^{-1}:\mu(U'')\to\R^{n-r}$ is a well-defined $C^{0,1-\eps}$-map. By shrinking $U''$ such that $\mu\in C^{0,1-\eps}$, we then finish the proof.

\medskip
To prove \eqref{Eqn::Further::MuInj}, for any open neighborhood $V_1\subseteq \R^{n-r}$ of $0$ such that  $\{e^{u\cdot X}(\Gamma(s),s):s\in V_1,|u|\le|\Gamma(s)|\}\Subset U$, we have
\begin{equation*}
    |e^{(\Gamma(s_2)-\Gamma(s_1))\cdot X}(\Gamma(s_1),s_1)-(\Gamma(s_2),s_2)|\ge |s_1-s_2|-\|\nabla \Gamma\|_{C^0(V_1;\R^n)}\|X\|_{C^0(U;\R^n)}|s_2-s_1|,\quad s_1,s_2\in V_1.
\end{equation*}
Since $\nabla\Gamma(0)=0$, by continuity we can choose $V_1$ small enough so that $\|\nabla \Gamma\|_{C^0(V_1)}\|X\|_{C^0(U)}\le\frac12$. Thus
\begin{equation}\label{Eqn::Further::PfMuInj::Tmp1}
    |e^{(\Gamma(s_2)-\Gamma(s_1))\cdot X}(\Gamma(s_1),s_1)-(\Gamma(s_2),s_2)|\ge\tfrac12|s_1-s_2|,\quad\forall s_1,s_2\in V_1.
\end{equation}

Let $W_1\Subset U$ be an open set containing $ \{(\Gamma(s),s):s\in V_1\}$. By Corollary \ref{Cor::ODE::MultFlow} \ref{Item::ODE::MultFlow::LLReg} with the chosen $\eps$ in the assumption, we can find $\delta>0$ such that $e^{u\cdot X}\in C^{0,1-\eps}(W_1;U)$ uniformly for $|u|<\delta$. Take $U'':=\{s\in V_1,|\Gamma(s)|<\frac12\delta\}$, we have for every $s_1,s_2\in U''$,
\begin{equation}\label{Eqn::Further::PfMuInj::Tmp2}
    |e^{(\Gamma(s_2)-\Gamma(s_1))\cdot X}(\Gamma(s_1),s_1)-(\Gamma(s_2),s_2)|\le \sup_{|u|<\delta}\|e^{u\cdot X}\|_{C^{0,1-\eps}(W_1;\R^n)} |e^{-\Gamma(s_1)\cdot X}(\Gamma(s_1),s_1)-e^{-\Gamma(s_2)}(\Gamma(s_2),s_2)|^{1-\eps}.
\end{equation}
Combining \eqref{Eqn::Further::PfMuInj::Tmp1} and \eqref{Eqn::Further::PfMuInj::Tmp2} we get \eqref{Eqn::Further::MuInj} with $c=(2\sup\limits_{|u|<\delta}\|e^{u\cdot X}\|_{C^{0,1-\eps}(W_1)})^{-\frac1{1-\eps}}>0$, finishing the proof.
\end{proof}

\begin{proof}[Proof of Theorem \ref{Thm::Further::PDEThm}]
By \eqref{Eqn::Further::PDEAssumption2} and that $\dim S=n-r$, we know that $L_1,\dots,L_r$ are linearly independent at every point in $S$. By continuity there is a neighborhood $\tilde U\subseteq M$ of $S$ such that they are still  linearly independent.

We can define $\V\le TM|_{\tilde U}$ as a tangent subbundle spanned by $L_1,\dots,L_r$. By construction $\V$ is a rank $r$ log-Lipschitz subbundle. By Proposition \ref{Prop::DisInv::CharInv1} \ref{Item::DisInv::CharInv1::Gen2}$\Rightarrow$\ref{Item::DisInv::CharInv1::Pair1}, $\V$ is involutive in the sense of distributions.

By passing to an invertible linear transform we can assume that $\Coorvec{y^1}|_p,\dots,\Coorvec{y^{n-r}}|_p$ span $T_pS\le T_p M$. Applying Lemma \ref{Lem::PfThm::CanGen} and shrinking $\tilde U$ if necessary, we can find a coordinate system $(x^1,\dots,x^r,y^1,\dots,y^{n-r})$ and a log-Lipschitz local basis $(X_1,\dots,X_r)$ for $\V$ on $\tilde U$ such that $[X_j,X_k]=0$ for $1\le j,k\le r$ and have expression $X_j=\Coorvec{x^j}+\sum_{k=1}^{n-r}b_j^k\Coorvec{y^k}$. Thus by shrinking $\tilde U$, in the coordinate $(x^1,\dots,x^r,y^1,\dots,y^{n-r})$ we can write
\begin{equation}\label{Eqn::Further::PfPDEThm::Tmp0}
    S\cap\tilde U=\{(\Gamma(s),s)\},\quad s\in U_0'',
\end{equation}
where $U_0''\subseteq\R^{n-r}$ is a small neighborhood of $0$ and $\Gamma:U_0''\to\R^r$ is a $C^1$-map.

By taking a translation we can assume $p=0\in\R^n$, thus $\Gamma$ satisfies $\Gamma(0)=0$ and $\nabla\Gamma(0)=0$.

Take $\tilde\eps>0$ be such that $\beta(1-\tilde\eps)^2\ge\beta-\eps$. Recall in the proof of Theorem \ref{Thm::MainThm1}, the map $\Phi(u,v)=e^{u\cdot X}(0,v)$ is defined in a neighborhood of $(u,v)=0$, such that
\begin{itemize}[parsep=-0.3ex]
    \item  $\Phi$ is a locally homeomorphism and $\Phi^{-1}$ is of the form $\Phi^{-1}=(x,\lambda)$, where $\lambda$ is given in \eqref{Eqn::PfThm::lambda}.
    \item There is a neighborhood $U_1\subseteq M$ of $0$, such that $\lambda:U_1\to\R^{n-r}$ is $C^{0,1-\tilde\eps}$.
    \item $\frac{\partial\Phi}{\partial u^j}=X_j\circ\Phi$ hold for $j=1,\dots,r$.
\end{itemize}

If $U\subseteq\tilde U$ and $f:U\subseteq M\to\R$ is a H\"older function satisfying $X_jf=0$ for $j=1,\dots,r$, then $\Coorvec{u^j}(f\circ\Phi)(u,v)=(X_jf)\circ\Phi(u,v)=0$. Therefore
\begin{equation}\label{Eqn::Further::PfPDEThm::Tmp1}
    f(\Phi(u,v))=f(\Phi(0,v))=f(0,v),\quad\text{for }(u,v)\text{ in the domain}.
\end{equation}

By Proposition \ref{Prop::Further::Mu} the map $\mu(s)$ given in \eqref{Eqn::Further::Mu} is a locally homeomorphism near $s=0$ and satisfies
\begin{equation}\label{Eqn::Further::PfPDEThm::Tmp2}
    \Phi(\Gamma(s),\mu(s))=e^{\Gamma(s)\cdot X}(0,\mu(s))=(\Gamma(s),s).
\end{equation}

Thus if $f|_{U\cap S}=h$, then in a neighborhood of $0$ we must have
\begin{align*}
    f(x,y)&=f(\Phi(x,\lambda(x,y)))=f(\Phi(x,\mu(\mu^{-1}\circ\lambda(x,y))))&(\Phi\circ(x,\lambda)=\id)
    \\
    &=f(\Phi(0,\mu((\mu^{-1}\circ\lambda)(x,y))))=f(\Phi(\Gamma((\mu^{-1}\circ\lambda)(x,y))),\mu((\mu^{-1}\circ\lambda)(x,y)))&(\text{by }\eqref{Eqn::Further::PfPDEThm::Tmp1})
    \\
    &=f(\Phi(\Gamma,\mu)\circ(\mu^{-1}\circ\lambda)(x,y))=f((\Gamma,\id_{\R^{n-r}})\circ(\mu^{-1}\circ\lambda)(x,y))&(\text{by }\eqref{Eqn::Further::PfPDEThm::Tmp2})
    \\
    &=h(\Gamma(\mu^{-1}\circ\lambda(x,y)),\mu^{-1}\circ\lambda(x,y)).&(\text{by }\eqref{Eqn::Further::PfPDEThm::Tmp0})
\end{align*}
Thus $f$ is completely determined by $h$. This completes the proof of the uniqueness.

\medskip
Let $U''_2\subseteq U''_0$ be a neighborhood of $0$ such that $\mu^{-1}:U''_2\to\R^{n-r}$ is $C^{0,1-\tilde\eps}$ and let $U:=\mu^{-1}(U''_2)\cap U_1$. By composition we have $\mu^{-1}\circ\lambda\in C^{0,(1-\tilde\eps)^2}(U;\R^{n-r})$. Since $h\in C^{0,\beta}(S)$, taking compositions we get $f\in C^{0,\beta(1-\tilde\eps)^2}(U)\subseteq C^{0,\beta-\eps}(U)$, finishing the proof of the existence and hence the whole proof.
\end{proof}

\subsection{Sharpness of the $C^{0,1-\eps}$-smoothness}
In Theorem \ref{Thm::MainThm1}, for a fixed $0<\eps<1$ we can find a parameterization $\Phi$ such that $\Phi$ and $\Phi^{-1}$ are both $C^{0,1-\eps}$.
The statement, that $\Phi$ is a ``$ C^{0,1-\eps}$-diffeomorphism'', is sharp, in the sense that there is a log-Lipschitz involutive subbundle $\V$ defined near a fixed point $p$, such that there is no topological parameterization $\Phi$ near $p$ that makes both $\Phi$ and $\Phi^{-1}$ to be $ C^{0,1^-}$.

More precisely, it is the following:
\begin{prop}\label{Prop::Further::SharpProp}Let $X(x,y)=(1,y\log |y|)$ be a vector field in $\R^2_{x,y}$. Let $\V\le T\R^2$ be the rank $1$ tangent subbundle spanned by $X$.

Suppose that $\Phi(u,v):\Omega\subseteq\R^2_{u,v}\to\R^2_{x,y}$ is a topological parameterization near $0\in\R^2_{x,y}$ such that $\frac{\partial\Phi}{\partial u}$ is continuous, and $\frac{\partial\Phi}{\partial u}(u,v)$ span $\V_{\Phi(u,v)}$ for every $(u,v)\in\Omega$.

Then either $\Phi\notin  C^{0,1^-}(\Omega;\R^2)$ or $\Phi^{-1}\notin C^{0,1^-}(\Phi(\Omega);\R^2)$.
\end{prop}
Here $X$  has Zygmund regularity, which is a bit more regular than the log-Lipschitz regularity (see \cite[Proposition 2.107]{BahouriCheminDanchin}). Clearly $\V$ is involutive because it has rank 1.

\begin{proof}Without loss of generality $\Phi(0,0)=(0,0)$. We assume by contrast that both $\Phi$ and $\Phi^{-1}$ are $C^{0,1^-}$ near the origin.

We write $\Phi(u,v)=:(\phi(u,v),\psi(u,v))$ and we define $G(u,v):=(\phi(u,v),v)$ for  $(u,v)\in\Omega$. By assumption $\Phi\in C^{0,1^-}$ so $G\in C^{0,1^-}_\loc(\Omega;\R^2)$. 

First we show that $G$ is bijective and has a $C^{0,1^-}$-inverse $G^{-1}:G(\Omega)\subset\R^2_{u,v}\to\R^2_{u,v}$, such that $\frac{\partial G^{-1}}{\partial u}\in C^0$.

\medskip
By assumption $\frac{\partial\phi}{\partial u}\in C^0(\Omega)$ is non-vanishing, so for each $v$, $\phi(\cdot,v)$ has a $C^1$-inverse. 

Define $h(u,v):=\phi(\cdot,v)^{-1}(u)$, so $G^{-1}(u,v)=(h(u,v),v)$ is well-defined on $G(\Omega)$.

By shrinking $\Omega$, we can assume that $|\phi_u|$ is bounded below. Denoting $C_0:=(\inf_\Omega|\phi_u|)^{-1}<\infty$, we have $\sup_{\Phi(\Omega)}|h_u|= C_0$. Therefore if $u_1,u_2,v_1,v_2\in\R$ satisfy $(u_1,v_1),(u_2,v_1),(u_2,v_2)\in\Phi(\Omega)$, then 
\begin{align*}
    |h(u_1,v_1)-h(u_2,v_2)|\le&|h(u_1,v_1)-h(u_2,v_1)|+|h(u_2,v_1)-h(u_2,v_2)|
    \\\le& C_0|u_1-u_2|+C_0|\phi(h(u_2,v_1),v_1)-\phi(h(u_2,v_2),v_1)|&(\sup|h_u|=(\inf|\phi_u|)^{-1}=C_0)
    \\
    =&C_0|u_1-u_2|+C_0|u_2-\phi(h(u_2,v_2),v_1)|&(\phi(h(\cdot,v_1),v_1)=\id)
    \\
    =&C_0|u_1-u_2|+C_0|\phi(h(u_2,v_2),v_2)-\phi(h(u_2,v_2),v_1)|&(\phi(h(\cdot,v_2),v_2)=\id)
    \\\le& C_0|u_1-u_2|+C_0\|\phi\|_{C^{1-\eps}}|v_1-v_2|^{1-\eps},&\forall \eps>0.\quad (\text{since }\phi\in C^{0,1^-})
\end{align*}
Therefore $h\in C^{0,1^-}(G(\Omega))$, which means $G^{-1}\in C^{0,1^-}(G(\Omega);\R^2)$. In particular, we know $G(\Omega)$ is an open set.

Note that $h_u(u,v)=\frac1{\phi_u(G^{-1}(u,v))}$, so $\frac{\partial G^{-1}}{\partial u}=(h_u,0)=(\frac1{\phi_u}\circ G^{-1},0)\in C^{0,1^-}\circ C^{0,1^-}\subset C^0(G(\Omega);\R^2)$. This proves the claim that $G^{-1},\frac{\partial G^{-1}}{\partial u}\in C^{0,1^-}$.

\medskip
Defining $\tilde \Phi:=\Phi\circ G^{-1}$ on $G(\Omega)$. We have  $$\tilde\Phi(u,v)=(u,\psi(G^{-1}(u,v)))=(u,\psi(h(u,v),v)),\quad\tilde\Phi_u(u,v)=(1,\psi_u(h(u,v),v)\cdot h_u(u,v)).$$
Therefore $\tilde\Phi\in C^{0,1^-}$, $\tilde \Phi^{-1}=G\circ\Phi^{-1}\in C^{0,1^-}$ and $\tilde\Phi_u\in C^0$.

We can now get a contradiction using the assumption that $\tilde\Phi,\tilde\Phi^{-1}$ are both $ C^{0,1^-}$ near $(0,0)$. 

\medskip
Define $\tilde\psi:=\psi\circ G^{-1}$, we have $\tilde\Phi(u,v)=(u,\tilde\psi(u,v))$.
By assumption $\V_{\tilde\Phi(u,v)}=\Span X(\tilde\Phi(u,v))= \Span\frac{\partial\tilde\Phi}{\partial u}(u,v)$ for $(u,v)\in G(\Omega)$, so $(1,\tilde\psi_u)=(1,\tilde\psi\log|\tilde\psi|)$. Solving the equation we have: $$\tilde\psi(u,v)=|\tilde\psi(0,v)|^{e^u}\sgn\tilde\psi(0,v),\quad\text{ for }u,v\text{ closed to }0.$$

By construction $\tilde\Phi(0,0)=0$. So for all $\eps>0$, there is a $C_\eps\approx\|\tilde\Phi^{-1}\|_{C^{0,1-\eps}}$ such that when $(0,v)\in G(\Omega)$, $$|v|=|v-0|=|\tilde\Phi^{-1}(\tilde\Phi(0,v))-\tilde\Phi^{-1}(\tilde\Phi(0,0))|\le C_\eps|\tilde\Phi(0,v)-\tilde\Phi(0,0)|^{1-\eps}=C_\eps|\tilde\psi(0,v)|^{1-\eps}.$$

Now for any $\delta>0$, take $u=-\frac\delta2$ and take $\eps>0$ such that $1-\eps>e^{-\delta/2}$. So $$\textstyle|\tilde\psi(-\frac\delta2,v)|=|\tilde\psi(0,v)|^{e^{-\delta/2}}\ge C_\eps^{-1}|v|^{\frac{e^{-\delta/2}}{1-\eps}},\quad\text{ whenever }(-\frac\delta2,v)\in G(\Omega).$$

So $\tilde\psi\notin C^{0,\frac{e^{-\delta/2}}{1-\eps}}$ in the domain, for such $\delta$ and $\eps$. Therefore $\tilde\Phi\notin C^{0,1^-}$ near $(0,0)$, contradiction!
\end{proof}
\begin{remark}\label{Rmk::Further::PreciseCountGenThm}
Using the similar argument in the proof of Proposition \ref{Prop::Further::SharpProp}, we can find a non-vanishing little log-Lipschitz vector field $X$ on $\R^2$ such that if $\Phi(u,v)$ is a corresponding topological parameterization to the subbundle $\V:=\Span X$ near $(0,0)$, then either $\Phi\notin \Clog$ or $\Phi^{-1}\notin \Clog$.

Indeed, we can consider $X(x,y)=(1,\mu(|y|))$ where
\begin{equation*}
    \mu(r):=\begin{cases}r\log^2\frac1r,&0\le r\le\frac1e,\\\frac1e&r\ge\frac1e.\end{cases}
\end{equation*}
By solving the ODE we have
\begin{equation*}
    e^{tX}(0,r)=(t,r^{(1+t\log\frac1r)^{-1}}),\quad\text{when } 0<r<\textstyle\frac1e\text{ and }0<t<1-\frac1{\log\frac1r}.
\end{equation*}

One can see that $r\mapsto |r|^{(1+t\log\frac1{|r|})^{-1}}$ is not log-Lipschitz near 0, no matter how small $t>0$ is. This implies that $\exp_X$ is not log-Lipschitz near $(0,(0,0))\in\R\times\R^2$.

We omit the details to reader.
\end{remark}

\subsection{Problems for the general Osgood cases}
Recall that a vector field $X$ is said to be Osgood if $|X(x)-X(y)|\le\mu(|x-y|)$ holds for some  modulus of continuity $\mu$ satisfying $\int_0^1\frac{dr}{\mu(r)}=+\infty$. By the Osgood's uniqueness theorem \cite{Osgood} the ODE flow $e^{tX}$ is locally well-defined.

We can naturally to ask whether the Frobenius theorem holds in the general Osgood setting:
\begin{prob}\label{Prob::OsgFro}
Let $\V^r\le TM^n$ be an Osgood subbundle which is involutive in the sense of distributions. For any $p\in M$, can we find a topological parameterization $\Phi(u,v):\Omega\subseteq\R^r\times\R^{n-r}\to M$ that satisfies the following?
\begin{enumerate}[nolistsep,label=(\roman*)]
    \item $\Phi(0,0)=p$.
    \item $\frac{\partial\Phi}{\partial u^1},\dots,\frac{\partial\Phi}{\partial u^r}\in C^0(\Omega;TM)$.
    \item For every $(u,v)\in\Omega$, $\frac{\partial\Phi}{\partial u^1}(u,v),\dots,\frac{\partial\Phi}{\partial u^r}(u,v)$ spans $\V_{\Phi(u,v)}\le T_{\Phi(u,v)}M$.
\end{enumerate}
\end{prob}

We know that if such $\Phi$ does exist, then by viewing $\Omega=\Omega'_u\times\Omega''_v$, and applying the argument similar to Proposition \ref{Prop::Further::SharpProp}, the regularity estimate $\Phi\in C^1_u(\Omega';C^0_v(\Omega'';M))$ is the best possible we can get. This can be done using Osgood's lemma \cite[Lemma 3.4]{BahouriCheminDanchin}.

Problem \ref{Prob::OsgFro} is reduced to the flow commuting problem, Problem \ref{Prob::Intro::FlowComProb}, when vector fields are merely Osgood. Unfortunately we do not know how to prove the Problem \ref{Prob::Intro::FlowComProb} for non log-Lipschitz vector fields in general. Precisely speaking, it is the following,
\begin{prob}\label{Prob::ConjOsgoodFlowCom}
Let $X,Y$ be two Osgood vector fields in $\R^n$ satisfying $[X,Y]=0$ in the sense of distributions. Do we must have $e^{tX}\circ e^{sY}=e^{sY}\circ e^{tX}$ locally for small $t$ and $s$?
\end{prob}

One attempt to deal with the Problem \ref{Prob::ConjOsgoodFlowCom} is to show whether an Osgood subbundle is asymptotically involutive (see  \cite[Definition 1.17]{ContFro}). Based on this, we can ask an even weaker question:
\begin{prob}
Is there a modulus of continuity $\mu$ such that the following hold?
\begin{itemize}[parsep=-0.3ex]
    \item $\lim\limits_{r\to0}\mu(r)^{-1}\cdot r\log\frac1r=0$.
    \item If $\V$ is a tangent subbundle with modulus of continuity $\mu$, and if $\V$ is involutive in the sense of distributions, then $\V$  is asymptotically involutive.
\end{itemize}
\end{prob}

A major difficulty in understanding Problem \ref{Prob::ConjOsgoodFlowCom} is to verify \eqref{Eqn::Intro::AsyInv1}. But the exponential term in \eqref{Eqn::Intro::AsyInv2} is somehow subtle.

\bibliographystyle{amsalpha}
\small
\bibliography{reference}

\center{\textit{The Ohio State University, Department of Mathematics, 231 W 18th Ave, Columbus, OH 43210}}

\center{\textit{yao.1015@osu.edu}}

\center{MSC 2020: 58A30 (Primary), 53C12 and 60L40 (Secondary)}
\end{document}